\newif\iffinal
\finalfalse	
\finaltrue	

\documentclass[letterpaper,12pt,reqno]{amsart}
\RequirePackage[utf8]{inputenc}
\usepackage[portrait,margin=2.8cm]{geometry}
\iffinal\else\usepackage[notref,notcite]{showkeys}\fi
\usepackage{mathrsfs,soul,mathtools}
\usepackage{hyperref}
\usepackage[foot]{amsaddr}
\usepackage{amssymb,amsthm,amsfonts,amsbsy,latexsym,dsfont}
\usepackage[textsize=small]{todonotes}       
\usepackage{graphicx}
\usepackage[numeric,initials,nobysame]{amsrefs}
\usepackage{upref,setspace}
\usepackage{tikz}
\usetikzlibrary{calc,intersections,through,backgrounds,shapes.geometric}
\usetikzlibrary{graphs}

\mathtoolsset{showonlyrefs=true}

\usepackage{caption}
\usepackage{subcaption}
\usepackage{fourier}
\usepackage{color}

\usepackage{enumerate}
\newenvironment{enumeratei}{\begin{enumerate}[\upshape (i)]}{\end{enumerate}}
\newenvironment{enumeratea}{\begin{enumerate}[\upshape (a)]}{\end{enumerate}}

\newenvironment{enumerateA}{\begin{enumerate}[\upshape (A)]}{\end{enumerate}}

\usepackage{paralist}

\newenvironment{inparaenuma}{\begin{inparaenum}[\upshape \bfseries (a) ]}{\end{inparaenum}}

\newenvironment{inparaenumi}{\begin{inparaenum}[\upshape (i) ]}{\end{inparaenum}}

\numberwithin{equation}{section}
\numberwithin{figure}{section}
\numberwithin{table}{section}

\newtheorem{thm}{Theorem}[section]
\newtheorem{lem}[thm]{Lemma}

\newtheorem{prop}[thm]{Proposition}
\newtheorem{defn}[thm]{Definition}

\newtheorem{ass}[thm]{Assumption}
\newtheorem*{ass*}{Assumption}
\newtheorem*{theorem*}{Theorem}

\newtheorem{conj}[thm]{Conjecture}
\newtheorem{lemma}[thm]{Lemma}
\newtheorem{problem}[thm]{Problem}

\newtheorem*{thm3.9*}{Theorem 3.9*}

\theoremstyle{definition}
\newtheorem{rem}{Remark}

\renewcommand{\leq}{\leqslant}
\renewcommand{\geq}{\geqslant}

\newcommand{\ind}{\1}
\newcommand{\eps}{\varepsilon}

\newcommand{\set}[1]{\left\{#1\right\}}

\newcommand{\probc}{\stackrel{\mathrm{P}}{\longrightarrow}}

\newcommand{\convas}{\stackrel{\mathrm{a.s.}}{\longrightarrow}}

\def\qed{ \hfill $\blacksquare$}
\newcommand{\Var}{{\mathrm{Var}}}

\newcommand*\xbar[1]{%
  \hbox{%
    \vbox{%
      \hrule height 0.5pt 
      \kern0.5ex
      \hbox{%
        \kern-0.5em
        \ensuremath{#1}%
        \kern-0.1em
      }%
    }%
  }%
}

\newcommand{\cA}{\mathcal{A}}\newcommand{\cB}{\mathcal{B}}\newcommand{\cC}{\mathcal{C}}
\newcommand{\cD}{\mathcal{D}}\newcommand{\cE}{\mathcal{E}}\newcommand{\cF}{\mathcal{F}}
\newcommand{\cG}{\mathcal{G}}\newcommand{\cH}{\mathcal{H}}\newcommand{\cI}{\mathcal{I}}
\newcommand{\cJ}{\mathcal{J}}\newcommand{\cL}{\mathcal{L}}
\newcommand{\cM}{\mathcal{M}}\newcommand{\cN}{\mathcal{N}}\newcommand{\cO}{\mathcal{O}}
\newcommand{\cP}{\mathcal{P}}\newcommand{\cR}{\mathcal{R}}
\newcommand{\cS}{\mathcal{S}}\newcommand{\cT}{\mathcal{T}}\newcommand{\cU}{\mathcal{U}}
\newcommand{\cV}{\mathcal{V}}
\newcommand{\cZ}{\mathcal{Z}}

\newcommand{\vA}{\mathbf{A}}
\newcommand{\vD}{\mathbf{D}}

\newcommand{\vM}{\mathbf{M}}

\newcommand{\vS}{\mathbf{S}}\newcommand{\vT}{\mathbf{T}}
\newcommand{\vV}{\mathbf{V}}
\newcommand{\vX}{\mathbf{X}}\newcommand{\vY}{\mathbf{Y}}\newcommand{\vZ}{\mathbf{Z}}
\newcommand{\vc}{\mathbf{c}}

\newcommand{\vl}{\mathbf{l}}

\newcommand{\vp}{\mathbf{p}}\newcommand{\vq}{\mathbf{q}}
\newcommand{\vt}{\mathbf{t}}\newcommand{\vu}{\mathbf{u}}
\newcommand{\vx}{\mathbf{x}}
\newcommand{\vy}{\mathbf{y}}

\newcommand{\mvU}{\boldsymbol{U}}
\newcommand{\mvX}{\boldsymbol{X}}

\newcommand{\mvv}{\boldsymbol{v}}
\newcommand{\mvw}{\boldsymbol{w}}\newcommand{\mvx}{\boldsymbol{x}}

\newcommand{\mvtheta}{\boldsymbol{\theta}}

\newcommand{\mvnu}{\boldsymbol{\nu}}

\newcommand{\mvPi}{\boldsymbol{\Pi}}

\newcommand{\mvUpsilon}{\boldsymbol{\Upsilon}}

\newcommand{\fN}{\mathfrak{N}}
\newcommand{\fP}{\mathfrak{P}}

\newcommand{\fm}{\mathfrak{m}}

\newcommand{\bE}{\mathbb{E}}
\newcommand{\bG}{\mathbb{G}}
\newcommand{\bL}{\mathbb{L}}
\newcommand{\bN}{\mathbb{N}}
\newcommand{\bQ}{\mathbb{Q}}\newcommand{\bR}{\mathbb{R}}
\newcommand{\bS}{\mathbb{S}}\newcommand{\bT}{\mathbb{T}}

\newcommand{\dA}{\mathfrak{G}} 

\newcommand{\dI}{\mathds{I}}

\newcommand{\dO}{\mathds{O}}

\newcommand{\dS}{\mathds{S}}

\newcommand{\sR}{\mathscr{R}}
\newcommand{\sS}{\mathscr{S}}

\newcommand{\sZ}{\mathscr{Z}}

\newcommand{\sP}{\mathfrak{P}}

\DeclareMathOperator{\E}{\mathbb{E}}
\DeclareMathOperator{\pr}{\mathbb{P}}

\DeclareMathOperator{\var}{Var}

\DeclareMathOperator{\dis}{dis}

\DeclareMathOperator{\GH}{GH}
\DeclareMathOperator{\GHP}{GHP}
\DeclareMathOperator{\diam}{diam}
\DeclareMathOperator{\height}{ht}
\DeclareMathOperator{\ord}{ord}
\DeclareMathOperator{\dist}{dist}
\DeclareMathOperator{\exec}{exc}
\DeclareMathOperator{\md}{mod}

\DeclareMathOperator{\mass}{mass}
\DeclareMathOperator{\con}{con}
\DeclareMathOperator{\nr}{nr}
\DeclareMathOperator{\NR}{NR}
\DeclareMathOperator{\he}{ht}

\DeclareMathOperator{\udim}{\underline{dim}}
\DeclareMathOperator{\odim}{\overline{dim}}
\DeclareMathOperator{\res}{res}

\newcommand{\sss}{\scriptscriptstyle}
\newcommand{\erdos}{Erd\H{o}s-R\'enyi }
\newcommand{\ldown}{\ell^2_{\downarrow}}
\newcommand{\ldownthree}{\ell^3_{\downarrow}}

\newcommand{\barV}{\bar{V}}
\newcommand{\barL}{\bar{L}}
\newcommand{\Ep}{\E_{\vp}}
\newcommand{\Et}{\E_{\mvtheta}}

\newcommand{\Epst}{\E_{\vp,\star}}
\newcommand{\icrt}{\cT_{\sss(\infty)}^{\mvtheta}}
\newcommand{\tilicrt}{\cT_{\sss(\infty)}^{\mvtheta,\star}}
\newcommand{\convd}{\stackrel{d}{\longrightarrow}}

\definecolor{aqua}{rgb}{0.0, 1.0, 1.0}

\newcommand{\SSS}[1]{\todo[color=aqua, inline]{SS: #1}}

\newcommand{\eqn}[1]{\begin{equation} #1 \end{equation}}
\newcommand{\eqan}[1]{\begin{align} #1 \end{align}}

\newcommand {\vep}{\varepsilon}

\newcommand{\cluster}{\mathscr{C}}
\newcommand{\Poi}{{\sf Poi}}

\newcommand{\prob}{\mathbb{P}}
\newcommand{\expec}{\mathbb{E}}

\def\1{{\mathchoice {1\mskip-4mu\mathrm l}      
{1\mskip-4mu\mathrm l}
{1\mskip-4.5mu\mathrm l} {1\mskip-5mu\mathrm l}}}
\newcommand{\indic}[1]{\1_{\{#1\}}}

\newcommand{\nn}{\nonumber}
\newcommand{\e}{{\mathrm{e}}}

\newcommand{\NRnw}{{\rm NR}_n(\boldsymbol{w})}
\newcommand{\NRnwl}{{\rm NR}_n(\boldsymbol{w}(\lambda))}
\newcommand{\normi}[1]{\|#1\|_{\infty}}
\newcommand{\barg}{\bar{\gamma}}
\newcommand{\RC}{\cR\cC}
\DeclareMathOperator{\shape}{shape}
\DeclareMathOperator{\point}{pt}
\DeclareMathOperator{\osc}{osc}
\DeclareMathOperator{\unordered}{uo}
\newcommand{\tildQ}{\tilde{Q}}
\newcommand{\tildR}{\tilde{R}}

\newcommand{\HH}{\mathcal{H}}
\newcommand{\II}{\mathcal{I}}
\newcommand{\DD}{\mathcal{D}}
\newcommand{\RR}{\mathcal{R}}
\renewcommand{\SS}{\mathcal{S}}
\newcommand{\ZZ}{\mathcal{Z}}
\newcommand{\Exp}{{\sf Exp}}
\newcommand{\ch}[1]{\textcolor[rgb]{0,0,0}{{#1}}}
\newcommand{\op}{o_{\sss \mathrm{P}}}

\usepackage[refpage,noprefix,intoc]{nomencl}							
\makenomenclature											
\begin{document}

\title[ICRTs and critical random graphs]{The multiplicative coalescent, inhomogeneous continuum random trees, and new universality classes for critical random graphs}

\date{}

\author[Bhamidi]{Shankar Bhamidi$^1$}
\address{$^1$Department of Statistics and Operations Research, 304 Hanes Hall, University of North Carolina, Chapel Hill, NC 27599}
\author[van der Hofstad]{Remco van der Hofstad$^2$}
\author[Sen]{Sanchayan Sen$^3$}
\address{$^2$Department of Mathematics and Computer Science, Eindhoven University of Technology, the Netherlands}
\address{$^3$Department of Mathematics and Statistics, McGill University, 805 Sherbrooke Street West, Montr\'{e}al, QC, Canada}
\email{bhamidi@email.unc.edu, r.w.v.d.hofstad@tue.nl, sanchayan.sen1@gmail.com}

\maketitle
\begin{abstract}
One major open conjecture in the area of critical random graphs, formulated by statistical physicists, and supported by a large amount of numerical evidence over the last decade \cite{BraBulCohHavSta03,wu2006transport,braunstein2007optimal,chen2006universal} is as follows: for a wide array of random graph models with degree exponent $\tau\in (3,4)$, distances between typical points both within maximal components in the critical regime as well as on the minimal spanning tree on the giant component in the supercritical regime scale like $n^{(\tau-3)/(\tau-1)}$.

In this paper we study the metric space structure of maximal components of the multiplicative coalescent, in the regime where the sizes converge to excursions of L\'evy processes ``without replacement'' \cite{aldous-limic}, yielding a completely new class of limiting random metric spaces. A by-product of the analysis yields the continuum scaling limit of one fundamental class of random graph models with degree exponent $\tau\in (3,4)$ where edges are rescaled by $n^{-(\tau-3)/(\tau-1)}$ yielding the first rigorous proof of the above conjecture. The limits in this case are compact ``tree-like'' random fractals with a dense collection of hubs (infinite degree vertices), a finite number of which are identified with leaves to form shortcuts. In a special case, we show that the Minkowski dimension of the limiting spaces equal $(\tau-2)/(\tau-3)$ a.s., in stark contrast to the \erdos scaling limit whose Minkowski dimension is $2$ a.s. It is generally believed that dynamic versions of a number of fundamental random graph models, as one moves from the barely subcritical to the critical regime can be approximated by the multiplicative coalescent. In work in progress,  the general theory developed in this paper is used to prove analogous limit results for other random graph models with degree exponent $\tau\in (3,4)$.

Our proof makes crucial use of inhomogeneous continuum random trees (ICRT), which have previously arisen in the study of the entrance boundary of the \emph{additive} coalescent. We show that \emph{tilted} versions of the same objects using the associated mass measure, describe connectivity properties of the \emph{multiplicative} coalescent. Since convergence of height processes of corresponding approximating $\vp$-trees is not known, we use general methodology in \cite{AthLorWin14} and develop novel techniques relying on first showing convergence in the Gromov-weak topology and then extending this to Gromov-Hausdorff-Prokhorov convergence by proving a global lower mass-bound.

\end{abstract}



\section{Introduction and results}
\label{sec:intro}

In the last two decades many results regarding scaling limits of large discrete random objects to continuum analogs have been proved. Examples range from Aldous's continuum random tree \cite{Aldo91a,Aldo93b,legall-survey}, Schramm-Loewner evolution and critical planar systems \cite{schramm2000scaling}, to what is most closely related to this paper: scaling limits of maximal components in the critical regime for random graphs as well as the minimal spanning tree on the giant component in the supercritical regime \cite{BBG-12,BBG-limit-prop-11,AddBroGolMie13}.

Motivated by empirical observations on real-world networks, in the last decade, researchers from a wide array of fields including computer science, the social sciences and statistical physics have proposed a large number of random graph models to explain various functionals of real world systems including power law degree distributions and small world scaling of distances between nodes in the network \cite{bollobas2001random,durrett-rg-book,Hofs17,CL-book,newman2003structure,newman2010networks,albert2002statistical,dorogovtsev2002evolution}. Many of these models have a parameter $t$ related to the edge density and a model-dependent critical point $t_c$. Writing $n$ for the number of vertices in the network, if $t< t_c$ then the maximal connected component  $\cC_1(n)$ has size that is negligible \ch{compared} to $n$, while if $t> t_c$ one has a giant component $\cC_1(n)\sim f(t) n$ for some positive model-dependent function $f(t) > 0$ for $t> t_c$. The ``$t=t_c$'' regime is often referred to as the {\em critical regime}. Just as a study of the classical critical \erdos random graph spurred enormous activity in probabilistic combinatorics in the 90s \cite{janson2011random,bollobas2001random,luczak1990component,luczak1994structure,aldous-crit}, the study of the critical regime of these new random graph models and new phenomena such as \emph{explosive percolation} \cite{achlioptas2009explosive,riordan2011explosive} have motivated a concerted effort to understand the critical regime of these new random graph models.

   In this context, for more than a decade \cite{BraBulCohHavSta03,wu2006transport,braunstein2007optimal,chen2006universal}, one of the fundamental open conjectures in this area (loosely stated) is as follows.	Consider distances between typical points in the maximal component either in the critical regime or the minimal spanning tree on the giant component in the supercritical regime scale
	\begin{enumeratea}
		\item If the random graph model has an asymptotic degree distribution with finite third moments, then distances scale like $n^{1/3}$.
		\item If the random graph model has a limiting degree distribution $\set{p_k}_{k\geq 1}$ with tail $p_k\sim C/k^{\tau}$ for $\tau \in (3,4)$, then distances scale like $n^{(\tau-3)/(\tau-1)}$.
	\end{enumeratea}

\noindent {\bf Contributions of this paper:} Since we will need to setup some notation before getting to the main results, let us give a general overview of the contributions of this paper:

\begin{enumeratei}
	\item {\bf General theory:} The fundamental aim of the paper is to develop a general theory one can use to prove (b) in the conjecture above for a wide class of random graphs and, in particular, derive a new class of continuum scaling limits. To do so, we consider the multiplicative coalescent with entrance boundary in the space $l_0$ as in \cite{aldous-limic} (see \eqref{eqn:lthree-lo} below). Viewing the maximal components as measured metric spaces (using graph distance and vertex weights), we show that these components with edges and associated measures properly rescaled converge to continuum random objects in the Gromov weak sense. These resulting objects are obtained via appropriate tilts and vertex identifications of inhomogeneous continuum random trees; {\it untilted} versions of the same objects have been used to describe the entrance boundary of the \emph{additive} coalescent \cite{aldous-pitman-entrance}.  These resulting random objects are ``tree-like'' but with a dense collection of ``hubs'' (corresponding to infinite-degree vertices).
	\item {\bf Proof techniques:}  The standard technique in proving such results is to study height processes of certain spanning trees of the components and to show that these processes converge to limiting excursions that code the limiting random real trees. In our context, the convergence of height processes of the corresponding approximating $\vp$-trees is not known. In \cite{AMP}, the height processes of $\vp$-trees were shown to converge to limiting excursions in certain regimes, but these results are not applicable to our situation.

Because of this, we develop new techniques relying on first showing convergence in Gromov-weak topology via a careful analysis of the tree spanning a finite collection of ``typical'' points in random ``tilted" $\vp$-trees.  In one fundamental class of random graph models, we then extend Gromov-weak convergence to Gromov-Hausdorff-Prokhorov convergence by proving a global lower mass-bound.
\item {\bf Special case:} As an example of the general theory, we study the special case of the Norros-Reittu model \cite{NorRei06} (which in the regime of interest has been proven \cite{janson-equiv} to be equivalent to the Chung-Lu model \cite{CL-connected} and the rank-one random graph \cite{BJR07}). In this case, we show that the limiting spaces are compact. We also show that the box-counting or Minkowski dimension equals $(\tau-2)/(\tau-3)$ a.s.
\end{enumeratei}

In work in progress \cite{SB-SD-vdH-SS}, we use the general theory in this paper to analyze another fundamental random graph model, the configuration model with degree distribution with exponent $\tau\in (3,4)$, and derive the continuum analogs of the maximal components of this model. We defer a more detailed discussion of related work and the relevance of the current study to Section \ref{sec:discussion}.
\medskip

\noindent{\bf Organization of the paper:} A reasonable amount of notation regarding the entrance boundary of the multiplicative coalescent is required to describe the main results (Theorems \ref{thm:mc-main-1}, \ref{thm:mc-main-2-tau}). To ease the reader into the paper, we start in Section \ref{sec:rank-one} with the special case of the Norros-Reittu model and in Theorem \ref{thm:rank-one} describe what the main results imply for this model.  Then in Section \ref{sec:mc-results} we define the multiplicative coalescent as well as the class of entrance boundaries of importance for the paper and then describe the two main results. The results use two notions of convergence of metric spaces; these are given a precise formulation in Section \ref{sec:metric-convg}. Section \ref{sec:p-tree-ICRT-def} describes an important class of random trees called $\vp$-trees and the corresponding inhomogenous continuum random trees that arise as scaling limits of these objects. These are then used in Section \ref{sec:descp-limit} to give a precise description of the scaling limits of maximal components. We discuss the relevance of the main results, relate these to existing work and give an overview of the proof in Section \ref{sec:discussion}. The proofs of the main results are contained in Sections \ref{sec:proof1} - \ref{sec:fractal-proofs}.
\medskip

\noindent{\bf Notation:} Throughout this paper, we make use of the following standard notation.
We let $\convd$ denote convergence in distribution, and
$\probc$ convergence in probability. 
\nomenclature[conv]{$\convd, \probc$}{Convergence in distribution and probability} For a sequence of random variables
$(X_n)_{n\geq 1}$, we write $X_n=\op(b_n)$ when $|X_n|/b_n\probc 0$ as $n\rightarrow\infty$.
For a non-negative function $n\mapsto g(n)$,
we write $f(n)=O(g(n))$ when $|f(n)|/g(n)$ is uniformly bounded, and
$f(n)=o(g(n))$ when $\lim_{n\rightarrow \infty} f(n)/g(n)=0$.
Furthermore, we write $f(n)=\Theta(g(n))$ if $f(n)=O(g(n))$ and $g(n)=O(f(n))$.
We say  that a sequence of events $({\mathcal E}_n)_{n\geq 1}$ occurs with high probability (whp) when $\pr({\mathcal E}_n)\rightarrow 1$.

\subsection{Rank-one random graph}
\label{sec:rank-one}
\subsubsection{Model formulation}
We start by describing a particular class of random graph models called the Poissonian random graph or the Norros-Reittu model \cite{NorRei06,BJR07}, sometimes also referred to as the rank-one random graph model \cite{BJR07}. In the regime of interest for this paper, as shown in \cite{janson-equiv}, this model is equivalent to the Chung-Lu model \cite{CL-connected,CL-distance,CL-distances-2,CL-book} and the Britton-Deijfen-Martin-L{\"o}f model \cite{BDL06}.  Start with vertex set $[n]:=\set{1,2,\ldots, n}$ and suppose each vertex $i\in [n]$ has a weight $w_i\geq 0$ attached to it; intuitively this measures the propensity or attractiveness of this vertex in the formation of links. Writing $\mvw=(w_1,\hdots, w_n)$, place an edge between $i$ and $j$ independently for each $i\neq j\in [n]$ with probability
\begin{equation}
\label{eqn:nr-connection}
	q_{ij}=q_{ij}(\mvw):= 1-\exp(-w_i w_j/\ell_n),
\end{equation}
where $\ell_n$ is the total weight given by
\begin{equation}
\label{eqn:ln-def}
	\ell_n:= \sum_{i\in [n]} w_i.\notag
\end{equation}
To complete the formulation, we need to specify how these vertex weights are chosen. Essentially we want the empirical distribution of weights $n^{-1} \sum_{i\in [n]} \delta \set{w_i}$ to converge to a fixed pre-specified distribution $F$ as $n\to\infty$. There are a number of ways to do this, but for this paper the following choice turns out to be convenient for a clear statement of the results. Let $(w_i)_{i\in [n]}$ be constructed by
\begin{equation}
\label{eqn:wi-construc}
	w_i:= [1-F]^{-1}(i/n), \qquad i\in [n],
\end{equation}
where $F$ is a cumulative distribution function on $[0,\infty)$ and $[1-F]^{-1}$ is the generalized inverse
\begin{equation}
\label{eqn:f-inverse}
	[1-F]^{-1}(u):= \inf\set{s:[1-F](s)\leq u}.\notag
\end{equation}
We assume there exists \nomenclature[tau]{$\tau$}{Tail exponent of the cdf of the weight sequence $\boldsymbol{w}$} $\tau \in (3,4)$ and $c_{\sss F} > 0$ such that
\begin{equation}
\label{eqn:tau-c-def}
	\lim_{x\to\infty} x^{\tau-1} [1-F(x)]:= c_{\sss F}.
\end{equation}
We will use $W$ for a random variable with distribution $F$. We will use $\NRnw$ \nomenclature[NR]{$\NRnw$}{Norros-Reittu random graph with weight sequence $\boldsymbol{w}$.} to denote the corresponding random graph.

\subsubsection{Motivation and known results}
As described in the introduction, one impetus for the formulation of a wide array of network models, is to capture the heterogeneous and heavy-tailed nature of the degree distribution of empirical networks. Write $N_k$ for the number of vertices with degree $k$ in $\NRnw$.  Under the assumptions in the previous section, one can show \cite[Theorem 3.13]{BJR07} that
\begin{equation}
\label{eqn:degree-dist}
	\frac{N_k}{n} \probc \E\left(\e^{-W } \frac{W^k}{k!}\right), \qquad k\geq 0,
\end{equation}
where $W\sim F$. In particular, the degree distribution also has tail exponent $\tau$. More important in the context of this paper is the connectivity threshold. For $i\geq 1$ write $\cC_i$ for the $i$th largest connected component and let $|\cC_i|$ denote its number of vertices. Now define the parameter
\begin{equation}
\label{eqn:nu-def}
	\nu:= \frac{\E(W^2)}{\E(W)},
\end{equation}
\nomenclature[nu]{$\nu$}{Asymptotic expected forward degree Norros-Reittu random graph}
and note that $\nu<\infty$ by \eqref{eqn:tau-c-def}. Then by \cite[Theorem 3.1 and Section 16.4]{BJR07}, we have the following criterion for the phase transition for the largest component:
\begin{enumeratea}
	\item {\bf Supercritical regime:} If $\nu >1$, then there exists $\rho\in (0,1) $ such that  $|\cC_1|/n\probc\rho$ whilst $|\cC_2|/n\probc 0$;
	\item {\bf Subcritical regime:} If $\nu<1$, then $|\cC_1|/n\probc 0$.
\end{enumeratea}

The main aim of this paper is to understand the critical regime $\nu=1$ where also $|\cC_1|/n\probc 0$. In this setting, there are different universality classes depending on the vertex weights. In the \erdos or weakly inhomogeneous universality class, critical clusters have size of order $n^{2/3}$ and their metric space structure was discovered by Addario-Berry, Broutin and Goldschmidt \cite{BBG-12}. Interestingly, when $\E(W^3)<\infty$, component sizes still scale like $n^{2/3}$ \cite{bhamidi2010scaling} while assuming finite $6+\eps$-moments the metric space structure of rank-1 inhomogeneous random graphs is (apart from a trivial rescaling of size and time) the {\em same} \cite{SBSSXW14}. However, in the strongly inhomogeneous regime where $\E(W^3)=\infty$, the scaling limits of critical clusters are dramatically different in the sense that their sizes are gives by $n^{(\tau-2)/(\tau-1)},$ where $\tau$ is the degree power-law exponent given by \eqref{eqn:tau-c-def} \cite{SBVHVJL12, Hofs09a}. In this paper, we focus on their {\em metric space structure,} obtained after rescaling edges by $n^{-(\tau-3)/(\tau-1)}$ and taking the limit as $n\rightarrow \infty.$ We show that this limiting metric space is compact and its Minkowski dimension equals $(\tau-2)/(\tau-3)$, whereas the \erdos scaling limit has Minkowski dimension $2$.

In this paper, we analyze the entire \emph{critical scaling window}. Let $\mvw$ denote the weight sequence as in \eqref{eqn:wi-construc} and fix $\lambda \in \bR$. Now consider the weight sequence $\mvw(\lambda):=(w_i(\lambda))_{i\in [n]}$ \nomenclature[mvw]{$\mvw(\lambda):=(w_i(\lambda))_{i\in [n]}$}{Weight sequence in the critical scaling window.} defined by
\begin{equation}
\label{eqn:wilambda-def}
	\mvw(\lambda):= \left(1+\frac{\lambda}{n^{(\tau-3)/(\tau-1)}}\right)\mvw.
\end{equation}
Write $\NRnwl$ for the corresponding random graph and let $\cC_i(\lambda)$ \nomenclature[cci]{$\cC_i(\lambda)$}{The $i$-th largest component in $\NRnwl$.} denote the corresponding $i$th largest component.  Then this critical scaling window was first identified and studied in \cite{Hofs09a} where it was shown that for every fixed $\lambda \in \bR$, $|\cC_1|/n^{(\tau-2)/(\tau-1)}$ as well as $n^{(\tau-2)/(\tau-1)}/|\cC_1|$ are tight. The entire distributional asymptotics of component sizes were derived in \cite{SBVHVJL12} where it was shown that in the product topology on $\bR^{\bN}$,
\begin{equation}
\label{eqn:size-asymp-rank-one}
	\left(\frac{|\cC_i(\lambda)|}{n^{(\tau-2)/(\tau-1)}}\colon i\geq 1\right) \convd (Z_i(\lambda)\colon i\geq 1),
\end{equation}
where $(Z_i(\lambda)\colon i\geq 1)$ are excursions away from zero of a special stochastic process described in more detail in Section \ref{sec:mc-results}.

\subsubsection{Our results} We make the following convention:
\begin{quote}
For any metric measure space $(\bS, d, \mu)$ and $a>0$, $a\bS$ denotes the metric measure space $(\bS, ad, \mu)$, i.e, the space where the distance is scaled by $a$ and the measure remains unchanged.
\end{quote}
Consider the random graph $\NRnwl$ and view each connected component $\cC$ as a connected metric space via the usual graph distance where  each edge has length one. Further, we can view each connected component $\cC$ as a metric measure space by assigning weight $w_i/(\sum_{j\in\cC}w_j)$ to vertex $i\in\cC$. Note that the normalization yields a probability measure on each connected component. Let $\sS$ denote the space of (equivalence classes) of {\bf compact} measured metric spaces equipped with the Gromov-Hausdorff-Prokhorov metric (see Section \ref{sec:ghp} for definition). View
\begin{align}\label{eqn:def-M-n}
\vM_n^{\nr}(\lambda):= \big(\cC_i(\lambda)\colon i\geq 1\big)
\end{align}
as a random element of $\sS^{\bN}$. \nomenclature[vmn]{$\vM_n^{\nr}(\lambda)$}{Components of $\NRnwl$ viewed as an element of $\sS^{\bN}$. }

Next recall that the lower and upper box counting dimensions of a compact metric space $\cM$ are given by
\[\udim(\cM):= \liminf_{\delta\downarrow 0}\frac{\log{[\cN(\cM,\delta)]}}{\log(1/\delta)},
\quad\text{and}\quad\odim(\cM):= \limsup_{\delta\downarrow 0}\frac{\log{[\cN(\cM,\delta)]}}{\log(1/\delta)}\]
\nomenclature[uodim]{$\udim,\odim$}{Lower and upper box counting dimensions.}
respectively, where $\cN(\cM,\delta)$ is the minimal number of open balls with radius $\delta$ required to cover $\cM$. Also let $\dim_h(\cM)$ denote the Hausdorff dimension of $\cM$. When $\udim(\cM)=\odim(\cM)=\dim$, then the box-counting or Minkowski dimension is $\dim$.
\nomenclature[gen-met]{$\cM$}{Symbol used to denote a generic metric space.}
\nomenclature[cn-cm-del]{$\cN(\cM,\delta)$}{Minimal number of open balls with radius $\delta$ required to cover a metric space $\cM$.}

Before stating our main result, we introduce a technical condition.
\begin{ass}\label{ass:density}
The support of the limiting distribution $F$ (defined just before \eqref{eqn:wi-construc}) is given by $[\iota, \infty)$ for some $\iota>0$. Further, $F$ has a continuous density $f$ on $[\iota, \infty)$ such that $xf(x)$ is non-increasing on $[\iota, \infty)$.
\end{ass}

Note that distributions $F$ that are exact power laws, i.e., of the form $F(x)=1-(\iota/x)^{\tau-1}$ for $x>\iota$ and some $\tau\in(3, 4)$, satisfy Assumption \ref{ass:density}. The main result of this section is as follows:

\begin{thm}[Scaling limits with degree exponent $\tau\in (3,4)$]
\label{thm:rank-one}
	Fix $\lambda \in \bR$ and consider the critical Norros-Reittu model $\NRnwl$, i.e, assume that $\nu=1$ where $\nu$ is as in \eqref{eqn:nu-def}. Assume that the limiting distribution $F$ satisfies Assumption \ref{ass:density}.

Then, there exists an appropriate limiting sequence of random compact metric measure spaces $\vM_\infty^{\nr}(\lambda):= (M_i^{\nr}(\lambda))_{i\geq 1}$ such that the components in the critical regime satisfy
	\begin{align}\label{eqn:conergence-rank-one}
\frac{1}{n^{(\tau-3)/(\tau-1)}}\vM_n^{\nr}(\lambda) \convd \vM_\infty^{\nr}(\lambda), \qquad \mbox{ as } n\to\infty.
\end{align}
	Here convergence is with respect to the product topology on $\sS^{\bN}$ induced by the Gromov-Hausdorff-Prokhorov metric on each coordinate $\sS$. For each $i\geq 1$, the limiting metric spaces have the following properties:
	\begin{enumeratea}
		\item $M_i^{\nr}(\lambda)$ is random compact metric measure space obtained by taking a random {\bf real tree} $\cT_i(\lambda)$ and identifying a random (finite) number of pairs of points (thus creating shortcuts).
		\item Call a point $u\in \cT_i(\lambda)$ a~ {\bf hub} point if deleting the $u$ results in infinitely many disconnected components of $\cT_i(\lambda)$. Then $\cT_i(\lambda)$ has infinitely many hub points which are {\bf everywhere dense} on the tree $\cT_i(\lambda)$.   \item The box-counting or Minkowski dimension of $M_i^{\nr}(\lambda)$ satisfies
	\begin{align}\label{eqn:dim-nr}
	\dim(M_i^{\nr}(\lambda))=\frac{\tau-2}{\tau-3} \qquad a.s.
	\end{align}
	Consequently, the Hausdorff dimension satisfies the bound $\dim_h(M_i^{\nr}(\lambda)) \leq (\tau-2)/(\tau-3)$ a.s.
	\end{enumeratea}
\end{thm}

\begin{conj}\label{conj:fractal-dimension}
We strongly believe that both the Hausdorff dimension and the packing dimension of $M_i^{\nr}(\lambda)$ equal $(\tau-2)/(\tau-3)$ a.s. See Section \ref{sec:disc-open} for a discussion.
\end{conj}

\subsection{Connectivity asymptotics for the multiplicative coalescent}
\label{sec:mc-results}
In this section we consider a slightly more general setting than in Section \ref{sec:rank-one}. The motivation is as follows: recall that for the rank-one model, two vertices were connected with essentially probability proportional to the product of the weight between these two vertices. For probabilists, this connectivity pattern is quite reminiscent of the famous multiplicative coalescent \cite{aldous-crit,aldous-limic,bertoin-coagulation}. Whilst interesting in its own right, its fundamental importance in the context of random graphs is as follows: A wide array of random graph models can be constructed in a dynamic fashion where as time progresses new edges are created between pre-existing clusters. Even though the merging dynamics between connected components tend to be quite different from that specified by the multiplicative coalescent, the mergers from the barely subcritical regime through the critical scaling window can be approximated by the multiplicative coalescent. This idea was exploited in \cite{SBSSXW-universal}  to prove universality of scaling limits in the critical regime for several random graphs models.

Thus components at criticality of a wide array of random graph models can be thought of consisting of two major parts:

\begin{enumeratea}
	\item  ``Blobs" that are components formed in the barely subcritical regime.
	\item Edges formed between such blobs as the system proceeds from the barely subcritical regime through the critical scaling window.
\end{enumeratea}

The results below (in particular Theorem \ref{thm:mc-main-1}) specify how to handle the second aspect. In a companion paper we show how one can use macroscopic averaging of distances within blobs in random graph models such as the configuration model to show that these models also have the same scaling limit in the critical regime as Theorem \ref{thm:rank-one} in the setting where degrees obey power-laws with exponents $\tau\in(3, 4)$. Further, it will follow from Theorem \ref{thm:mc-main-1} that the convergence in \eqref{eqn:conergence-rank-one} holds with respect to the product topology induced by Gromov-weak topology on each coordinate. Therefore, Theorem \ref{thm:rank-one} can be recovered partially from the more general Theorem \ref{thm:mc-main-1} at the expense of working with a weaker topology.

Before stating the result we will need to define the multiplicative coalescent.  The natural domain of this Markov process is the space
\begin{equation}
\label{eqn:ldown}
	\ldown:= \Big\{\vx = (x_1, x_2, \ldots)\colon x_1\geq x_2\geq \cdots \geq 0, ~ \sum_i x_i^2 <\infty\Big\},
\end{equation}
equipped with the metric $d(\vx, \vy):= \sqrt{\sum_{i\geq 1} (x_i-y_i)^2}$. \nomenclature[ldown]{$\ldown$}{Space describing component sizes for the multiplicative coalescent} We will work in the simpler setup where the Markov process starts with a finite number of clusters, i.e, the process starts with $\vx \in \ldown$ such that $\exists n < \infty$ such that $x_i =0$ for $i> n$. Write $\ldown(n)$ for the collection of such vectors.  Now the Markov process $(\vX(t))_{t\geq 0}$ with initial state $\vX(0) = \vx$ evolves as follows. Write $\vX(t) = (X_i(t))_{i\geq 1}$. Then for $i\neq j$, clusters $i$ and $j$ merge at rate $X_i(t)\cdot X_j(t)$ into a single cluster of size $X_i(t)+X_j(t)$.

Note that for any fixed time $t>0$, it is easy to find the distribution of masses $\vX(t)$ via the following random graph:
\begin{defn}[Random graph $\cG_n(\vx,t)$]
	\label{def:cgvx-q}
	 Consider the vertex set $[n]:=\set{1,2,\ldots, n}$ and assign weight $x_i$ to vertex $i$. Now connect each pair of vertices $i, j$ with $i\neq j$ independently with probability
	\begin{equation}
	\label{eqn:mc-connection}
		q_{ij}:= 1-\exp(- tx_i x_j).
	\end{equation}
	Call this random graph $\cG_n(\vx,t)$. For a connected component $\cC\subseteq \cG_n(\vx,t)$, let $\mass(\cC):= \sum_{i\in \cC} x_i$. Let $(\cC_i(t))_{i\geq 1}$ denote the connected components arranged in decreasing order of their masses.
\end{defn}
The following is obvious from the definition of the multiplicative coalescent:

\begin{lemma}\label{lem:mc-comp-sizes}
	For each fixed $t\geq 0$, the masses of the multiplicative coalescent at time $t$ started with finite number of initial clusters with masses $\vx$ satisfies
	\[(X_i(t)\colon i\geq 1) \stackrel{d}{=} \big(\mass(\cC_i(t))\colon i\geq 1\big). \]
\end{lemma}
Analogous to \eqref{eqn:ldown}, consider the two spaces
\begin{equation}
\label{eqn:lthree-lo}
	\ldownthree:= \bigg\{\vc:= (c_1, c_2,\ldots)\ :\ c_1\geq c_2\geq \cdots \ch{\geq}~ 0,~ \sum_{i\geq 1} c_i^3 < \infty\bigg\}, \qquad l_0:= \ldownthree \setminus \ldown.
\end{equation}
\nomenclature[ell-three]{$\ldownthree$}{Decreasing positive vectors with finite $\ell^3$-norm}
These spaces turn out to be crucial in describing the entrance boundary of the eternal multiplicative coalescent  in \cite{aldous-limic}. In the context of this paper, we are interested in studying scaling limits of connected components of the random graph $\cG_n(\vx,t)$ when (suitably normalized) asymptotics of the weight vector $\vx$ are described by a vector $\vc\in l_0$. Let
\begin{equation}
\label{eqn:moments-def}
	\sigma_r(\vx):= \sum_i x_i^r, \qquad 1\leq r\leq 3.
\end{equation}
\nomenclature[sigmr]{$\sigma_r(\vx)$}{$r$th moment of the weight sequence $\vx$.}
We will make the following assumptions about the weight vector $\vx:=\vx(n)$ used to form the graph $\cG_n(\vx,t)$. These place the associated graph in a particular entrance boundary of the associated eternal multiplicative coalescent \cite[Proposition 7]{aldous-limic}.

\begin{ass}
	\label{ass:mc-assumptions}
	For each $n\geq 1$, let\hskip4pt $\vx^{\sss(n)} = (x_i^{\sss(n)}: 1\leq i\leq n)$ be an initial finite-length vector belonging to $\ldown(n)$. Suppose that as $n\to\infty$ there exists $\vc\in l_0$ such that
	\begin{align}
		\frac{\sigma_3(\vx^{\sss(n)})}{(\sigma_2(\vx^{\sss(n)}))^3} &\to \sum_j c_j^3, \label{eqn:sigma2-sigma3}\\
		\frac{x_j^{\sss(n)}}{\sigma_2(\vx^{\sss(n)}))} &\to c_j\ \text{ for } j\geq 1,\ \text{ and} \label{eqn:xj-to-cj}\\
		\sigma_2(\vx^{\sss(n)}) &\to 0. \label{eqn:sigma2-zero}
	\end{align}
\end{ass}
Now let $\set{\xi_j:j\geq 1}$ be a sequence of independent exponential random variables where $\xi_j$ has rate $c_j$ for each $j\geq 1$. For a fixed $\lambda \in \bR$, consider the process
\begin{equation}
\label{eqn:vc-process-def}
	V^{\vc}_\lambda(s):= \lambda s + \sum_j (c_j \ind\set{\xi_j \leq s} - c_j^2 s), \qquad s\geq 0.
\end{equation}
\nomenclature[vcl]{$V^{\vc}_\lambda(\cdot)$}{Levy process ``without replacement''. The corresponding process reflected at zero is $\tilde{V}^{\vc}_\lambda(\cdot)$. }
It turns out that this process is well defined precisely if $\vc\in \ldownthree$ \cite{aldous-limic}. Consider the ``reflected at zero'' process
\begin{equation}
\label{eqn:vc-reflection}
	\tilde{V}^{\vc}_\lambda(s):= V^{\vc}(s) - \min_{0\leq s^\prime \leq s} \tilde{V}^{\vc}(s^\prime),
\end{equation}
and the excursions of $\tilde{V}^{\vc}_\lambda(\cdot)$ from zero.  Then Aldous and Limic \cite{aldous-limic} showed that the lengths of these excursions are a.s. in $l^2$ precisely when $\vc\in l_0$, and thus can be arranged in decreasing order. Write
\begin{equation}
\label{eqn:excursions-def}
	\sZ(\lambda):= (\cZ_i(\lambda)\colon i\geq 1)
\end{equation}
for these excursions in decreasing order of their length. Let $Z_i(\lambda):= |\cZ_i(\lambda)|$ denote the length of the $i$th largest excursion and let
\begin{equation}
\label{eqn:vz-excursions-def}
	\vZ(\lambda):= (Z_i(\lambda)\colon i\geq 1) \in \ldown \qquad \mbox{a.s.}
\end{equation}
\nomenclature[vz]{$\vZ(\lambda)$}{Lengths of excursion of $\tilde{V}^{\vc}_\lambda(\cdot)$ from zero.}
Then Aldous and Limic  \cite{aldous-limic} proved the following result:

\begin{thm}[{\cite[Proposition 7]{aldous-limic}}]
\label{thm:aldous-limic}
	Fix $\lambda \in \bR$ and consider the time scale $t_n:= \lambda + [\sigma_2(\vx^{\sss(n)})]^{-1}$.
	 Under Assumptions \eqref{eqn:sigma2-sigma3}, \eqref{eqn:xj-to-cj}, \eqref{eqn:sigma2-zero}, the masses of the connected components of the graph $\cG_n(\vx,t_n)$ satisfy
	 \[\left(\mass\left[\cC_i\left(\lambda + \frac{1}{\sigma_2(\vx^{\sss(n)})}\right)\right]\colon i\geq 1\right) \convd \vZ(\lambda), \qquad \mbox{ as } n\to \infty,\]
with respect to the topology in $\ldown$, where $\vZ(\lambda)$ is as in \eqref{eqn:vz-excursions-def}.
\end{thm}

Now consider the connected components in $\cG_n(\vx,t)$, and as before, view each component $\cC$ as a connected metric space via the usual graph distance where each edge has length one. Further, view each component $\cC$ as a measured metric space by assigning mass $x_i/\mass(\cC)$ to each vertex $i\in\cC$. Let $\sS_{*}$ denote the space of (equivalence classes) of measured metric spaces equipped with Gromov-weak topology (see Section \ref{sec:gromov-weak} for definition) and view
\[\vM_n(\lambda):= \left(\cC_i\bigg(\lambda + \frac{1}{\sigma_2(\vx^{\sss(n)})}\bigg)\colon i\geq 1\right)\]
as a random element in $\sS_{*}^{\bN}$. Then our next result is about Gromov-weak convergence of $\vM_n(\lambda)$.

\begin{thm}\label{thm:mc-main-1}
	Fix $\lambda \in \bR$. Then under Assumption \ref{ass:mc-assumptions}, there exist an appropriate limiting sequence of metric spaces $\vM_{\infty}^{\vc}(\lambda):= (M_i^{\vc}(\lambda)\colon i\geq 1)$ such that
	\[\sigma_2(\vx^{\sss(n)})\vM_n(\lambda) \convd \vM_\infty^{\vc}(\lambda), \qquad \mbox{ as } n\to \infty.\]
	Here weak convergence is on $\sS_{*}^{\bN}$ which is equipped with the natural product topology induced by the Gromov-weak topology on each coordinate $\sS_{*}$.
\end{thm}
\begin{rem}
	A full description of the limit objects is given in Section \ref{sec:descp-limit}. The limit objects use tilted versions of inhomogeneous continuum random trees and checking compactness even of the original versions at this level of generality turns out to be quite intractable. However as the next theorem shows, in the special case of relevance to the rank-one model, one can prove much more.
\end{rem}

Consider the special sequence $\vc= \vc(\alpha,\tau):=(c_i(\alpha,\tau):\ i\geq 1)\in l_0$ with $\tau \in (3,4)$ and $\alpha > 0$, where
\begin{equation}
\label{eqn:c-tau-lamb-def}
	c_i(\alpha,\tau):= \frac{\alpha}{i^{1/(\tau-1)}}, \qquad i\geq 1.
\end{equation}

Then we have the following result about the limiting metric spaces:
\begin{thm}\label{thm:mc-main-2-tau}
	\label{cor:fractal-dim}
	Fix $\alpha >0$, $\tau\in (3,4)$ and let $\vc = \vc(\alpha,\tau)$ as in \eqref{eqn:c-tau-lamb-def}. Consider the limiting metric spaces $\vM_\infty^{\vc}(\lambda):=(M_i^{\vc}(\lambda)\colon i\geq 1)$.

Then almost surely $M_i^{\vc}(\lambda)$ is compact for every $i\geq 1$. Further, the Minkowski dimension of $M_i^{\vc}(\lambda)$ satisfies
	\begin{align}\label{eqn:dim-mc}
	\dim(M_i^{\vc}(\lambda)) = \frac{\tau-2}{\tau-3} \qquad a.s.
	\end{align}
	Consequently, the Hausdorff dimension satisfies the bound $\dim_h(M_i^{\vc}(\lambda)) \leq (\tau-2)/(\tau-3)$ a.s.
\end{thm}
\begin{rem}
Since we are dealing with equivalence classes of metric spaces (see Sections \ref{sec:ghp} and \ref{sec:gromov-weak}), Theorem \ref{thm:mc-main-2-tau} should be understood as claiming the existence of representative spaces $M_i^{\vc}(\lambda)$ that are compact, and satisfy the said conditions about the fractal dimensions. We will only work with these representative spaces throughout this paper.
\end{rem}

\section{Definitions and limit objects}

\subsection{Convergence of metric spaces}
\label{sec:metric-convg}
Proper notions of convergence of (measured) metric spaces is one of the central themes in this paper. Here we define the two topologies used in the statement of our results. We mainly follow \cite{Winter-gromov-weak,EJP2116,metric-geometry-book,gromov-book}.

\subsubsection{Gromov-Hausdorff-Prokhorov metric}\label{sec:ghp}
In this section, all metric spaces under consideration will be compact metric spaces with associated probability measures. Let us first recall the Gromov-Hausdorff distance $d_{\GH}$ between metric spaces.  Fix two metric spaces $(X_1,d_1)$ and $(X_2, d_2)$. For a subset $C\subseteq X_1 \times X_2$, the {\it distortion} of $C$ is defined as
\begin{equation}
	\label{eqn:def-distortion}
	\dis(C):= \sup \set{|d_1(x_1,y_1) - d_2(x_2, y_2)|: (x_1,x_2) , (y_1,y_2) \in C}.
\end{equation}
\nomenclature[dis]{$\dis(C)$}{Distortion of correspondence $C\subseteq X_1 \times X_2$. }
A {\it correspondence} $C$ between $X_1$ and $X_2$ is a measurable subset of $X_1 \times X_2$ such that for every $x_1 \in X_1$ there exists at least one $x_2 \in X_2$ such that $(x_1,x_2) \in C$ and vice-versa. The {\it Gromov-Hausdorff distance} between the two metric spaces  $(X_1,d_1)$ and $(X_2, d_2)$ is defined as
\begin{equation}
\label{eqn:dgh}
	d_{\GH}(X_1, X_2) = \frac{1}{2}\inf \set{\dis(C): C \mbox{ is a correspondence between } X_1 \mbox{ and } X_2}.
\end{equation}
\nomenclature[dgh]{$d_{\GH}(X_1, X_2)$}{Gromov-Hausdorff distance between two metric spaces  $(X_1,d_1)$ and $(X_2, d_2)$. See same section for pointed Gromov-Hausdorff distance $d_{\GH}^{\point}$. }

Suppose $(X_1, d_1)$ and $(X_2, d_2)$ are two metric spaces and $p_1\in X_1$, and $p_2\in X_2$. Then the {\it pointed Gromov-Hausdorff distance} between $\mvX_1:=(X_1, d_1, p_1)$ and $\mvX_2:=(X_2, d_2, p_2)$ is given by
\begin{align}
\label{eqn:dgh-pointed}
	d_{\GH}^{\point}(\mvX_1, \mvX_2) = \frac{1}{2}\inf \set{\dis(C): C \mbox{ is a correspondence between }X_1 \mbox{ and } X_2\mbox{ and }(p_1, p_2)\in C}.
\end{align}

We will need a metric that also keeps track of associated measures on the corresponding spaces. A compact measured metric space $(X, d , \mu)$ is a compact metric space $(X,d)$ with an associated probability measure $\mu$ on the Borel sigma algebra $\cB(X)$. Given two compact measured metric spaces $(X_1, d_1, \mu_1)$ and $(X_2,d_2, \mu_2)$ and a measure $\pi$ on the product space $X_1\times X_2$, the {\it discrepancy} of $\pi$ with respect to $\mu_1$ and $\mu_2$ is defined as
\begin{equation}
	\label{eqn:def-discrepancy}
	D(\pi;\mu_1, \mu_2):= ||\mu_1-\pi_1|| + ||\mu_2-\pi_2||,
\end{equation}
where $\pi_1, \pi_2$ are the marginals of $\pi$ and $||\cdot||$ denotes the total variation distance between probability measures. Then the {\it Gromov-Haussdorf-Prokhorov distance} between $X_1$ and $X_2$ is defined as
\begin{equation}
\label{eqn:dghp}
	d_{\GHP}(X_1, X_2):= \inf\set{ \max\left(\frac{1}{2} \dis(C),~D(\pi;\mu_1,\mu_2),~\pi(C^c)\right) },
\end{equation}
where the infimum is taken over all correspondences $C$ and measures $\pi$ on $X_1 \times X_2$.
\nomenclature[dghp]{$d_{\GHP}(X_1, X_2)$}{Gromov-Hausdorff distance between two measured metric spaces  $(X_1,d_1,\mu_1)$ and $(X_2, d_2, \mu_2)$.}

Similar to \eqref{eqn:dgh-pointed}, we can define a ``{\it pointed Gromov-Hausdorff-Prokhorov distance}'', $d_{\GHP}^{\point}$ between two metric measure spaces $X_1$ and $X_2$ having two distinguished points $p_1$ and $p_2$ respectively by taking the infimum in \eqref{eqn:dghp} over all correspondences $C$ and measures $\pi$ on $X_1 \times X_2$ such that $(p_1, p_2)\in C$.

Write $\sS$ for the collection of all measured compact metric spaces $(X,d,\mu)$.
\nomenclature[sps]{$\sS$}{Space of all measured compact metric spaces. $\bar \sS$ is the corresponding space of isometry equivalent classes under $d_{\GHP}$}
 The function $d_{\GHP}$ is a pseudometric on $\sS$, and defines an equivalence relation $X \sim Y \Leftrightarrow d_{\GHP}(X,Y) = 0$ on $\sS$. Let $\bar \sS := \sS / \sim $ be the space of isometry equivalent classes of measured compact metric spaces and $\bar d_{\GHP}$ the induced metric. Then by \cite{EJP2116}, $(\bar \sS, \bar d_{\GHP})$ is a complete separable metric space. To ease notation, we will continue to use $(\sS, d_{\GHP})$ instead of $(\bar \sS, \bar d_{\GHP})$ and $X = (X, d, \mu)$ to denote both the metric space and the corresponding equivalence class.

\subsubsection{Gromov-weak topology}\label{sec:gromov-weak}
Here we mainly follow \cite{Winter-gromov-weak}. Introduce an equivalence relation on the space of complete and separable metric spaces that are equipped with a probability measure on the associated Borel $\sigma$-algebra by declaring two such spaces $(X_1, d_1, \mu_1)$ and $(X_2, d_2, \mu_2)$ to be equivalent when there exists an isometry $\psi:\mathrm{support}(\mu_1)\to\mathrm{support}(\mu_2)$ such that $\mu_2=\psi_{\ast}\mu_1:=\mu_1\circ\psi^{-1}$, i.e., the push-forward of $\mu_1$ under $\psi$ is $\mu_2$. Write $\sS_{*}$ for the associated space of equivalence classes. As before, we will often ease notation by not distinguishing between a metric space and its equivalence class.

Fix $m\geq 2$, and a complete separable metric space $(X, d)$. Then given a collection of points $\vx:=(x_1, x_2, \ldots, x_m)\in X^m$, let $\vD(\vx):= (d(x_i, x_j))_{i,j\in [m]}$ denote the symmetric matrix of pairwise distances between the collection of points. A function $\Phi\colon \sS_* \to \bR$ is called a polynomial of degree $m$ if there exists a bounded continuous function $\phi\colon \bR_+^{n^2}\to \bR$ such that
\begin{equation}\label{eqn:polynomial-func-def}
	\Phi((X,d,\mu)):= \int  \phi(\vD(\vx)) \mu^{\otimes m}(d(\vx)).
\end{equation}
Here $\mu^{\otimes m}$ is the $m$-fold product measure of $\mu$. Let $\mvPi$ denote the space of all polynomials on $\sS_*$.

\begin{defn}[Gromov-weak topology]
	\label{def:gromov-weak}
	A sequence $(X_n, d_n, \mu_n)_{n\geq 1} \in \sS_*$ is said to converge to $(X, d, \mu) \in \sS_*$ in the Gromov-weak topology if and only if $\Phi((X_n, d_n, \mu_n))\to \Phi((X, d, \mu))$ for all $\Phi\in \mvPi$.
\end{defn}
In \cite[Theorem 1]{Winter-gromov-weak} it is shown that $\sS_*$ is a Polish space under the Gromov-weak topology.
\nomenclature[sps]{$\sS_*$}{Space of measured metric spaces under the Gromov weak topology.} 
It is also shown that, in fact, this topology can be completely metrized using the so-called Gromov-Prokhorov metric.

\subsubsection{Spaces of trees with edge lengths, leaf weights and root-to-leaf measures}
\label{sec:space-of-trees}
In the proof of the main results we need the following two spaces built on top of the space of discrete trees. The first space $\vT_{IJ}$ was formulated in \cite{aldous-pitman-edge-lengths,aldous-pitman-entrance} where it was used to study trees spanning a finite number of random points sampled from an inhomogeneous continuum random tree (as described in the next section). We use the same notation in this paper.

\nomenclature[vtij]{$\vT_{IJ}$}{Space of tress with $I$ leaves all labelled, and $J$ other labeled ``hub'' vertices and further every edge has strictly positive edge length.}
\noindent{\bf The space $\vT_{IJ}$:} Fix $I\geq 0$ and $J\geq 1$. Let $\vT_{IJ}$ be the space of trees having the following properties:
\begin{enumeratea}
	\item There are exactly $J$ leaves labeled $1+, \ldots, J+$, and the tree is rooted at another labeled vertex $0+$.
	\item There may be extra labeled vertices (called hubs) with distinct labels in $\set{1,2,\ldots, I}$. (It is possible that only some, and not all labels in $\set{1,2,\ldots, I}$ are used.)
	\item Every edge $e$ has a strictly positive edge length $l_e$.
\end{enumeratea}
A tree $\vt\in \vT_{IJ}$ can be viewed as being composed of two parts:\\
(1) $\shape(\vt)$ describing the shape of the tree (including the labels of leaves and hubs) but ignoring edge lengths. The set of all possible shapes $\vT_{IJ}^{\shape}$ is obviously finite for fixed $I, J$.\\
(2) The edge lengths $\vl(\vt):= (l_e:e\in \vt)$. Consider the product topology on $\vT_{IJ}$ consisting of the discrete topology on $\vT_{IJ}^{\shape}$ and the product topology on $\bR^m$ where \ch{$m$ is the number of edges of $\vt$}.

\nomenclature[vtijs]{$\vT_{IJ}^*$}{The pace $\vT_{IJ}$ where in addition the trees are equipped with leaf weights and root-to-leaf measures. }
\noindent{\bf The space $\vT_{IJ}^*$:} We will need a slightly more general space. Along with the three attributes above in $\vT_{IJ}$, the trees in this space have the following two additional properties. Let $\cL(\vt):= \set{1+, \ldots, J+}$ denote the collection of non-root leaves in $\vt$. \nomenclature[ltv]{$\cL(\vt)$}{The collection of non-root leaves in a tree $\vt$.} Then every leaf $v\in \cL(\vt) $ has the following attributes:

\begin{enumeratea}
	\item[(d)] {\bf Leaf weights:} A strictly positive number $A(v)$. Write $\vA(\vt):=(A(v): v\in \cL(\vt))$.
	\item[(e)] {\bf Root-to-leaf measures:} A probability measure $\nu_{\vt,v}$ on the path $[0+,v]$ connecting the root and the leaf $v$. Here the path is viewed as a line segment pointed at $0+$ and has the usual Euclidean topology. Write $\mvnu(\vt):= (\nu_{\vt,v}: v\in \cL(\vt))$ for this collection of probability measures.
\end{enumeratea}
In addition to the topology on $\vT_{IJ}$, the space $\vT_{IJ}^*$ with these additional two attributes inherits the product topology on $\bR^{J}$ owing to leaf weights and $(d_{\GHP}^{\point})^J$ owing to the root-to-leaf measures.

For consistency, we add to the spaces $\vT_{IJ}$ and $\vT_{IJ}^*$ a conventional state $\partial$. Its use will be clear later on.

\subsection{Random $\vp$-trees and inhomogeneous continuum random trees (ICRTs)}
\label{sec:p-tree-ICRT-def}

For fixed $m \geq 1$, write $\bT_m$ and $\bT_m^{\ord}$ for the collection of all rooted trees with vertex set $[m]$ and rooted ordered trees with vertex set $[m]$ respectively. \nomenclature[tmtmor]{$\bT_m,\bT_m^{\ord}$}{Collection of all rooted (respectively rooted ordered) trees  with vertex set $[m]$.} Here we will view a rooted tree as being directed with the root being the original progenitor and each edge being directed from child to parent. An ordered rooted tree is a tree where children of each individual are assigned an order (meant to describe for example orientation in a planar embedding, e.g.,  right to left or some notion of age, e.g., oldest to youngest).

In this section, we define a family of random tree models called $\vp$-trees \cite{pitman-camarri,pitman-random-mappings}, and their corresponding limits, the so-called inhomogeneous continuum random trees, which play a key role in describing the limit metric spaces as well as in the proof.   Fix $m \geq 1$, and a probability mass function $\vp = (p_1, p_2,\ldots, p_m)$ with $p_i > 0$ for all $i\in [m]$. A $\vp$-tree is a random tree in $\bT_m$, with law as follows. For any fixed $\vt \in \bT_m$ and $v\in \vt$, write $d_v(\vt)$ for the number of children of $v$ in the tree $\vt$. Then the law of the $\vp$-tree, denoted by $\pr_{\text{tree}}$, is defined as:
\begin{equation}
\label{eqn:p-tree-def}
	\pr_{\text{tree}}(\vt) = \pr_{\text{tree}}(\vt; \vp) = \prod_{v\in [m]} p_v^{d_v(\vt)}, \qquad \vt \in \bT_m.
\end{equation}
\nomenclature[ptree]{$\pr_{\text{tree}}(\cdot; \vp)$}{Distribution of a $\vp$-tree with driving pmf $\vp$. }
Generating a random $\vp$-tree $\cT\sim \pr_{\text{tree}}$ and then assigning a uniform random order on the children of every vertex $v\in \cT$ gives a random element with law $\pr_{\ord}(\cdot ; \vp)$ given by
\begin{equation}
\label{eqn:ordered-p-tree-def}
	\pr_{\ord}(\vt) = \pr_{\ord}(\vt; \vp) = \prod_{v\in [m]} \frac{p_v^{d_v(\vt)}}{(d_v(\vt)) !}, \qquad \vt \in \bT_m^{\ord}.
\end{equation}
Obviously a $\vp$-tree can be constructed by first generating an ordered $\vp$-tree with the above distribution and then forgetting about the order.

In a series of papers \cite{aldous-pitman-entrance,AMP,aldous-pitman-edge-lengths} it was shown that $\vp$-trees, under various assumptions, converge to inhomogeneous continuum random trees that we now describe. Recall the space $\ldown$ in \eqref{eqn:ldown}. Consider the subset $\Theta\subset \ldown$ given by
\begin{equation}
\label{eqn:Theta-def}
	\Theta:= \bigg\{\mvtheta:=(\theta_i\colon i\geq 1)\in \ldown: \sum_{i=1} \theta_i =\infty,~  \sum_{i=1}^\infty \theta_i^2 = 1\bigg\}.
\end{equation}
\nomenclature[Theta]{$\Theta$}{Space of tenable parameters giving rise to ICRTs. }
Now recall from \cite{legall-survey,evans-book} that a real tree is a metric space $(\cT,d)$ that satisfies the following for every pair $a,b\in \cT$:
\begin{enumeratea}
	\item There is a {\bf unique} isometric map $f_{a,b}\colon [0,d(a,b)]\to \cT$ such that $f_{a,b}(0)=a,~ f_{a,b}(d(a,b)) =b$.
	\item For any continuous one-to-one map $g:[0,1]\to \cT$ with $g(0)=a$ and $g(1)=b$, we have $g([0,1]) = f_{a,b}([0,d(a,b)])$.
\end{enumeratea}

\noindent {\bf Construction of the ICRT:}
Given $\mvtheta\in \Theta$, we will now define the inhomogeneous continuum random tree $\cT^{\mvtheta}_{\sss(\infty)}$.
\nomenclature[ICRT]{$\cT^{\mvtheta}_{\sss(\infty)}$}{An ICRT constructed using $\mvtheta\in \Theta$.}
 We mainly follow the notation in \cite{aldous-pitman-entrance}. Assume that we are working on a probability space $(\Omega, \cF,\pr_{\mvtheta})$ rich enough to support the following:
\begin{enumeratea}
	\item For each $i\geq 1$, let $\cP_i:= (\xi_{i,1}, \xi_{i,2}, \ldots)$ be a rate $\theta_i$ Poisson process, independent for different $i$. The first point of each process $\xi_{i,1}$ is special and is called a \emph{joinpoint}, whilst the remaining points $\xi_{i,j}$ with $j\geq 2$ will be called \emph{$i$-cutpoints} \cite{aldous-pitman-entrance}.
	\item Independent of the above, let $\mvU=(U_j^{\sss(i)}:j\geq 1,\ i\geq 1)$ be a collection of i.i.d.\ uniform $(0,1)$ random variables. These are not required to construct the tree but will be used to define a certain function on the tree.
\end{enumeratea}
  The {\bf random} real tree (with marked vertices) $\icrt$ is then constructed as follows:
\begin{enumeratei}
	\item  Arrange the cutpoints $\set{\xi_{i,j}\colon i\geq 1, j\geq 2}$ in increasing order as $0< \eta_1 < \eta_2 < \cdots$. The assumption that $\sum_i \theta_i^2 <\infty$ implies that this is possible. For every cutpoint $\eta_k=\xi_{i,j}$, let $\eta_k^*:=\xi_{i,1}$ be the corresponding joinpoint.
	\item Next, build the tree inductively. Start with the branch $[0,\eta_1]$. Inductively assuming we have completed step $k$, attach the branch $(\eta_k, \eta_{k+1}]$ to the joinpoint $\eta_k^*$ corresponding to $\eta_k$.
\end{enumeratei}
Write $\cT_0^{\mvtheta}$ for the corresponding tree after one has used up all the branches $[0,\eta_1], \set{(\eta_k, \eta_{k+1}]: k\geq 1}$.
Note that for every $i\geq 1$, the joinpoint $\xi_{i,1}$ corresponds to a vertex with infinite degree. Label this vertex $i$. The ICRT $\icrt$ is the completion of the marked metric tree $\cT^{\mvtheta}_0$. As argued in \cite[Section 2]{aldous-pitman-entrance}, this is a real-tree as defined above which can be viewed as rooted at the vertex corresponding to zero. We call the vertex corresponding to joinpoint $\xi_{i,1}$  {\bf hub} $i$. Since $\sum_i \theta_i = \infty$, one can check that hubs are almost everywhere dense on $\icrt$.

\begin{figure}[htbp]
	\centering
		\begin{tikzpicture}
		[scale=.65, background rectangle/.style=
		     {draw=blue!80,fill=blue!10,rounded corners=2ex},
		   show background rectangle]
		\draw[very thick,blue, |->] (0,0) -- (16,0);
		   \node at (-1,0) [draw = red!100, fill = blue!20] {$\mathcal{P}_1$};
			\node at ( 4,0) [circle,draw=red!50,fill=red!100, label = above: $\xi_{11}$] {};
			\node at ( 6,0) [circle,draw=black!50,fill=black!100, label = above: $\xi_{12}$] {};
			\node at ( 11,0) [circle,draw=black!50,fill=black!100, label = above: $\xi_{13}$] {};
			\node at ( 14,0) [circle,draw=black!50,fill=black!100, label = above: $\xi_{14}$] {};
			\node at ( 15,0) [circle,draw=black!50,fill=black!100, label = above: $\xi_{15}$] {};

			   \node at (-1,-2) [draw = red!100, fill = blue!20] {$\mathcal{P}_2$};
			\draw[very thick,blue, |->] (0,-2) -- (16,-2);
			\node at ( 1,-2) [star,draw=red!50,fill=red!100, label = above: $\xi_{21}$] {};
			\node at ( 5,-2) [star, draw=black!50,fill=black!100, label = above: $\xi_{22}$] {};
			\node at ( 9.5, -2) [star,draw=black!50,fill=black!100, label = above: $\xi_{23}$] {};

		   \node at (-1,-4) [draw = red!100, fill = blue!20] {$\mathcal{P}_3$};
			\draw[very thick,blue, |->] (0,-4) -- (16,-4);
			\node at ( 8,-4) [dart, draw=red!50,fill=red!100, label = above: $\xi_{31}$] {};

			   \node at (-1,-6) [draw = red!100, fill = blue!20] {$\mathcal{P}_4$};
			\draw[very thick,blue, |->] (0,-6) -- (16,-6);
	
				\node at ( 12,-6) [regular polygon, regular polygon sides= 6, draw=red!50,fill=red!100, label = below: $\xi_{41}$] {};
				\node at ( 13,-6) [regular polygon, regular polygon sides= 6, draw=black!50,fill=black!100, label = below: $\xi_{42}$] {};

			\draw[very thick,blue, |->] (0,-8) -- (16,-8);

			\node at ( 4,-8) [circle,draw=red!50,fill=red!100] {};
			\node at ( 6,-8) [circle,draw=black!50,fill=black!100] {};
			\node at ( 11,-8) [circle,draw=black!50,fill=black!100] {};
			\node at ( 14,-8) [circle,draw=black!50,fill=black!100] {};
			\node at ( 15,-8) [circle,draw=black!50,fill=black!100] {};
	
			\node at ( 1,-8) [star,draw=red!50,fill=red!100] {};
			\node at ( 5,-8) [star, draw=black!50,fill=black!100] {};
			\node at ( 9.5, -8) [star,draw=black!50,fill=black!100] {};
			\node at ( 8,-8) [dart, draw=red!50,fill=red!100] {};
	
			\node at ( 12,-8) [regular polygon, regular polygon sides= 6, draw=red!50,fill=red!100] {};
			\node at ( 13,-8) [regular polygon, regular polygon sides= 6, draw=black!50,fill=black!100] {};

		\end{tikzpicture}	
		
	\caption{An illustration of the ICRT construction with four point process $\set{\cP_i:1\leq i\leq 4}$. The red points represent the \emph{joinpoint} of the corresponding point process and the blue points the corresponding cutpoints. The last line contains the union of the four point processes. See Figure \ref{fig:tinf-theta} for the corresponding tree.  }
	\label{fig:point-process-icrt}
\end{figure}
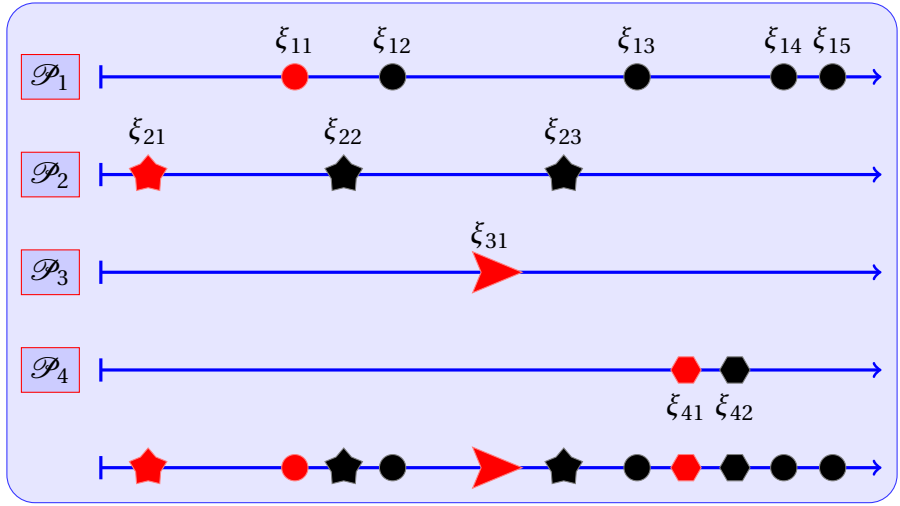

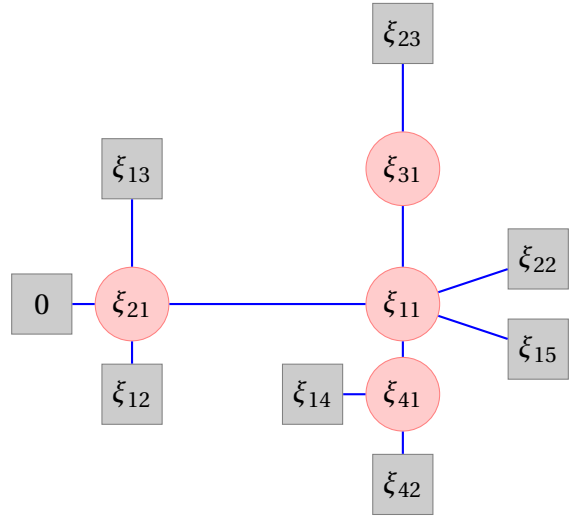
\begin{figure}[htbp]
	\begin{tikzpicture}[scale=.6]
	
	  \node (zero) [rectangle, draw=black!50,fill=black!20 , minimum size=8mm] at (0,0)   {$0$};
	  \node (x21) [circle, draw=red!50,fill=red!20 ] at (2,0) {$\xi_{21}$};
	  \node (x11) [circle, draw=red!50,fill=red!20] at (8,0) {$\xi_{11}$};
	  \node (x22) [rectangle, draw=black!50,fill=black!20 , minimum size=8mm] at (11,1) {$\xi_{22}$};
	  \node (x13) [rectangle, draw=black!50,fill=black!20 , minimum size=8mm] at (2,3) {$\xi_{13}$};
	  \node (x12) [rectangle, draw=black!50,fill=black!20, minimum size=8mm] at (2,-2) {$\xi_{12}$};
	  \node (x31) [circle, draw=red!50,fill=red!20] at (8,3) {$\xi_{31}$};
	  \node (x23) [rectangle, draw=black!50,fill=black!20, minimum size=8mm] at (8,6) {$\xi_{23}$};
	  \node (x41) [circle, draw=red!50,fill=red!20] at (8,-2) {$\xi_{41}$};
	  \node (x42) [rectangle, draw=black!50,fill=black!20, minimum size=8mm] at (8,-4) {$\xi_{42}$};
	  \node (x14) [rectangle, draw=black!50,fill=black!20, minimum size=8mm] at (6,-2) {$\xi_{14}$};
	  \node (x15) [rectangle, draw=black!50,fill=black!20 , minimum size=8mm] at (11,-1) {$\xi_{15}$};
    \draw[thick,blue] (zero) -- (x21) -- (x11) -- (x22);
    \draw[thick,blue] (x21) -- (x13);
    \draw[thick,blue] (x21) -- (x12);
    \draw[thick,blue] (x11) -- (x31);
    \draw[thick,blue] (x31) -- (x23);
    \draw[thick,blue] (x11) -- (x41);
    \draw[thick,blue] (x41) -- (x42);
    \draw[thick,blue] (x41) -- (x14);
    \draw[thick,blue] (x11) -- (x15);
	\end{tikzpicture}
	\caption{The tree constructed via the stick-breaking construction from Figure \ref{fig:point-process-icrt}.}
	\label{fig:tinf-theta}
\end{figure}

\begin{rem}\label{rem:planar-embedding-U-ij}
The uniform random variables $(U_j^{\sss(i)}:j\geq 1,\ i\geq 1)$ give rise to a natural ordering on $\icrt$ (or a planar embedding of $\icrt$) as follows. For $i\geq 1$, let $(\cT_j^{\sss(i)}:j\geq 1)$ be the collection of subtrees hanging off of the $i$th hub. Associate $U_j^{\sss(i)}$ with the subtree $\cT_j^{\sss(i)}$, and think of $\cT_{j_1}^{\sss(i)}$ appearing ``to the right of" $\cT_{j_2}^{\sss(i)}$ if $U_{j_1}^{\sss(i)}< U_{j_2}^{\sss(i)}$. This is the natural ordering on $\icrt$ when it is being viewed as a limit of ordered $\vp$-trees. We can think of the pair $(\icrt, \mvU)$ as the {\bf ordered ICRT}.
\end{rem}

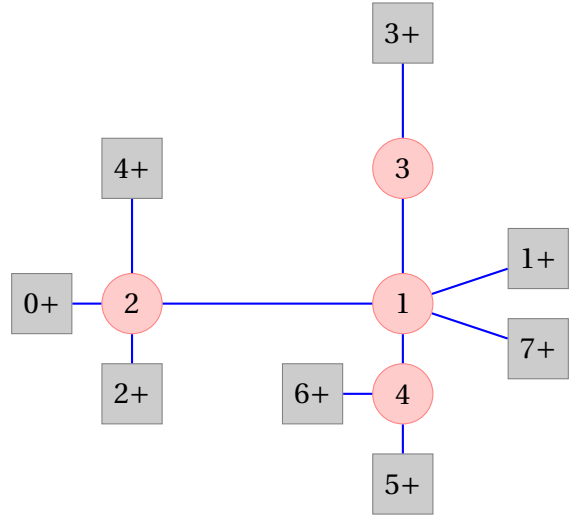
\begin{figure}[htbp]
	
	\begin{tikzpicture}[scale=.6]
	
	  \node (zero) [rectangle, draw=black!50,fill=black!20 , minimum size=8mm] at (0,0)   {$0+$};
	  \node (x21) [circle, draw=red!50,fill=red!20, minimum size=8mm ] at (2,0) {$2$};
	  \node (x11) [circle, draw=red!50,fill=red!20, minimum size=8mm] at (8,0) {$1$};
	  \node (x22) [rectangle, draw=black!50,fill=black!20 , minimum size=8mm] at (11,1) {$1+$};
	  \node (x13) [rectangle, draw=black!50,fill=black!20 , minimum size=8mm] at (2,3) {$4+$};
	  \node (x12) [rectangle, draw=black!50,fill=black!20, minimum size=8mm] at (2,-2) {$2+$};
	  \node (x31) [circle, draw=red!50,fill=red!20, minimum size=8mm] at (8,3) {$3$};
	  \node (x23) [rectangle, draw=black!50,fill=black!20, minimum size=8mm] at (8,6) {$3+$};
	  \node (x41) [circle, draw=red!50,fill=red!20, minimum size=8mm] at (8,-2) {$4$};
	  \node (x42) [rectangle, draw=black!50,fill=black!20, minimum size=8mm] at (8,-4) {$5+$};
	  \node (x14) [rectangle, draw=black!50,fill=black!20, minimum size=8mm] at (6,-2) {$6+$};
	  \node (x15) [rectangle, draw=black!50,fill=black!20 , minimum size=8mm] at (11,-1) {$7+$};
     \draw[thick,blue] (zero) -- (x21) -- (x11) -- (x22) ;
     \draw[thick,blue] (x21) -- (x13);
     \draw[thick,blue] (x21) -- (x12);
     \draw[thick,blue] (x11) -- (x31);
     \draw[thick,blue] (x31) -- (x23);
     \draw[thick,blue] (x11) -- (x41);
     \draw[thick,blue] (x41) -- (x42);
     \draw[thick,blue] (x41) -- (x14);
     \draw[thick,blue] (x11) -- (x15);
	\end{tikzpicture}

	\caption{Reduced tree \ch{$r_{47}^{\sss(\infty)}$} corresponding to the tree in Figure \ref{fig:tinf-theta}.}
	\label{fig:reduced-tree}
\end{figure}

\noindent{\bf Reduced tree $r_{IJ}^{\sss(\infty)}$:} Fix $I\geq 0$ and $J\geq 1$. Now let $\eta_0 = 0$ and for $j\geq 0$ call vertex $\eta_j$ the $j$th sampled leaf and label this as $j+$ to differentiate this from hub $j$. Note that the subtree of $\icrt$ spanned by $\set{0+,1+, \ldots, J+}$ (namely the part of the tree constructed from the interval $[0,\eta_J]$) is a tree in the usual sense with random edge lengths. For all hubs $i$, if $i\leq I$, retain its label and remove the label otherwise. This gives a random element of $\vT_{IJ}$ (recall the definiton Section \ref{sec:space-of-trees}), which we denote by $r_{IJ}^{\sss(\infty)}$. See Figure \ref{fig:reduced-tree} corresponding to the stick-breaking construction in Figures \ref{fig:point-process-icrt} and \ref{fig:tinf-theta}.

\noindent{\bf Mass measure:} For every vertex $v\in \icrt$, define the \emph{degree} of $v$ to be the number of connected components of $\icrt\setminus \set{v}$. Vertices with degree one are called \emph{leaves} of $\icrt$ and all other vertices form the \emph{skeleton} of the tree. Let $\cL(\icrt)$ denote the set of leaves of $\icrt$. In \cite{aldous-pitman-entrance}, it was shown that one can associate to $\icrt$, a natural probability measure called the {\bf mass measure} satisfying $\mu(\cL(\icrt))=1$. \nomenclature[ICRTL]{$\cL(\icrt)$}{Set of leaves of $\icrt$.}

\noindent{\bf Root-to-vertex path measures: }
Now using the collection of uniform random variables above, we will define a function $\dA_{\sss(\infty)}$ on the tree as well as a collection of measures on paths emanating from the root.
Recall that the hubs in $\icrt$ have infinite degrees.  Let $(\cT_j^{\sss(i)}:j\geq 1)$ be the collection of subtrees of hub $i$ in $\icrt$ (labeled in some fashion).
For each $y\in \icrt$, let
\begin{equation}
\label{eqn:limit-da}
	\dA_{\sss(\infty)}(y)=\sum_{i\geq 1}\theta_{i}\left[\sum_{j\geq 1}U_j^{\sss(i)}\times\ind\set{y\in \cT_j^{\sss(i)}}\right].
\end{equation}
We will show in our proof that $\dA_{\sss(\infty)}(y)$ is finite for almost every realization of $\icrt$ and for $\mu$-almost every $y\in\icrt$ (see Lemma \ref{lem:tail-fexec-inf} and Theorem \ref{thm:jt-convg-E} below). For $y\in \icrt$, let $[\rho,y]$ denote the path from the root $\rho$ to $y$. For every $y$, define a probability measure on $[\rho,y]$ as
\begin{equation}
\label{eqn:right-end-prob-inft}
Q_{y}^{\sss(\infty)}(v):= \frac{\theta_i U_{j}^{\sss(i)}}{ \dA_{\sss(\infty)}(y)}, \qquad\mbox{ if } v\mbox{ is the }i\mbox{th hub and } y\in \cT_j^{\sss(i)}\mbox{ for some }j. 	
\end{equation}
\nomenclature[ICRTm]{$\dA_{\sss(\infty)}(y),Q_{y}^{\sss(\infty)}$}{Root-to-vertex weights and measures in $\icrt$ defined in \eqref{eqn:limit-da} and \eqref{eqn:right-end-prob-inft}. See analogous objects for finite trees in Section \ref{sec:gmvpa-explicit}}
Thus, this probability measure is concentrated on the hubs on the path from $y$ to the root.
\begin{rem}\label{rem:suppress}
Note that both $\dA_{\sss(\infty)}(\cdot)$ and $Q_{y}^{\sss(\infty)}(\cdot)$ depend on the realization of the pair $(\icrt, \mvU)$, but we chose to suppress them to avoid cumbersome notation.
\end{rem}

\noindent{\bf Random tree $\cR_{IJ}^{\sss(\infty)}$:} Recall the tree $r_{IJ}^{\sss(\infty)}$ above. Recall that $\eta_j$ is the vertex in the tree $\icrt$ corresponding to leaf $j+$ for $1\leq j\leq J$.  To each of these $J$ leaves, associate the value $\dA_{\sss(\infty)}(\eta_j)$, and associate the probability measure $Q_{\eta_j}^{\sss(\infty)}$ to the path $[0+, j+]$. This tree is a random element of the space $\vT_{IJ}^{*}$ (see Section \ref{sec:space-of-trees}), which we denote by $\cR_{IJ}^{\sss(\infty)}$.

\subsection{Continuum limits of components}
\label{sec:descp-limit}
The aim of this section is to give an explicit description of the limiting (random) metric spaces in Theorem \ref{thm:mc-main-1}.  We start by constructing a specific tilted version of the ICRT in Section \ref{sec:tilt-icrt}. Then in Section \ref{sec:limits-comp-descp} we describe the limits of maximal components.

\subsubsection{Tilted ICRTs and vertex identification}
\label{sec:tilt-icrt}
Let $(\Omega, \cF, \pr_{\theta})$ and $\icrt$ be as in Section \ref{sec:p-tree-ICRT-def} and let $\gamma >0$ a constant. Informally, the construction goes as follows:
We will first tilt the distribution of the original ICRT $\icrt$ using the functional
\nomenclature[tilt-icrt]{$L_{\sss(\infty)}(\icrt, \mvU)$}{Tilt functional to construct tilted ICRT.}
	\begin{equation}
	\label{eqn:ltheta-def}
		L_{\sss(\infty)}(\icrt, \mvU):= \exp\left(\gamma\int_{y\in \icrt} \dA_{\sss(\infty)}(y)\mu(dy)\right)
	\end{equation}
	to get a tilted tree $\tilicrt$.
	\nomenclature[tilICRT]{$\tilicrt$}{Tilted ICRT with distribution ${\pr}_{\theta}^\star$.}
	 We then generate a random but finite number $N_{\sss(\infty)}^\star$ of pairs of points $\set{(x_k, y_k):1\leq k\leq N_{\sss(\infty)}^\star}$.
	\nomenclature[shortcut]{$N_{\sss(\infty)}^\star$}{Number of shortcuts in $\tilicrt$. }
	 The final metric space is obtained by creating ``shortcuts" by identifying the points $x_k$ and $y_k$. Formally the construction proceeds in four steps:
\begin{enumeratea}
	\item {\bf Tilted ICRT:}  
Define ${\pr}_{\theta}^\star$ on $\Omega$ by
\[\frac{d {{\pr}}_{\theta}^\star}{d{{\pr}}_{\theta}}=\frac{\exp\left(\gamma\int_{y\in\cT^{\theta}}\dA_{(\infty)}(y)\mu(dy) \right)}{\E\left[\exp\left(\gamma\int_{x\in \cT^{\theta}} \dA_{\sss(\infty)}(x)\mu(dx) \right)\right]}. \]
The expectation in the denominator is with respect to the original measure ${\pr}_{\theta}$. In our proof we will show that this object is finite. Write $(\tilicrt, \mu^\star)$ and $\mvU^{\star}=(U_j^{(i), \star}: i,j\geq 1)$ for the tree and the mass measure on it, and the associated random variables under this change of measure.
 \item {\bf Poisson number of identification points:} Conditionally on $((\tilicrt, \mu^\star), \mvU^{\star})$, generate $N_{\sss(\infty)}^\star$ having a $\mathrm{Poisson}(\Lambda_{\sss(\infty)}^\star)$ distribution, where
\begin{align*}
\Lambda_{\sss(\infty)}^\star:= \gamma\int_{y\in \tilicrt}\dA_{\sss(\infty)}(y)\mu^\star(dy)
=\gamma\sum_{i\geq 1}\theta_{i}\left[\sum_{j\geq 1}U_j^{(i), \star}\mu^\star(\cT_j^{\sss(i), \star})\right].
\end{align*}
Here, $(\cT_j^{\sss(i), \star} : j\geq 1)$ denotes the collection of subtrees of hub $i$ in $\tilicrt$. (As mentioned before in Remark \ref{rem:suppress}, $\dA_{\sss(\infty)}(\cdot)$ depends on the realization of the ordered ICRT. $U_j^{(i), \star}$ appears in the expression above as the function $\dA_{\sss(\infty)}$ acts on $y\in\tilicrt$ for which the associated order is described by $\mvU^{\star}$.)

\item {\bf ``First'' endpoints (of shortcuts): } Conditionally on (a) and (b), sample $x_k$ from $\tilicrt$  with density proportional to $\dA_{\sss(\infty)}(x)\mu^\star(dx)$ for $1\leq k\leq N_{\sss(\infty)}^\star$.
\item {\bf ``Second'' endpoints (of shortcuts) and identification:} Having chosen $x_k$, choose $y_k$ from the path $[\rho, x_k]$ joining the root $\rho$ and $x_k$ according to the probability measure $Q_{x_k}^{(\infty)}$ as in \eqref{eqn:right-end-prob-inft} but with $U_j^{(i),\star}$ replacing $U_j^{\sss (i)}$.
(Note that $y_k$ is always a {\bf hub} on $[\rho, x_k]$.)  Identify $x_k$ and $y_k$, i.e., form the quotient space by introducing the equivalence relation $x_k\sim y_k$ for $1\leq k\leq N_{\sss(\infty)}^\star$.
\end{enumeratea}
\begin{defn}\label{def:limiting-space} Fix $\gamma\geq 0$ and $\mvtheta\in \Theta$ as in \eqref{eqn:Theta-def}. Let $\cG_{\infty}(\mvtheta,\gamma)$ be the metric measure space constructed via the four steps above equipped with the measure inherited from the mass measure on $\tilicrt$.	
\end{defn}
In our proofs, we will always think of the leaf end (of a shortcut or a surplus edge) as the first endpoint, and the second endpoint will be selected from the skeleton.

\subsubsection{Limits of the components}
\label{sec:limits-comp-descp}
Fix $\lambda \in \bR$ and $\vc\in l_0 $ as in \eqref{eqn:lthree-lo} and consider the setting of Theorem \ref{thm:mc-main-1}. We will need 2 main objects:
\begin{enumeratea}
	\item The process $\tilde{V}^{\vc}_{\lambda}(\cdot)$ in \eqref{eqn:vc-reflection}. Recall that the excursions of this process from zero could be arranged in increasing order of lengths as $\sZ(\lambda)$. Let $\Xi^{\sss(i)} = (c_j: \xi_j \in \cZ_i)$ denote the point process of jumps of the process $\tilde{V}^{\vc}_{\lambda}(\cdot)$ corresponding to the excursion $\cZ_i(\lambda)$. Abusing notation we will write  $\Xi^{\sss(i)} = (c_j: j \in \cZ_i)$.
	\item The actual lengths of these excursions $(Z_i(\lambda)\colon i\geq 1)$ as in \eqref{eqn:vz-excursions-def}.
\end{enumeratea}
From these objects, for each fixed $i\geq 1$, define the random variable $\barg^{\sss(i)}$ and the point process $\mvtheta^{\sss(i)} = (\theta_j^{\sss(i)}:j\in \cZ_i(\lambda))$ as
\begin{equation}
\label{eqn:barg-mvtheta-def-comps}
	\barg^{\sss(i)}:= Z_i(\lambda)\sqrt{\sum_{v\in \cZ_i(\lambda)} c_v^2}, \qquad \mvtheta^{\sss(i)}:=  \left(\frac{c_j}{\sqrt{\sum_{v\in \cZ_i(\lambda)} c_v^2}}: j\in \cZ_i(\lambda)\right).
\end{equation}
Our proof (see Proposition \ref{prop:weights-good}) will imply that $\mvtheta^{\sss(i)} \in \Theta$ as in \eqref{eqn:Theta-def} a.s.
Define
\begin{align}\label{eqn:def-Gamma-i}
\Gamma_i(\lambda):=Z_i(\lambda)\left(\sum_{v\in\cZ_i(\lambda)}c_v^2\right)^{-1/2},
\end{align}
and generate the random metric measure spaces
\begin{equation}
\label{eqn:mi-c-def}
	M_i^{\vc}(\lambda):= \Gamma_i(\lambda)\cdot\cG_{\infty}(\mvtheta^{\sss(i)}, \barg^{\sss(i)}),
\end{equation}
where $\cG_{\infty}(\mvtheta,\barg)$ is as described in Section \ref{sec:tilt-icrt} and the metric spaces are {\it conditionally independent} across $i$ given the driving parameters in \eqref{eqn:barg-mvtheta-def-comps}. Let $\vM_{\infty}^{\vc}(\lambda) = (M_i^{\vc}(\lambda)\colon i\geq 1)$. Then this is the limiting collection of metric spaces in Theorem \ref{thm:mc-main-1}.

To describe the sequence of spaces $\vM_{\infty}^{\nr}(\lambda)$ appearing in Theorem \ref{thm:rank-one}, define
\begin{align}\label{eqn:def-c-nr}
\vc^{\nr}:=(c_j^{\nr}: j\geq 1),\qquad\text{ where }\qquad c_j^{\nr}= \frac{1}{\bE W}\left(\frac{c_F}{j}\right)^{1/(\tau-1)},
\end{align}
\begin{align}\label{eqn:t-nr}
\zeta:=-\left(\frac{c_F^{2/(\tau-1)}}{\E W}\right)\sum_{i=1}^{\infty}\left[\int_{i-1}^i\frac{du}{u^{2/(\tau-1)}}-\frac{1}{i^{2/(\tau-1)}}\right], 
\qquad \text{ and }\qquad  t^{\nr}_{\lambda}:=\frac{(\lambda+\zeta)}{\E W}.
\end{align}
Here $W$ is a random variable with distribution $F$ as in \eqref{eqn:tau-c-def}. Then
\begin{align}\label{eqn:def-M-nr}
\vM_{\infty}^{\nr}(\lambda)=\frac{1}{\E W}\cdot\vM_{\infty}^{\vc^{\nr}}\left(t^{\nr}_{\lambda}\right).
\end{align}


\section{Discussion}\label{sec:discussion}
We describe the two major motivations for developing the general theory of this paper in Sections \ref{sec:disc-dom} and \ref{sec:disc-mst}. In Sections \ref{sec:disc-icrt} and \ref{sec:disc-overview}, we include a brief discussion about ICRTs as well as give an overview of the order in which the proofs are carried out.

\subsection{Universality and domains of attraction of critical random graph models}
\label{sec:disc-dom}
One natural question the reader might ask at this point is why the general theory in Section \ref{sec:mc-results}, why not just stick to the rank-one random graph model as in Section \ref{sec:rank-one}. As we have described in the introduction, the aim of this paper is the development of general theory applicable to a wide array of models. What does one mean by this? It turns out that many different random graph models can be constructed in a dynamic fashion as a graph-valued process $\set{\cG_n(t): t\geq 0}$ where edges are added as time advances thus resulting in mergers of components as $t\uparrow t_c$. In this construction, there is a critical time $t_c$ (model-dependent) such that the giant component emerges after time $t_c$.

Now for most random graph models (including the configuration model) the dynamics of mergers of components starting at time zero {\bf do not} look like the multiplicative coalescent. However if one were to zoom in at the critical time $t_c$, for many models, there exists $\eps_n\downarrow 0$ such that if one were to look at the interval $[t_c-\eps_n, t_c+\eps_n]$, then mergers of components {\bf can} be approximated by the multiplicative coalescent. Here $t_c -\eps_n$ often corresponds to the \emph{barely subcritical} regime of the random graph.  Thus if one had good control over component functionals at the barely subcritical time $t_c-\eps_n$ and in particular if one was able to show that component sizes appropriately normalized satisfied Assumption \ref{ass:mc-assumptions}, then one can use Theorem \ref{thm:mc-main-1} to derive convergence at the critical time $t_c$ of the maximal components. Note that one does not expect component sizes at time $t_c-\eps_n$ to satisfy assumptions of the Norros-Reittu model in \eqref{eqn:degree-dist}. Rather in most cases, at time $t_c-\eps_n$,  the expected size of the component of a randomly selected vertex $V_n$ would scale like $n^{\delta_1}$ while the maximal component would scale like $n^{\delta_2}$ (ignoring logarithmic corrections) where $\delta_1 < \delta_2$ are related to various scaling exponents of the system. In work in progress \cite{SB-SD-vdH-SS},  Theorem \ref{thm:mc-main-2-tau} coupled with delicate estimates of various scaling exponents for the configuration model in the barely subcritical regime, proves analogous results for the configuration model with degree exponent $\tau \in (3,4)$. Sizes of maximal components in the critical regime including the heavy-tailed regime for this model was previously analyzed in \cite{joseph2014component}.  Further as was done in \cite{SBSSXW-universal}, where a number of sufficient conditions for the domain of attraction of the critical \erdos scaling limits were derived, we hope to derive similar general conditions for a random graph model to belong to the same domain of attraction as the rank-one model with $\tau\in(3,4)$, established in this paper.


\begin{figure}[htbp]
	\centering
		\includegraphics[trim=0cm 3cm 0cm 5cm, clip=true, scale=.4]{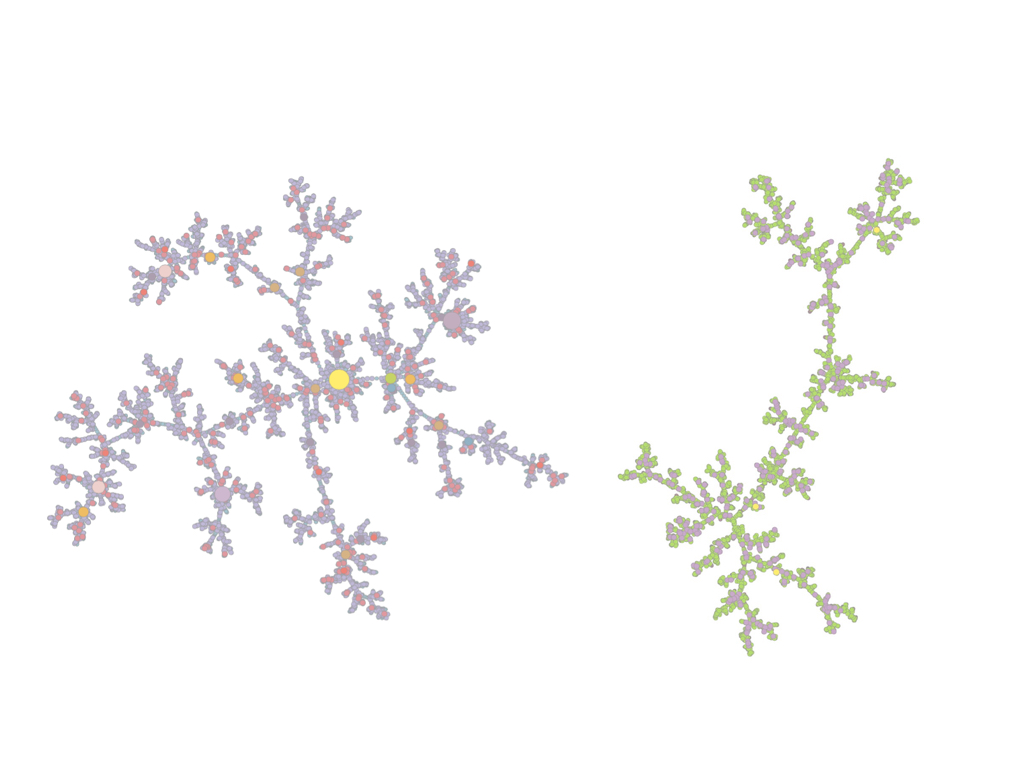}
	\caption{On the left, an approximation of an ICRT (using $\vp$-trees on approximately $20000$ vertices) corresponding to $\theta_i\propto i^{-1/(\tau-2)}$ where $\tau=3.01$. The reason behind this choice of $\theta_i$ is explained in Section \ref{sec:disc-open}. On the right, an approximation of a Brownian CRT (using a uniform random tree on the same number of vertices). Vertex sizes are proportional to the degree of the vertex.}
	\label{fig:crt-icrt}
\end{figure}

\subsection{Minimal spanning tree on inhomogeneous random graphs}
\label{sec:disc-mst}
As described in the introduction, a second major motivation for the technical analysis in this paper is the minimal spanning tree. To fix ideas, consider the Norros-Reittu model in the supercritical regime (the parameter in \eqref{eqn:nu-def} $\nu > 1$). To each edge attach a random edge weight i.i.d.\ across edges, assumed to be derived from a continuous distribution. Consider the minimal spanning tree (MST) of the giant component. A large amount of simulation-based evidence from statistical physics \cite{BraBulCohHavSta03,wu2006transport,braunstein2007optimal,chen2006universal} suggests that when the degree exponent $\tau\in (3,4)$ then the distances in this object scale like $n^{(\tau-3)/(\tau-1)}$, the same distance scaling shown in this paper for the maximal components in the \emph{critical} regime (Theorem \ref{thm:rank-one}).

This is not a coincidence. As has been shown in a series of fundamental papers \cite{BBG-12,BBG-limit-prop-11,AddBroGolMie13} for the complete graph and the supercritical \erdos random graph, a major ingredient in the analysis of the MST problem is the scaling of maximal components in the critical regime which then provides crucial input for the scaling limit of the MST. Till date we have no rigorous results for the scaling of the MST on any ``inhomogeneous'' random graph model.  This paper provides the first step in answering this question in the heavy-tailed regime. Further this program should enable one to analyze the MST for random graph models other than the rank-one model which belong to the same ``domain of attraction'' in the critical regime.

\subsection{Inhomogeneous continuum random trees}
\label{sec:disc-icrt}
As evident from Section \ref{sec:p-tree-ICRT-def},  ICRTs play a major role in the description of our limiting objects. Despite a lot of work on these objects in the last decade \cite{AMP,aldous-pitman-entrance,pitman-camarri}, a number of questions regarding these continuum objects are still open, ranging from sufficient conditions for compactness to the dependence of the fractal properties of this object on the driving parameter $\mvtheta$. Our proof shows that in some special cases,  ICRTs are compact metric spaces when $\mvtheta$ is sampled according to an appropriate size-biased distribution. This can be seen as an {\em annealed} result on compactness of the ICRT. Whether compactness is true for non-random sequences $\mvtheta\in\Theta$ has been open problem for more than a decade \cite{AMP}. Similar questions hold for its fractal dimensions. See Section \ref{sec:disc-open} for a more detailed account of these problems.


\subsection{Overview of the proof}\label{sec:disc-overview}
In Section \ref{sec:proof1}, we study the random graph $\cG_n(\vx,t)$ as in Definition \ref{def:cgvx-q}. We start with the simple observation that conditional on the vertex set of components of $\cG_n(\vx,t)$, a fixed component $\cC$ has the same distribution as $\cG_n(\vx,t)$ \emph{conditional} on being connected.  This section studies asymptotics for such distributions assuming specific regularity properties of
vertex weights in the component in the large network limit, showing Gromov-weak convergence of the associated graph under proper normalization of edge lengths and vertex weights. Section \ref{sec:proof-convg-gromov} uses the size-biased exploration of the process $\cG_n(\vx,t)$ \cite{aldous-crit} to show that maximal connected components satisfy the hypothesis required in Section \ref{sec:proof1}. Section \ref{sec:proof-convg-dghp} studies the special entrance boundary in \eqref{eqn:c-tau-lamb-def} proving both compactness of the limiting objects as well as strengthening the convergence in the Gromov-weak topology to convergence in $d_{\GHP}$. In Section \ref{sec:fractal-proofs}, we derive the box-counting or Minkowski dimension. In Section \ref{sec:disc-open}, we conclude by describing a number of open problems.

\section{Proofs: Asymptotics conditional on being connected }
\label{sec:proof1}
The aim of this Section is to study large connected components of $\cG_n(\vx,t)$ assuming vertex weights satisfy a few regularity properties.

\subsection{Tilted $\vp$-trees and connected components of $\cG(\vx,t)$}

Recall the random graph $\cG(\vx,t)$ from Definition \ref{def:cgvx-q}. Here for any $t\geq 0$, $(\cC_i(t)\colon i\geq 1)$ \ch{denotes} the components in decreasing order of their mass sizes. In this section we will describe results from \cite{SBSSXW14} which gave a method of constructing {\bf connected} components of $\cG(\vx,t)$ conditional on the vertices of the components. This construction involved tilted versions of $\vp$-trees introduced in Section \ref{sec:p-tree-ICRT-def}. Since these trees are parametrized via a driving \ch{probability mass function (pmf)} $\vp$, it will be easy to parametrize various random graph constructions in terms of pmfs as opposed to vertex weights $\vx$. Proposition \ref{prop:generate-nr-given-partition} will relate vertex weights to pmfs.

Fix $n\geq 1$ and $\cV\subset [n]$ and write $\bG_{\cV}^{\con}$ for the space of all simple connected graphs with vertex set $\cV$.
For fixed $a > 0$, and probability mass function $\vp = (p_v: v \in \cV)$, define probability distributions $\pr_{\con}(\cdot; \vp, a, \cV)$ on $\bG_{\cV}^{\con}$ as follows: Define for $i,j \in \cV$,
\begin{equation}
\label{eqn:qij-def-vp}
	q_{ij}:= 1-\exp(-a p_i p_j).
\end{equation}
Then
\begin{equation}
	\label{eqn:pr-con-vp-a-cV-def}
	\pr_{\con}(G; \vp, a, \cV): = \frac{1}{Z(\vp,a)} \prod_{(i,j)\in E(G)} q_{ij} \prod_{(i,j)\notin E(G)} (1-q_{ij}), \mbox{ for } G \in \bG_{\cV}^{\con},
\end{equation}
where $Z(\vp,a)$ is the normalizing constant
$$Z(\vp,a) := \sum_{G \in \bG_{\cV}^{\con}}{\prod_{(i,j)\in E(G)} q_{ij} \prod_{(i,j)\notin E(G)}(1-q_{ij})}.$$

Now let $\cV^{\sss(i)} := V(\cC_i(t))$ be the vertex set of $\cC_i(t)$ for $i \geq 1$ and note that $\set{\cV^{\sss(i)}\colon i\geq 1}$ denotes a random finite partition of the full vertex set $[n]$. The following result is obvious from the construction of $\cG(\vx,t)$:  

\begin{prop}[{\protect{\cite[Proposition 6.1]{SBSSXW14}}}]
	\label{prop:generate-nr-given-partition}
	Conditional on the partition $\set{\cV^{\sss(i)}\colon i\geq 1}$ define
	\begin{equation*}
		\vp_n^{\sss(i)} := \left( \frac{x_v}{\sum_{v \in \cV^{\sss(i)}}x_v } : v \in \cV^{\sss(i)} \right), \qquad  a_n^{\sss(i)}:= t\left(\sum_{v\in \cV_{\sss(i)}} x_v\right)^2, \qquad i\geq 1.
	\end{equation*}
	For each fixed $i \geq 1$, let $G_i \in  \bG_{\cV^{\sss(i)}}^{\con}$ be a connected simple graph with vertex set $\cV^{\sss(i)}$. Then
	\begin{equation*}
		\pr\left(\cC_i(t) = G_i, \;\; \forall i \geq 1\ \big|\ \set{\cV^{\sss(i)}\colon i\geq 1} \right) = \prod_{i\geq 1} \pr_{\con}( G_i; \vp_n^{\sss(i)}, a_n^{\sss(i)}, \cV^{\sss(i)}).
	\end{equation*}
\end{prop}
Thus the random graph $\cG(\vx,t)$ can be generated in two stages:
\begin{enumeratei}
	\item {\bf Stage I:} Generate the partition of the vertices
into different components, i.e., generate $\set{\cV^{\sss(i)}\colon i\geq 1}$.
\item {\bf Stage II:} Conditional on the partition, generate the internal structure of each component following the law of $\pr_{\con}(\cdot ; \vp^{\sss(i)}, a^{\sss(i)}, \cV^{\sss(i)})$, independently across different components.
\end{enumeratei}
Let us now describe an algorithm to generate such {\bf connected} components using distribution \eqref{eqn:pr-con-vp-a-cV-def}. To ease notation, let $\cV = [m]$ for some $m\geq 1$ and fix a probability mass function $\vp$ on $[m]$ and a constant $a>0$ and write $\pr_{\con}(\cdot):= \pr_{\con}(\cdot;\vp,a,[m])$ on $\bG_m^{\con}:= \bG_{[m]}^{\con}$. We will first need to set up some notation before describing this result.

\medskip
\noindent{\bf Depth-first exploration of ordered trees.} Recall that we used $\bT_m^{\ord}$ for the space of ordered (or planar) trees with vertex set $[m]$. Given a tree $\vt \in \bT_m^{\ord}$, one can use the associated order to explore the tree in a depth-first manner. More precisely we start with $v(1)$ being the root of $\vt$. At each stage $1\leq i\leq m$, we will keep track of three types of vertices: the set of \emph{active} vertices--$\cA(i)$, the set of \emph{explored} vertices--$\cO(i)$, and the set of \emph{unexplored} vertices--$\cU(i)$. The set of active vertices will in fact be viewed as a vertical stack (not just a set) with $\cA(i)$ representing the state of this stack at the end of step $\cA(i)$. Initialize the process with $\cA(1) = \set{v(1)}$ (the root of $\vt$), $\cO(1) = \emptyset$ and $\cU(1) = [m]\setminus \set{v(1)}$. At step $i\geq 1$, we let
\begin{enumeratei}
	\item $v(i)$ denote the vertex at the top of the stack $\cA(i)$ and let $\cD(i)\subset \cU(i)$ denote the set of children of $v(i)$. Delete $v(i)$ from $\cA(i)$ and arrange the vertices of $\cD(i)$ from oldest to youngest at the \emph{top} of the stack to \ch{form} $\cA(i+1)$;
	\item $\cO(i+1) = \cO(i) \cup \set{v(i)}$;
	\item $\cU(i+1) = \cU(i)\setminus \cD(i)$.
\end{enumeratei}

 Write $\sP(\vt)$ for set of pairs of vertices $\set{u,v}$ such that $u,v\in \cA(i)$ for some $1\leq i\leq m$; namely both vertices are active but have not yet been explored. Using terminology from \cite{BBG-12}, call this collection \emph{the set of permitted edges}. \nomenclature[perm]{$\sP(\vt)$}{Set of permitted edges in a tree $\vt$.} Thus,
 \begin{equation}
 \label{eqn:permitted-set}
 	\sP(\vt):= \set{(v(i),u)\ \big|\ 2\leq i\leq m,\ u\in \cA(i-1)\setminus \set{v(i)}}.
 \end{equation}
 Write $E(\vt)$ for the edge set of $\vt$.  Now define the function  $L : \bT_m^{\ord} \to \bR_+$ by
\begin{equation}
\label{eqn:ltpi-def}
	\displaystyle L(\vt)=\displaystyle L_{\sss (m)}(\vt):= \prod_{(k,\ell)\in E(\vt)} \left[\frac{\exp(a p_k p_{\ell})- 1}{ap_k p_{\ell}} \right] \exp\left(\sum_{(k,\ell) \in \sP(\vt)} a p_k p_{\ell}\right), \qquad \vt \in \bT_m^{\ord}.
\end{equation}
Recall the (ordered) $\vp$-tree distribution from \eqref{eqn:ordered-p-tree-def}. Using $L(\cdot)$ to tilt this distribution results in the distribution
\begin{equation}
	\label{eqn:tilt-ord-dist-def}
	\pr_{\ord}^\star( \vt) := \pr_{\ord}(\vt) \cdot \frac{L(\vt)}{\E_{\ord}[ L(\cT^{\vp}_m)]}, \qquad \vt \in \bT_m^{\ord}.
\end{equation}
\nomenclature[til-pt]{$L(\cdot),\pr_{\ord}^\star$}{Tilt functional and associated tilted $\vp$-tree distribution}
For future reference we fix notation for the various objects required in the proof below.
\nomenclature[gmcon]{$\tilde{\cG}_m(\vp, a)$}{Random graph with distribution $\pr_{\con}(\cdot,\vp,a,[m])$ defined in \eqref{eqn:pr-con-vp-a-cV-def}. See Section \ref{sec:gmvpa-explicit} for the modified random graph ${\cG}_m^{\md}(\vp,a)$.}
\begin{defn}
	\label{def:con-graph-tilde-p}
	Fix $m\geq 1$, $a> 0$, and a probability mass function $\vp$ on $[m]$. We will write $\tilde{\cG}_m(\vp, a)$ to denote a random graph with distribution $\pr_{\con}(\cdot,\vp,a,[m])$. ${\cT}^{\vp, \star}_m$ will denote a random planar tree with the tilted $\vp$-tree distribution \eqref{eqn:tilt-ord-dist-def}, and ${\cT}^{\vp}_m$ will denote a random tree with the original $\vp$ tree distribution \eqref{eqn:ordered-p-tree-def}.
\end{defn}
\nomenclature[ptree]{${\cT}^{\vp}_m,{\cT}^{\vp, \star}_m$}{Random $\vp$-tree, respectively tilted $\vp$-tree using $L(\cdot)$}
\begin{prop}[{\cite[Proposition 7.4]{SBSSXW14}}]
	\label{prop:SBSSXW}
	Fix $m\geq 1$, a probability mass function $\vp$ on $[m]$, and $a>0$. Consider a random connected graph on $[m]$ constructed as follows:
	\begin{enumeratea}
		\item First generate a rooted planar random tree ${\cT}^{\vp, \star}_m$ with distribution $\tilde \pr_{\ord}(\cdot)$ as in \eqref{eqn:tilt-ord-dist-def}.
		\item Let $\sP({\cT}^{\vp, \star}_m)$ denote the permitted edge set of this random tree. Add each such edge $\set{u,v} \in \sP({\cT}^{\vp, \star}_m)$ with probability $q_{uv}$ as in \eqref{eqn:qij-def-vp}, independent across permitted edges.
	\end{enumeratea}
	Then, the resulting random graph has distribution $\pr_{\con}$ on $\bG_m^{\con}$, i.e, has the same distribution as $\tilde{\cG}_m(\vp, a)$.
\end{prop}

\subsection{Convergence of connected components under weight assumptions}
The aim of this section is to prove Gromov-weak convergence for the connected graph $\tilde\cG_m(\vp,a)$ under regularity conditions on $a$ and $\vp$ as $m\to\infty$. We will assume that we have ordered the index set $[m]$ so that $p_1\geq p_2\geq \cdots \geq p_m >0$. Let
\begin{equation}
\label{eqn:sigma-vp-def}
	\sigma(\vp):= \sqrt{\sum_i p_i^2}.
\end{equation}

\begin{ass}
	\label{ass:vp-a}
	 As $m\to\infty$, the following hold:
	\begin{enumeratei}
		\item  $\sigma(\vp)\to 0$ and further for each fixed $i\geq 1$, $p_i/\sigma(\vp)\to \theta_i$ where $\mvtheta:= (\theta_1, \theta_2, \ldots)$ is an element of $~\Theta$ as in \eqref{eqn:Theta-def}.
		\item There is a constant $\gamma >0$ such that $a\sigma(\vp)\to \gamma$.
	\end{enumeratei}	
\end{ass}

The following theorem is the main result of this section.
\begin{thm}\label{thm:connected-convg}
	Consider the connected random graph $\tilde{\cG}_m(\vp,a)$ viewed as a metric measure space via the graph distance where each vertex $v$ is assigned measure $p_v$.  Under Assumption \ref{ass:vp-a},
	\[\sigma(\vp) \tilde{\cG}_m(\vp,a) \convd \cG_{\infty}(\mvtheta,\gamma), \]
	where $\cG_{\infty}(\mvtheta,\gamma)$ is the random metric space defined in Definition \ref{def:limiting-space} and convergence is in the Gromov-weak topology on metric spaces.
\end{thm}
The rest of this section proves this result. We will throughout assume that $\tilde{\cG}_m(\vp,a)$ has been
 constructed using Proposition \ref{prop:SBSSXW}.
 \subsubsection{Two constructions of $\vp$-trees:  Exploration process and the birthday construction}
 \label{sec:ptree-const}
 We start by describing an explicit construction of the (untilted) $\vp$-tree $\cT_m^{\vp}$ first developed in \cite{AMP}. At the end of this section we describe a second construction used later in the paper.

\begin{figure}[htbp]
	\centering
	\begin{tikzpicture}[scale=1]
	
		\draw[blue, very thick, ->] (0,-3) -- (0,2);
		\draw[blue, very thick, ->] (0,0) -- (10,0);
		\draw[black, thick] (0,0) -- (.5,-.5);
		\draw[black, thick] (.5,0) -- (1,-.5);
		\draw[black, thick] (1,.5) -- (2.5,-1);
		\draw[black, thick] (2.5,0) -- (4,-1.5);
		\draw[black, thick] (4,-1) -- (5,-2);
		\draw[black, thick] (5,-1) -- (7,-3);
		\draw[black, thick] (7,-1) -- (8,-2);
		\draw[black, thick] (8,1) -- (9.5,-.5);
		\draw[black, thick] (9.5,.5) -- (10,0);


		\draw[very thick,blue, |->] (0,-4.5) -- (10,-4.5);
	    \node at (.5,-4.5) [circle,draw=red!50,fill=red!50, label = above: $u_8$] {};
	 \node at (1,-4.5) [circle,draw=red!50,fill=red!50, label = above: $u_6$] {};
	 \node at (2.5,-4.5) [circle,draw=red!50,fill=red!50, label = above: $u_3$] {};
	 \node at (5,-4.5) [circle,draw=red!50,fill=red!50, label = above: $u_1$] {};
	 \node at (7,-4.5) [circle,draw=red!50,fill=red!50, label = above: $u_2$] {};
	 \node at (8,-4.5) [circle,draw=red!50,fill=red!50, label = above: $u_4$] {};
	 \node at (4,-4.5) [circle,draw=red!50,fill=red!50, label = above: $u_5$] {};
	 \node at (9.5,-4.5) [circle,draw=red!50,fill=red!50, label = above: $u_7$] {};

	 \path[thick, red, ->] (7.05,-5) edge            node[below]  {$p_2$}     (8.95,-5);
	 \path[thick, red, ->] (9.05,-5) edge                 (10,-5);
	 \path[thick, red, ->] (0,-5) edge            node[below]  {$p_4$}     (1.95,-5);
	 \path[thick, red, ->] (2.05,-5) edge            node[below]  {$p_7$}     (2.95,-5);
	 \path[thick, red, ->] (3.05,-5) edge            node[below]  {$p_3$}     (4.25,-5);
	 \path[thick, red, ->] (4.35,-5) edge            node[below]  {$p_5$}     (4.65,-5);
	 \path[thick, red, ->] (4.75,-5) edge            node[below]  {$p_8$}     (5.35,-5);
	 \path[thick, red, ->] (5.45,-5) edge            node[below]  {$p_1$}     (6.05,-5);
	 \path[thick, red, ->] (6.15,-5) edge            node[below]  {$p_6$}     (6.95,-5);

	 \path[thick, blue, -] (-.05, 1.8) edge            node[left]  { ${\scriptstyle .2}$}     (.05,1.8);
	 \path[thick, blue, -] (-.05, -1.8) edge            node[left]  {${\scriptstyle -.2}$}     (.05,-1.8);

	\end{tikzpicture}
	
	\vspace{.4in}
	
	\begin{tikzpicture}[scale=.3]

	    \node (1) [circle,draw=red!50,fill=red!50, label = below: Root] at (0,0) {$2$};
	    \node (2) [circle,draw=red!50,fill=red!50] at (0,5) {$4$};
	    \node (3) [circle,draw=red!50,fill=red!50] at (-4,9) {$7$};
	    \node (4) [circle,draw=red!50,fill=red!50] at (-7,12) {$3$};
	    \node (5) [circle,draw=red!50,fill=red!50] at (-10,15) {$5$};
	    \node (6) [circle,draw=red!50,fill=red!50] at (4,9) {$6$};
	    \node (8) [circle,draw=red!50,fill=red!50] at (0,9.5) {$8$};
	    \node (7) [circle,draw=red!50,fill=red!50] at (0,13.5) {$1$};

        \draw[blue, very thick] (1)-- (2) -- (3) -- (4) -- (5);
        \draw[blue, very thick] (2) -- (6);
        \draw[blue, very thick] (2) -- (8) -- (7);
	
	\end{tikzpicture}

	\caption{The function $F^{\vp}$ and the corresponding tree $\psi_{\vp}$. }
	\label{fig:Fvp-psivp}
\end{figure}
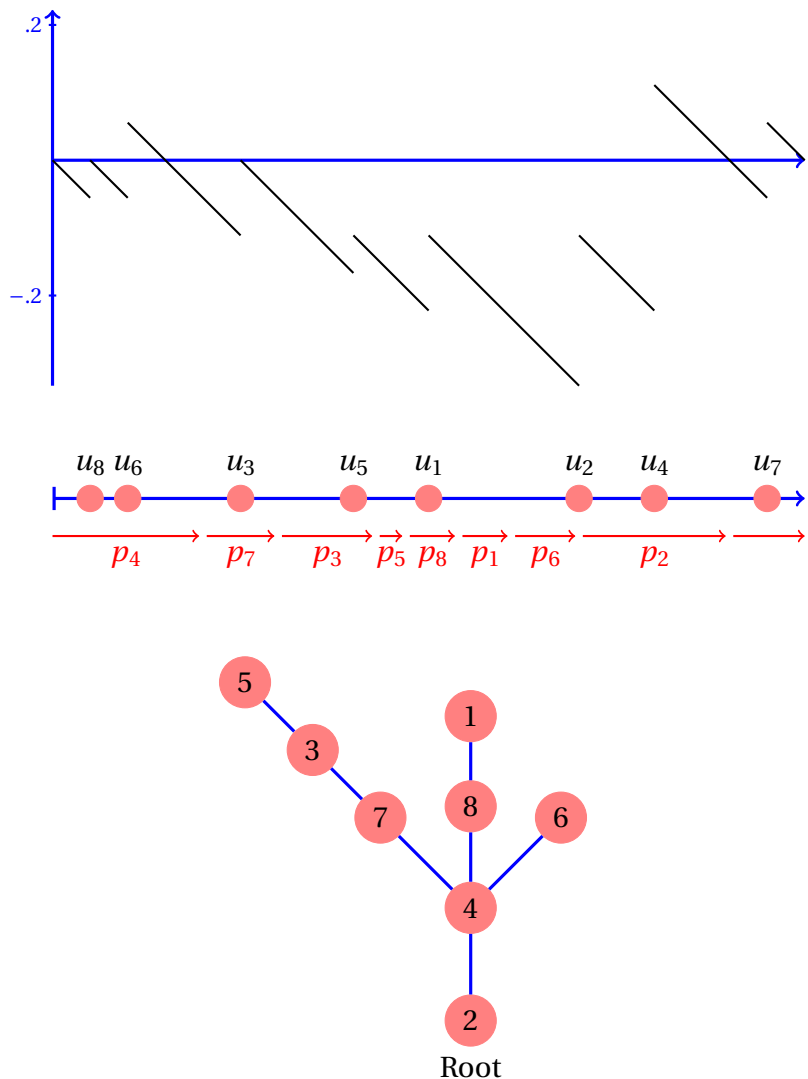

\noindent{\bf Exploration process construction: }
 The first construction is initiated by setting up a map $\psi_{\vp}:[0,1]^m\to \bT^{\ord}$ as follows. Let $\vu:=(u_v:v\in [m])$ be a collection of distinct points in $(0,1)$. Define
\begin{equation}
\label{eqn:fp-defn}
	F^{\vp}(s) := - s + \sum_{v=1}^m p_v \ind\set{u_v\leq s}, \qquad s\in [0,1].
\end{equation}
Assume that there exists a unique point $v^* \in [m] $ such that $F^{\vp}(u_{v^*}-) = \ch{\min_{s\in [0,1]}} F^{\vp}(s)$. Set $v^*$ to be the root of the tree $\psi_{\vp}(\vu)$. Define $y_i := u_i - u_{v^*}$ \ch{$\mbox{ mod } 1$} for $i \in [m]$, and
\begin{equation*}
\label{eqn:fexc-def}
	F^{\exec,\vp}(s):= F^{\vp}( u_{v^*} + s \mbox{ mod } 1) - F^{\vp}(u_{v*}-), \qquad 0\leq s < 1.
\end{equation*}
{Then $F^{\exec,\vp}(1-) = 0$ and $F^{\exec,\vp}(s) > 0$ for $s \in [0,1)$. Extend the definition of $F^{\exec,\vp}$ to $s \in [0,1]$ by define $F^{\exec,\vp}(1) = 0$.} We use $F^{\exec,\vp}$ to construct a depth-first-search of an ordered tree whose exploration in this depth-first manner is encoded by the function $F^{\exec,\vp}$. This in turn defines the tree $\psi_{\vp}(\vu)$. As before, in this construction we carry along a set of explored vertices $\cO(i)$, active vertices $\cA(i)$ and unexplored vertices $\cU(i) = [m]\setminus (\cA(i)\cup \cO(i))$, for $0\leq i \leq m$. We view $\cA(i)$ as the state of a vertical stack $\cA$ after the $i$th step in the depth-first-search. Initialize with $\cO(0) = \emptyset$, $\cA(0) = \set{v^*}$, $\cU(0) = [m] \setminus \set{v(1)}$, and define $y^*(0) = 0$. At step $i \in [m]$, let $v(i)$ be the value that is on the top of the stack $\cA(i-1)$ and define $y^*(i) := y^*(i-1)+p_{v(i)}$. Define $\cD(i) := \set{ \ch{j \in [m]} : y^*(i-1) < \ch{y_j} < y^*(i) }$. Suppose $\cD(i) = \set{ u(j) : \ch{1\leq j \leq k}}$ where we have ordered these vertices in the sequence that they are found in this interval, i.e.,
\[y^*(i-1) < y_{u(1)} <... < y_{u(k)} < y^*(i).\]
Update the stack $\cA$ as follows:
\begin{enumeratei}
	\item Delete $v(i)$ from $\cA$.
	\item Push $u(j)$, $1\leq j\leq k$, to the top of $\cA$ sequentially (so that $u(k)$ will be on the {\bf top} of the stack at the end).
\end{enumeratei}
Let $\cA(i)$ be the state of the stack after the above operations.
Update $\cO(i) := \cO(i-1) \cup \set{v(i)}$ and $\cU := \cU(i-1)\setminus \cD(i) $. See Figure \ref{fig:Fvp-psivp} for a pictorial description of this construction.

The tree $\psi_{\vp}(\vu) \in \bT_m^{\ord}$ is constructed by putting the edges $\set{ (v(i),v): i \in [m], v \in \cD(i)}$ and using the order prescribed in the above exploration to make the tree an ordered tree. The fact that this procedure actually produces a tree is proved in \cite{AMP}.  
\begin{lemma}[{\cite[Section 3.2]{AMP}}]
	\label{lem:AMP-Fex-ptree}
	Consider the map $\psi_{\vp}$. Let $\vX:=(X_v: v\in [m])$ be i.i.d.\ random variables distributed uniformly on $(0,1)$. Then the random tree $\psi_{\vp}(\vX)$ has distribution \eqref{eqn:ordered-p-tree-def}, i.e., $\psi_{\vp}(\vX) \stackrel{d}{=} \cT^{\vp}$.
\end{lemma}

For future reference, coupled with the above construction, define $\sS(i):=\cA(i-1)\setminus \set{v(i)}$ for $i\in [m]$. Define the function $A_m(\cdot)$ on $[0,1]$ via
\begin{equation}
\label{eqn:ant-def}
	A_m(u):= \sum_{v\in \cS(i)}  p_v, \qquad \mbox{ for } u \in (y^*(i-1), y^*(i)], i \in [m].
\end{equation}
\nomenclature[Fp-Am]{$F^{\vp}(\cdot),A_m(\cdot)$}{Functions in the depth first construction of a $\vp$-tree.}
Further let $\bar A_m(u) := a A_m (u)$, $u \in [0,1]$, where $a$ is the scaling constant in \eqref{eqn:qij-def-vp}.
\ \\
\noindent{\bf Birthday construction:} We now describe a second construction of $\vp$-trees, first formulated in \cite{pitman-camarri}. We urge the reader to skim this portion and return to it once she has reached Section \ref{sec:completing-proof-thm-con}. Let $\vY:=(Y_0, Y_1, \ldots)$ be an infinite sequence of \emph{i.i.d.} random variables with distribution $\vp$. Let $R_0=0$ and for $l\geq 1$, let $R_l$ denote the $l$-th repeat time, i.e.,
\begin{align*}
R_l=\min\bigg\{k>R_{l-1}: Y_k\in\{Y_0,\hdots,Y_{k-1}\}\bigg\}.
\end{align*}

Now consider the directed graph formed via the edges
\[\cT(\vY):= \set{(Y_{j-1}, Y_j): Y_j\notin \set{Y_0, \ldots, Y_{j-1}}, j\geq 1}.\]
It is easy to check that this gives a tree which we view as rooted at $Y_0$. Intuitively the process of constructing a tree is as follows: the tree ``grows'' via the addition of new vertices sampled using $\vp$ till it stumbles across a ``repeat'' (a vertex already found) when it goes back to the first occurrence of this ``repeat'' and starts growing from that position. The following striking result was shown in \cite{pitman-camarri}.

\begin{thm}[{\cite[Lemma 1 and Theorem 2]{pitman-camarri}}]
	\label{thm:pit-cam-birthday}
	The random tree $\cT(\vY)$ viewed as an object in $\bT_m$ is distributed as a $\vp$-tree with distribution \eqref{eqn:p-tree-def} {\bf independently} of $Y_{R_1-1}, Y_{R_2-1}, \ldots$ which are i.i.d.\ with distribution $\vp$.
\end{thm}

 \begin{rem}
 	\label{rem:p-tree-dist} The independence between the sequence $Y_{R_1-1}, Y_{R_2-1}, \ldots$ and the constructed $ \vp$ tree $\cT(\vY)$ is truly remarkable. In particular, suppose $\cS$ is a $\vp$-tree with distribution as in \eqref{eqn:p-tree-def} and for fixed $r\geq 1$, let  $\tilde{Y_1}, \tilde{Y_2}, \ldots \tilde{Y_r} $ be i.i.d. with distribution $\vp$. Write $\cS_r\subset \cS$ for the tree spanned by these vertices and the root. Let $\cT_r^{\cB}\subset \cT(\vY)$ denote the subtree with vertex set $\set{Y_0, Y_1, \ldots, Y_{R_r-1}}$, namely the tree constructed in the first $R_r$ steps. Here $\cB$ is a mnemonic for ``birthday tree'' and also to distinguish this construction from a generic random tree model with $r$ vertices.   Then the above result (formalized as \cite[Corollary 3]{pitman-camarri}) implies that these can be jointly constructed as
	\begin{equation}
	\label{eqn:p-tree-joint-sample}
		(\tilde{Y_1}, \tilde{Y_2},\ldots, \tilde{Y_r}; \cS_r)\stackrel{d}{=} (Y_{R_1-1}, Y_{R_2-1}, \ldots Y_{R_r -1}; \cT_r^{\cB}).
	\end{equation}
	We use this fact often in Section \ref{sec:completing-proof-thm-con}.
 \end{rem}

\subsection{Uniform integrability of the tilt}

The first use of the above construction of the $\vp$-tree is to prove the following:

\begin{prop}
	\label{prop:tilt-tightness}
	Fix $s\geq 1$ and consider the tilt $L(\cdot)$ as in \eqref{eqn:ltpi-def}. Under Assumptions \ref{ass:vp-a}, there is a constant $K:=K(s) < \infty$ such that
	\[\sup_{m\geq 1}\ \E_{\ord}\left(\left[L(\cT^{\vp}_m)\right]^s\right) \leq K.\]
	In particular, the collection of random variables $\set{L(\cT^{\vp}_m): m\geq 1}$ is uniformly integrable.
\end{prop}

\noindent{\bf Proof:} Writing out the tilt $L(\cdot)$ explicitly, we have
\begin{equation}
\label{eqn:lt-it-barl-decomp}
	\displaystyle L(\vt):= \prod_{(k,\ell)\in E(\vt)} \left[\frac{\exp(a p_k p_{\ell})- 1}{ap_k p_{\ell}} \right] \exp\bigg(\sum_{(k,\ell) \in \sP(\vt)} a p_k p_{\ell}\bigg) = \dI(\vt) \barL(\vt),
\end{equation}
say, where,
\begin{equation}
\dI(\vt):= \prod_{(i,j)\in E(\vt)} \frac{\exp(a p_i p_j)- 1}{ap_ip_j}
\leq \exp\bigg( a\sum_{(i,j)\in E(\vt)} p_i p_j  \bigg)
\leq \exp( a \ch{p_1}).
	\label{eqn:dit-bound}
\end{equation}
Here we have used $(\e^x-1)/x \leq \e^x$ for $x >0$ for the first inequality and the second inequality follows using the fact that $\vt$ is a tree, so that for each $(i,j) \in E(\vt)$ such that $i$ is the parent of $j$, we have $p_ip_j \leq \ch{p_1} p_j$. By Assumption \ref{ass:vp-a}, we have \ch{$ap_1 \to \gamma \theta_1$}. In particular, there is a constant $C> 0$ such that for all $m\geq 1$, and $\vt\in{\bT}_m^{\ord}$,
\begin{equation}
\label{eqn:lvp-bound}
	\dI(\vt)\leq C\ \text{ and }\ L(\vt)\leq C \exp\bigg(\sum_{(k,\ell) \in \sP(\vt)} a p_k p_{\ell}\bigg).
\end{equation}
Now recall the functions $A_m$ and $\bar{A}_m:= aA_m$ from \eqref{eqn:ant-def}. Using the equivalent characterization of the permitted edge set from \eqref{eqn:permitted-set} and comparing this with \eqref{eqn:ant-def}, it is easy to check that
\begin{align*}
	\sum_{(i,j) \in \sP(\cT_m^{\vp})} a p_i p_j = a \sum_{i \in [m]} \sum_{ j \in \cS(i)} p_i p_j= \int_0^1 \bar A_m(s) ds.
\end{align*}
Now by the definition of $F^{\exec,\vp}$,
\begin{equation}
	\label{eqn:sum-of-wts-active}	
	F^{\exec,\vp}(y^*(i)) = \sum_{v \in \cA(i) } p_v, \;\; \mbox{ for } i \in [m].
\end{equation}
By \eqref{eqn:ant-def},
\begin{equation*}
	A_m(t) = \sum_{v \in \cS(i)} p_v = \sum_{v \in \cA(i-1)}p_v - p_{v(i)}, \;\; \mbox{ for } t \in (y^*(i-1), y^*(i)].
\end{equation*}
Thus
\begin{align}\label{eqn:88}
\normi{A_m}\leq \normi{F^{\exec,\vp}}.
\end{align}
By Assumption \ref{ass:vp-a}(ii) and \eqref{eqn:lvp-bound}, for any $s\geq 0$, there exists $K=K(s) < \infty$ such that
\begin{equation}
\label{eqn:lvp-final}
	\left[L(\cT^{\vp}_m)\right]^s \leq C^s \exp\left(K \frac{\normi{F^{\exec,\vp}}}{\sigma(\vp)}\right).
\end{equation}
Now the following lemma completes the proof of Proposition \ref{prop:tilt-tightness}.\qed

\begin{lemma}\label{lem:tail-fexec-inf}
	There exists a positive constant $c > 0$ such that for every $m\geq 1$ and $x\geq \e$,
	\[\pr\left(\normi{F^{\exec,\vp}} \geq x\sigma(\vp)\right)\leq \exp\big(-c x\log(\log{x})\big). \]
\end{lemma}

\noindent{\bf Proof:} Write $\cR(m):= \normi{F^{\exec,\vp}}/\sigma(\vp)$ and as before, let $\vX=(X_v:v\in [m])$ be the collection of uniform random variables used to construct $F^{\vp}$. Write $\bQ[0,1]$ for the set of rationals in $[0,1]$. Then note that
\begin{equation}
\label{eqn:949}
	\cR(m)= \sup_{q\in \bQ[0,1]} \frac{F^{\vp}(q)}{\sigma(\vp)}- \inf_{q\in \bQ[0,1]} \frac{F^{\vp}(q)}{\sigma(\vp)}:= \cR_1(m)+ \cR_2(m).
\end{equation}
We start by analyzing $\cR_1(m)$. For fixed $q\in \bQ[0,1]$, define the collection of $m$ functions
\begin{equation}
\label{eqn:sqj-def}
	s_q^j(x):= \frac{p_j}{\sigma(\vp)}\left(\ind\set{x\leq q} -q\right), \qquad 1\leq j\leq m.\notag
\end{equation}
Note that  for all $j\in [m]$, $s_q^j:[0,1]\to [-1,1]$, with  $\E(s_q^j(X_j)) =0$ and further
\begin{equation}
\label{eqn:1013}
	\sR_1(m)= \sup_{q\in \bQ[0,1]} \left(s_q^1(X_1)+\cdots + s_q^m(X_m)\right).\notag
\end{equation}
Also note that
\begin{equation}
\label{eqn:var-bd}
	\sup_{q\in \bQ[0,1]} \var\left(s_q^1(X_1)+\cdots + s_q^m(X_m)\right)=\sup_{q\in \bQ[0,1]} q(1-q)= \frac{1}{4}.\notag
\end{equation}
If we can show that
\begin{equation}
\label{eqn:nts-r1}
	\kappa:= \sup_{m\geq 1} \E(\cR_1(m)) < \infty,
\end{equation}
then standard concentration inequalities for the maxima in empirical processes \cite[Theorem 1.1(b)]{klein-rio05} will imply the existence of a constant $c_1>0$ such that for all $m\geq 1$ and $x >0$,
\begin{equation}
\label{eqn:tail-r1}
	\pr(\cR_1(m)\geq \E(\cR_1(m))+x)\leq \exp\left(-\frac{x}{4} \log\left[1+ 2\log\left(1+\frac{x}{2\kappa+\frac{1}{4}}\right)\right]\right).
\end{equation}
Let us now prove \eqref{eqn:nts-r1}. In fact we will show the stronger result:
\begin{equation}
\label{eqn:mod-r1-3}
	\sup_{m\geq 1}\E\left(\sup_{q\in \bQ[0,1]} \left|\sum_{j=1}^m s_q^j(X_j)\right|\right) < \infty.\notag
\end{equation}
Let $X_{\sss(1)} < X_{\sss(2)} < \cdots < X_{\sss(m)}$ denote the order statistics of $\vX$ and let $\pi$ denote the corresponding permutation of $[m]$, namely $X_{\sss(i)} = X_{\pi(i)}$.  Note that
\begin{align*}
\sup_{q\in \bQ[0,1]} \left|\sum_{j=1}^m s_q^j(X_j)\right|:= \max_{1\leq i\leq m} |\vartheta_i|, \qquad \mbox{ where } \qquad \vartheta_i:= \frac{-X_{\sss(i)} +\sum_{j=1}^i p_{\pi(j)}}{\sigma(\vp)}.
\end{align*}
Hence
\begin{align}
\label{eqn:mod-r1-3a}
	\max_{i\in [m]} |\vartheta_i|&\leq \max_{i\in [m]}~ [\sigma(\vp)]^{-1}{\left|-X_{\sss(i)} +\frac{i}{m}\right|} ~+~ \max_{i\in [m]} [\sigma(\vp)]^{-1}{\left|\sum_{j=1}^{i} p_{\pi(j)} -\frac{i}{m}\right|}\notag\\
	&:= \cR_{11}(m)+\cR_{12}(m).\notag
\end{align}
We first analyze $\cR_{11}(m)$.  By the DKW inequality \cite{massart-DKW},
\begin{equation}
\label{eqn:dkw}
	\pr\left(\max_{i\in [m]}~ {\left|-X_{\sss(i)} +\frac{i}{m}\right|} \geq \sigma(\vp) x\right) \leq 2\exp\left(-2m\cdot\left(\sigma(\vp)x\right)^2\right)\notag
\end{equation}
By Cauchy-Schwartz, $m\sigma^2(\vp)\geq (\sum_i p_i)^2=1$. Thus $\sup_{m\geq 1} \E(\cR_{11}(m)) < \infty$. We now analyze $\cR_{12}(m)$. Since
\begin{equation}
\label{eqn:pi-pj}
\E(p_{\pi(j)}) = \frac{1}{m},\quad\text{and}\quad\E(p_{\pi(i)}p_{\pi(j)}) = \frac{\sum_{k\neq\ell\in [m]} p_k p_{\ell}}{m(m-1)}= \frac{1-\sigma^2(\vp)}{m(m-1)}\quad\text{for}\quad i\neq j\in [m],\notag
\end{equation}
for any $i\in [m]$ we have
 \begin{align}
 	\E\left(\left(\sum_{j=1}^{i} p_{\pi(j)} -\frac{i}{m}\right)^2\right) &= \frac{i\sigma^2(\vp)}{m} + \frac{i(i-1)}{m(m-1)}(1-\sigma^2(\vp)) -\frac{i^2}{m^2} \leq \frac{i}{m}\sigma^2(\vp) \label{eqn:rmm-4}
 \end{align}
by simply expanding the square. Now note that since $\pi$ is a uniform random permutation of the vertex set $[m]$, for any fixed $i\geq 1$ we also have
\begin{align*}
\sum_{j=1}^i p_{\pi(j)} - \frac{i}{m} \stackrel{d}{=} \sum_{j=0}^{i-1} p_{\pi(m-j)} - \frac{i}{m} = \left(\frac{m-i}{m} -\sum_{j=1}^{m-i} p_{\pi(i)}\right).
\end{align*}
Thus
\begin{equation}
\label{eqn:1521}
	\E(\cR_{12}(m)) \leq 2 \E\left(\max_{i\in [m/2]} [\sigma(\vp)]^{-1}{\left|\sum_{j=1}^{i} p_{\pi(j)} -\frac{i}{m}\right|}\right).
\end{equation}
Now assuming that we construct $\pi$ by sequentially sampling without replacement from $[m]$, let $\cF_k$ denote the $\sigma$-field generated by $(\pi(1), \pi(2), \ldots, \pi(k))$ for $0\leq k\leq m-1$. Let $M_0 = 0$ and consider the sequence
\[M_k:= \frac{\sum_{j=1}^k p_{\pi(j)} - k/m}{m-k}, \qquad 0\leq k\leq m-1. \]
It is easy to check that $\set{M_k:0\leq k\leq m-1}$ is a martingale with respect to the filtration $\set{\cF_k: 0\leq k\leq m-1}$. Then \eqref{eqn:1521} and Doob's $\bL^2$-maximal inequality yield
\begin{equation}
\label{eqn:1526}
	\E(\cR_{12}(m)) \leq \frac{2m}{\sigma(\vp)} \sqrt{\E(\left[M_{m/2}\right]^2)}\notag
\end{equation}
Using \eqref{eqn:rmm-4} with $i=m/2$ then gives $\E(\cR_{12}(m)) \leq 16$ for all $m\geq 1$. Thus we have shown that $\sup_{m\geq 1} \max(\E(\cR_{11}(m)),\E(\cR_{12}(m))) < \infty$. This proves \eqref{eqn:nts-r1} and thus \eqref{eqn:tail-r1}.

To complete the proof of the lemma, we need to get a tail bound on $\cR_2(m)$ appearing in \eqref{eqn:949}. As before, using \cite{klein-rio05}, it is enough to show  $\sup_{m\geq 1} \E(\cR_2(m))< \infty$. However, note that
\begin{align*}
	\cR_2(m)&= \max_{i\in [m]}\frac{\left|\sum_{j=1}^{i-1} p_{\pi(j)} - X_{\sss(i)}\right|}{\sigma(\vp)} \leq \cR_1(m) + \frac{\ch{p_1}}{\sigma(\vp)}.
\end{align*}
We now use \eqref{eqn:nts-r1} together with Assumption \ref{ass:vp-a} to complete the proof. \qed
\subsection{Another construction of $\tilde{\mathcal{G}}_m(\mathbf{p}, a)$ and a modification}
\label{sec:gmvpa-explicit}
In this section, we start by giving a more explicit description of the algorithm described in Proposition \ref{prop:SBSSXW} via adding permitted edges to a tilted $\vp$-tree. We first set up some notation.  As a matter of convention, we will view ordered rooted trees via their planar embedding, using the associated ordering to determine the relative locations of siblings of an individual. We think of the left most sibling as the ``oldest''. Further, in a depth-first exploration, we explore the tree from left to right. Now given a planar rooted tree $\vt\in \bT_m$, let $\rho$ denote the root and for every vertex $v\in [m]$, let $[\rho,v]$ denote the path connecting $\rho$ to $v$ in the tree. Given this path and a vertex $i\in [\rho,v]$, write $\RC(i,[\rho,v])$ for the set of all children of $i$ which fall to the right of $[\rho,v]$. \nomenclature[rc]{$\RC(i,[\rho,v])$}{For vertex $i$ in a  path $[\rho,v]$, set of all children of $i$ which fall to the right of $[\rho,v]$.} Thus in the depth-first exploration of the tree, when we get to $v$,
	\begin{equation}
	\label{eqn:permit-vct}
		\fP(v,\vt):= \cup_{i\in [m]} \RC(i,[\rho,v])
	\end{equation}
	denotes the set of endpoints of all permitted edges emanating from $v$. Define
	\begin{equation}
	\label{eqn:amv-def}
		\dA_{\sss(m)}(v):= \sum_{i\in [\rho,v]} \sum_{j\in [m]} p_j \ind\set{j\in \RC(i,[\rho,v])}.
	\end{equation}
The function $A_m(\cdot)$ defined in \eqref{eqn:ant-def} is intimately connected to $\dA_{\sss(m)}(\cdot)$. More precisely, let $(v(1), v(2), \ldots, v(m))$ denote the order in the depth-first exploration of the tree. Let $y^*(0)=0$ and $y^*(i) = y^*(i-1) + p_{v(i)}$.  Define
	\begin{equation}
	\label{eqn:ant-def-new}
		A_{\sss(m)}(u) = \dA_{\sss(m)}(u),\quad\text{for}\quad u\in (y^*(i-1), y^*(i)],\quad\text{and}\quad \bar{A}_{\sss(m)}(\cdot):= a A_{\sss(m)}(\cdot).
	\end{equation}
Then the function $A_{\sss(m)}(\cdot)$ associated with an ordered $\vp$-tree has the same distribution as the function $A_{m}(\cdot)$ associated with the tree $\psi_{\vp}(\vX)$, where $\vX=(X_v: v\in[m])$ are i.i.d. random variables uniformly distributed on $(0,1)$ .

Finally, define the function
		\begin{equation}
		\label{eqn:Lambda-tree}
			\Lambda_{\sss(m)}(\vt) := a\sum_{v\in [m]} p_v \dA_{\sss(m)}(v).
		\end{equation}

While all of these objects depend on the tree $\vt$, we suppress this dependence to ease notation. Now Proposition \ref{prop:SBSSXW} implies we can construct $\tilde\cG_m(\vp,a)$ via the following five steps:

\begin{enumeratei}
	\item {\bf Tilted $\vp$-tree:}
	Generate a tilted ordered $\vp$-tree $\cT^{\vp,\star}_m$ with distribution \eqref{eqn:tilt-ord-dist-def}. Now consider the (random) objects $\fP(v,\cT^{\vp,\star}_m)$ for $v\in [m]$ and the corresponding (random) functions $\dA_{\sss(m)}(\cdot)$ on $[m]$ and $A_{\sss(m)}(\cdot)$ on $[0,1]$.
		\item {\bf Poisson number of \emph{possible} surplus edges:} Let $\cP$ denote a rate one Poisson process on $\bR_+^2$ and define
	\begin{equation}
	\label{eqn:poisson-pp}
		\bar{A}_{\sss(m)}\cap {\cP}:= \set{(s,t)\in \cP: s\in [0,1], t\leq \bar{A}_{\sss(m)}(s)}.
	\end{equation}
	 Write $\bar{A}_{\sss(m)}\cap {\cP}:= \set{(s_j,t_j):1\leq j\leq N_{\sss(m)}^\star}$ where $N_{\sss(m)}^\star = |\bar{A}_{\sss(m)}\cap {\cP}|$.

We will now use the set $\set{(s_j, t_j):1\leq j\leq N_{\sss(m)}^\star}$ to generate pairs of points $\set{(\cL_j,\cR_j): 1\leq j\leq N_{\sss(m)}^\star}$ in the tree that will be joined to form the surplus edges.
	 \item {\bf ``First'' endpoints:} Fix $j$ and suppose $s_j \in (y^*(i-1), y^*(i)]$ for some $i\geq 1$, where $y^*(i)$ is as given right above \eqref{eqn:ant-def-new}. Then the \emph{first endpoint} of the surplus edge corresponding to $(s_j, t_j)$ is $\cL_j:= v(i)$.
	 \item {\bf ``Second'' endpoints:} Note that in the interval $(y^*(i-1), y^*(i)]$, the function $\bar{A}_{\sss(m)}$ is of constant height $a\dA_{\sss(m)}(v(i))$. We will view this height as being partitioned into sub-intervals of length $a p_u$ for each $u\in \fP(v(i),\cT^{\vp,\star}_m)$, the collection of endpoints of permitted edges emanating from $\cL_k$. (Assume that this partitioning is done according to some preassigned rule, e.g., using the order of the vertices in $\fP(v(i),\cT^{\vp,\star}_m)$.) Suppose $t_j$ belongs to the interval corresponding to $u$. Then the \emph{second endpoint} is $\cR_j = u$. Form an edge between $(\cL_j, \cR_j)$.
\item In this construction, it is possible that one created more than one surplus edge between two vertices. Remove any multiple surplus edges.
\end{enumeratei}

\begin{lemma}
	\label{lem:lk-rk-equivalent}
	The above construction gives a random graph with distribution $\tilde{\cG}_m(\vp,a)$ as in Definition \ref{def:con-graph-tilde-p}. Further, conditional on $\cT^{\vp,\star}_m$:
	\begin{enumeratea}
		\item $N_{\sss(m)}^\star$ has Poisson distribution with mean $\Lambda_{\sss(m)}(\cT_m^{\vp,\star})$ where $\Lambda_{\sss(m)}$ is as in \eqref{eqn:Lambda-tree}.
		\item Conditional on $\cT_m^{\vp,\star}$ and $N_{\sss(m)}^\star=k$, the first endpoints $(\cL_j: 1\leq j\leq k)$ can be generated in an i.i.d. fashion by sampling from the vertex set $[m]$ with probability distribution
		\[{\cJ}^{\sss(m)}(v) \propto p_v \dA_{\sss(m)}(v), \qquad v\in [m].\]
		\item Conditional on $\cT_m^{\vp,\star}$, $N_{\sss(m)}^\star=k$ and the first endpoints $(\cL_j: 1\leq j\leq k)$, the second endpoints can be generated in an i.i.d. fashion where the probability that $\cR_j = u$ is proportional to $p_u$ if $u$ is a right child of some individual $y\in [\rho,\cL_j]$.
	\end{enumeratea}
	
\end{lemma}

\noindent{\bf Proof:} The assertions follow from Proposition \ref{prop:SBSSXW} and standard properties of Poisson processes. \phantom{m}\hfill\qed

\medskip
\noindent{\bf The modified space $\cG_m^{\md}(\vp,a)$:} We construct a modified graph $\cG_m^{\md}(\vp,a)$ as follows:
\begin{enumeratei}
\item[(i$^\prime$)] Generate a tilted ordered $\vp$-tree $\cT_m^{\vp,\star}$ with distribution \eqref{eqn:tilt-ord-dist-def}.
\item[(ii$^\prime$)] Conditional on $\cT_m^{\vp,\star}=k$, generate $N_{\sss(m)}^\star\sim\Poi(\Lambda_{\sss(m)}(\cT_m^{\vp,\star}))$.
\item[(iii$^\prime$)] Conditional on $\cT_m^{\vp,\star}$ and $N_{\sss(m)}^\star=k$, generate the first endpoints $(\cL_j: 1\leq j\leq k)$ in an i.i.d. fashion by sampling from the vertex set $[m]$ with probability distribution
		\[{\cJ}^{\sss(m)}(v) \propto p_v \dA_{\sss(m)}(v), \qquad v\in [m].\]
\item[(iv$^\prime$)] Conditional on $\cT_m^{\vp,\star}$, $N_{\sss(m)}^\star=k$ and the first endpoints $(\cL_j: 1\leq j\leq k)$, generate the second endpoints in an i.i.d. fashion where conditional on $\cL_j = v$, the probability distribution of $\cR_j$ is given by
	\begin{equation}
	\label{eqn:right-end-pt-prob}
		Q_{v}^{\sss(m)}(y):= \begin{cases}
			\ch{\sum_{u}} p_u \ind\set{u\in \RC(y,[\rho,v])}/\dA_{\sss(m)}(v) & \text{ if } y\in [\rho,v],\\
			0 & \text{ otherwise }.
		\end{cases}
	\end{equation}	
 Identify $\cL_j$ and $\cR_j$ for $1\leq j\leq k$.
\end{enumeratei}

Thus, instead of adding an edge between $\cL_j$ and one of the right children on the path $[\rho,\cL_j]$ as in Lemma \ref{lem:lk-rk-equivalent}(c), we identify it to the parent of this vertex which is on $[\rho,\cL_j]$. Also, we do not remove any multiple surplus edges. This construction turns out to be easier to work with. $\cG_m^{\md}(\vp,a)$ will be viewed as a metric measure space via the graph distance where vertex $v$ has mass $\sum p_u$ where the sum is taken over all $u\in[m]$ which have been identified with $v$.
Intuitively it is clear that $\sigma(\vp)\tilde\cG_m(\vp,a)$ and $\sigma(\vp)\cG_m^{\md}(\vp,a)$ are ``close''. This is formalized in  Lemma \ref{lem:ghp-distance-negligible-1}.
\begin{rem}
	At this point we urge the reader to go back to Section \ref{sec:tilt-icrt} and remind themselves of the four steps in the construction of the limit metric space $\cG_{\infty}(\mvtheta,\gamma)$, and note the similarities to the construction above. In particular, we make note of the following:
\begin{enumeratea}
\item For finite $m$, we essentially tilt the $\vp$-tree distribution via the functional $\bar L(\cT_m^{\vp}) = \exp(a\E[\dA_{\sss(m)}(V_1)\ |\ \cT_{m}^{\vp}])$ (the term $\dI(\cT_m^{\vp})$ as in \eqref{eqn:lt-it-barl-decomp} can be ignored as we will see in Lemma \ref{lem:it-go-one}), and the number of shortcut points selected, namely $N_{\sss(m)}^\star$, has a Poisson distribution with mean $a\E(\dA_{\sss(m)}(V_1)\ |\ \cT_{m}^{\vp,\star})$. Here $V_1$ has distribution $\vp$.
 \item For the limit object, we tilt the measure using the functional $L_{\sss(\infty)}(\icrt, \mvU) = \exp(\gamma\E[\dA_{\sss(\infty)}(V_1)\ |\ \icrt, \mvU])$, and the number of shortcuts, namely $N_{\sss(\infty)}^\star$, follows a Poisson distribution with mean $\gamma \E(\dA_{\sss(\infty)}(V_1)\ |\ \tilicrt, \mvU^{\star})$. Here $V_1$ is distributed according to the mass measure $\mu^\star$ on $\tilicrt$.
 \end{enumeratea}
\end{rem}

As a brief warm-up to the kind of calculations in the next section, we now prove a simple lemma on tightness of the number of surplus edges. We will prove distributional convergence of this object in the next section.

\begin{lemma}\label{lem:nm-tight}
Under Assumption \ref{ass:vp-a}, the sequence $\set{N_{\sss(m)}^\star:m\geq 1}$ is tight, where $N_{\sss(m)}^\star$ is as given below \eqref{eqn:poisson-pp}.
\end{lemma}

\noindent {\bf Proof:} Fix $r > 1$.
First note that conditional on $\cT_m^{\vp,\star} =\vt$, $N_{\sss(m)}^\star$ has a Poisson distribution with mean $\Lambda_{\sss(m)}(\vt)$. Thus, there exists a constant $C = C(r) $ such that
\begin{equation}
\label{eqn:poisson-moment-bd}
	\E([N_{\sss(m)}^\star]^r|\cT_m^{\vp,\star} =\vt) \leq C [\Lambda_{\sss(m)}(\vt)]^r.\notag
\end{equation}
 Further, note that the tilt $L(\vt)$ in \eqref{eqn:ltpi-def} satisfies
\begin{equation}
\label{eqn:lvt-decomp}
	L(\vt):= \dI(\vt) \exp\left(\sum_{(k,\ell) \in \sP(\vt)} a p_k p_{\ell}\right) = \dI(\vt) \exp(\Lambda_{\sss(m)}(\vt)),\notag
\end{equation}
where $1\leq \dI(\vt) \leq C^\prime$ for a fixed constant $C^\prime$ independent of $m$ by \eqref{eqn:dit-bound}. Thus, Proposition \ref{prop:tilt-tightness} shows that
\begin{equation}
\label{eqn:sup-lm-bd}
	\sup_{m\geq 1}\ \E\big(\exp(\gamma \Lambda_{\sss(m)}(\cT^{\vp}_m))\big) < \infty\notag
\end{equation}
for any $\gamma>0$. In particular,
\begin{align}
\sup_{m\geq 1}\ \E([N_{\sss(m)}^\star]^r)& \leq\sup_{m\geq 1}\ C\E\left([\Lambda_{\sss(m)}(\cT^{\vp,\star}_m)]^r\right)
	 =C\sup_{m\geq 1}\ \frac{\E([\Lambda_{\sss(m)}(\cT^{\vp}_m)]^r L(\cT^{\vp}_m)))}{\E(L(\cT^{\vp}_m))}<\infty,\notag
\end{align}
which proves tightness of $\set{N_{\sss(m)}^\star:m\geq 1}$.\qed

We conclude this section by proving a lemma which essentially says that it is enough to work with the modified space $\cG_m^{\md}(\vp,a)$.
\begin{lem}\label{lem:ghp-distance-negligible-1}
Recall the five-step construction of $\tilde{\cG}_m(\vp, a)$. Construct $\cG_m^{\md}(\vp,a)$ on the same space by coupling it with $\tilde{\cG}_m(\vp, a)$ in the obvious way.
Then, under Assumption \ref{ass:vp-a},
\begin{align*}
d_{\GHP}\left(\sigma(\vp)\tilde\cG_m(\vp,a),\ \sigma(\vp)\cG_m^{\md}(\vp,a)\right)\probc 0.
\end{align*}
\end{lem}
\noindent {\bf Proof:} Define the event
\begin{align*}
F:=\set{N_{\sss(m)}^\star \text{ equals the number of surplus edges in }\tilde{\cG}_m(\vp, a)}.
\end{align*}
In other words, $F$ describes the event in which $\tilde{\cG}_m(\vp, a)$ does not have multiple surplus edges. It is easy to check that
\begin{align*}
d_{\GHP}\left(\tilde\cG_m(\vp,a),\ \cG_m^{\md}(\vp,a)\right)\leq N_{\sss(m)}^\star\ \text{ on the set }\ F.
\end{align*}
Thus, Lemma \ref{lem:nm-tight} combined with the assumption $\sigma(\vp)\to 0$ yields the result provided we show that $\pr(F^c)\to 0$.
To this end, note that
\begin{align*}
&\pr\left(\exists\text{ multiple surplus edges between }u\text{ and }v\big|\cT^{\vp,\star}_m=\vt\right)\\
&\hskip180pt=\pr\left(\Poi(ap_u p_v)\geq 2\right)
\leq c(ap_u p_v)^2
\end{align*}
for every $u\in[m]$, $v\in\fP(u,\vt)$, and some universal positive constant $c$. Hence
\begin{align*}
\pr\left(F^c\mid \cT^{\vp,\star}_m=\vt\right)
&\leq ca^2\sigma(\vp)^2\sum_{u\in[m]}p_u^2\sum_{v\in\fP(u,\vt)}\left(p_v/\sigma(\vp)\right)^2\\
&\leq c(a\sigma(\vp))^2\sum_{u\in[m]}p_u^2=c(a\sigma(\vp))^2 \sigma(\vp)^2.
\end{align*}
Since $\sigma(\vp)\to 0$ and $a\sigma(\vp)\to\gamma$, $\pr(F^c)\to 0$ as desired.
\qed

\subsection{Completing the proof of Theorem \ref{thm:connected-convg}.}
\label{sec:completing-proof-thm-con} At this point we urge the reader to remind themselves of (a) the four steps in the construction of the limit object in Section \ref{sec:descp-limit}, (b) the birthday construction of $\vp$-trees at the end of Section \ref{sec:ptree-const} and (c) the definition of Gromov-weak topology in Section \ref{sec:gromov-weak}  of complete separable measured metric spaces $\sS_*$. Fix $\ell\geq 1$ and a bounded continuous function  $\phi:\bR_+^{\ell^2}\to \bR$. Let $\Phi$ be as in \eqref{eqn:polynomial-func-def}. To simplify notation, we will write $\Phi(X)$ instead of $\Phi(X,d,\mu)$.  To prove Theorem \ref{thm:connected-convg}, we need to show that for every fixed $\ell\geq 1$ and functions $\phi$ and $\Phi$ as above,
\begin{align*}
	\E\left[\Phi\left(\sigma(\vp)\cdot\tilde{\cG}_m(\vp,a)\right)\right] \to \E\left[\Phi\left(\cG_{\infty}(\mvtheta,\gamma)\right)\right]\qquad \mbox{ as } m\to\infty,
\end{align*}
where we sample $\ell$ points according to $\vp$ in $\tilde{\cG}_m(\vp,a)$ while we sample $\ell$ points according to the measure on $\cG_{\infty}(\mvtheta,\gamma)$ inherited from the mass measure.  Now recall the explicit five step construction of $\tilde{\cG}_m(\vp,a)$ in Section \ref{sec:gmvpa-explicit} starting from the tilted $\vp$-tree $\cT_m^{\vp,\star}$ and the Poisson number of surplus edges $N_{\sss(m)}^\star$. Fix $K\geq 1$ and note that
\begin{align*}
\left|\E\left[\Phi\left(\sigma(\vp)\tilde{\cG}_m(\vp,a)\right)\right] - \sum_{k=0}^K \E\left[\Phi\left(\sigma(\vp)\tilde{\cG}_m(\vp,a)\right)\ind\set{N_{\sss(m)}^\star =k}\right]\right| \leq ||\phi||_{\infty} \pr(N_{\sss(m)}^\star\geq K+1).
\end{align*}
Using Lemma \ref{lem:nm-tight}, we can choose $K$ large (independent of $m$) to make the bound on the right arbitrarily small. Further, in view of Lemma \ref{lem:ghp-distance-negligible-1}, we can work with ${\cG}_m^{\md}(\vp,a)$ instead of $\tilde{\cG}_m(\vp,a)$. Hence it suffices to prove the following convergence for every fixed $k\geq 0$:
\begin{align}\label{eqn:999}
\E\left[\Phi\left(\sigma(\vp){\cG}_m^{\md}(\vp,a)\right)\ind\set{N_{\sss(m)}^\star =k}\right]\to
\E\left[\Phi\left({\cG}_{\infty}(\mvtheta, \gamma)\right)\ind\set{N_{\sss(\infty)}^\star =k}\right]\ \text{ as }\ m\to\infty.
\end{align}
To analyze this term, we first need to setup some notation.

 Note that both the finite $m$ and the limit object are obtained by starting with a discrete tree for finite $m$ and a \emph{real tree} in the limit, and sampling a random number of pairs to create ``shortcuts''. Recall the space $\vT_{IJ}^*$ in Section \ref{sec:space-of-trees}. Fix $k\geq 0$ and let $\vt$ be an element in $\vT_{I,(k+\ell)}^*$ for some $I\geq 0$. ``$I$" will not play a role in the definition below. Write $\rho$ for the root and denote the leaves by 
\begin{equation}
\label{eqn:leaves-def}
\vx_{k,k+\ell}:= (x_1, x_2, \ldots,x_{k}, x_{k+1}, \ldots, x_{k+\ell}).
\end{equation}
Also recall that for each $i$, there is a probability measure $\nu_{\vt,i }(\cdot)$ on the path  $[\rho, x_i]$ for $1\leq i\leq k+\ell$. For $1\leq i\leq k$, sample $y_i$ according to the distribution  $\nu_{\vt,i}(\cdot)$ independently for different $i$ and connect $x_i$ and $y_i$. Let $\vt'$ denote the (random) tree thus obtained and let $d_{\vt'}$ denote the graph distance on $\vt'$. Define the function $g^{\sss(k)}_\phi:\vT_{I,(k+\ell)}^*\to \bR$ by
\begin{align}
\label{eqn:gphi-def}
g^{\sss(k)}_{\phi}(\vt):=
\left\{
\begin{array}{l}
\E\left[\phi\left(d_{\vt'}(x_i, x_j): k+1\leq i\leq k+\ell\right)\right],\text{ if }\vt\neq\partial,\\
0,\text{ if }\vt=\partial.
\end{array}
\right.
\end{align}
\nomenclature[gphi]{$g^{\sss(k)}_{\phi}(\vt)$}{For a tree $\vt\in \vT_{I,(k+\ell)}^*$, the functional defined in \eqref{eqn:gphi-def}.}
In words, we look at the expectation of $\phi$ applied to the pairwise distances between the last $\ell$ leaves after sampling $y_i$ on the path $[\rho, x_i]$ for $1\leq i\leq k$ and connecting $x_i$ and $y_i$. Note that here the expectation is only taken over the choices of $y_i$.

Next, given $\vt\in\bT_m^{\ord}$ and $\mvv:=(v_1, \hdots, v_{r})$ with $v_i\in[m]$, set $\vt(\mvv)$ to be the subtree of $\vt$ spanning the vertices $\mvv$ and the root provided $v_1, \hdots, v_{r}$ are all distinct and none of them is an ancestor of another vertex in $\mvv$. When this condition fails, set $\vt(\mvv)=\partial$.

Now, conditional on $\cT_m^{\vp,\star}$, construct a tree \nomenclature[span-tree]{$\cT_m^{\vp,\star}(\widetilde{\vV}_{k,k+\ell}^{\sss (m)})$}{Spanning subtree of tilted $\vp$-tree $\cT_m^{\vp,\star}$ using sampled vertex set $\widetilde{\vV}_{k,k+\ell}^{\sss (m)}$.} $\cT_m^{\vp,\star}(\widetilde{\vV}_{k,k+\ell}^{\sss (m)})$ where
\begin{enumeratei}
\item $\widetilde{\vV}_{k,k+\ell}^{\sss (m)}:= (\barV_1^{\sss (m)}, \ldots, \barV_k^{\sss (m)},V_{k+1}^{\sss (m)},\ldots V_{k+\ell}^{\sss (m)})$;
\item $\barV_i^{\sss (m)}$, $1\leq i\leq k$ are i.i.d. with the distribution $\cJ^{\sss(m)}(\cdot)$ as in Lemma \ref{lem:lk-rk-equivalent}(b); and
\item $V_{k+1}^{\sss (m)}, \ldots V_{k+\ell}^{\sss (m)}$ are i.i.d. with distribution $\vp$. Further, $\barV_1^{\sss (m)}, \ldots, \barV_k^{\sss (m)},V_{k+1}^{\sss (m)},\ldots V_{k+\ell}^{\sss (m)}$ are jointly independent.
\end{enumeratei}

We will drop the superscript and simply write \ch{$V_i$, $\barV_i$} etc. when there is no scope of confusion.
Note that $\cT_m^{\vp,\star}(\widetilde{\vV}_{k,k+\ell})=\partial$ whenever $\barV_1, \ldots, \barV_k,V_{k+1},\ldots V_{k+\ell}$ are not all distinct or one of them is an ancestor of another vertex in $\widetilde{\vV}_{k,k+\ell}$. In either of these two case, the subtree spanned by the root and $\widetilde{\vV}_{k,k+\ell}$ will have less than $k+\ell$ leaves. We made the convention of setting $\cT_m^{\vp,\star}(\widetilde{\vV}_{k,k+\ell})=\partial$ to make sure that we are always working with a bona fide element in $\vT_{I,(k+\ell)}^*$. However, this makes no difference at all since by \cite[Corollary 15]{pitman-camarri}, \begin{align*}
\lim_m\ \pr\left(\cT_m^{\vp}(V_1,\hdots,V_{k+\ell})=\partial\right)=0
\end{align*}
where $V_1,\hdots,V_{k+\ell}$ are i.i.d. $\vp$ random variables. Now $\cT_m^{\vp,\star}$ is obtained by tilting the distribution of $\cT_m^{\vp}$, where the tilt $L(\cdot)$  is uniformly integrable (Proposition \ref{prop:tilt-tightness}). Further, $\barV_i$, $1\leq i\leq k$ are i.i.d. with the distribution $\cJ^{\sss(m)}(v)\propto p_v\dA_{\sss(m)}(v)$ where $\max_v \dA_{\sss(m)}(v)$ is stochastically dominated by $\normi{F^{\exec, \vp}}$ (see \eqref{eqn:88} and the discussion below \eqref{eqn:ant-def-new}). It thus follows that
\begin{align}\label{eqn:conventional-state-doesnt-matter}
\lim_m\ \pr\left(\cT_m^{\vp,\star}(\widetilde{\vV}_{k,k+\ell})=\partial\right)=0.
\end{align}

Using \eqref{eqn:conventional-state-doesnt-matter}, we see that
\begin{align}
\label{eqn:n-eq-k-gphi}
&\E\left[\Phi\left(\sigma(\vp){\cG}_m^{\md}(\vp,a)\right)\ind\set{N_{\sss(m)}^\star =k}\right] \\ &\hskip80pt=\E\left\{\Epst\left[g_\phi^{\sss(k)}\left(\sigma(\vp)\cT_m^{\vp,\star}(\widetilde{\vV}_{k,k+\ell})\right)\right]
\ind\set{N_{\sss(m)}^\star =k}\right\}+o(1),\notag	
\end{align}
where $\Epst(\cdot):= \E(\cdot|\cT_{m}^{\vp,\star})$. At this point, we also define $\Ep(\cdot):=\E(\cdot|\cT_{m}^{\vp})$ where $\cT_{m}^{\vp}$ has the original ordered $\vp$-tree distribution \eqref{eqn:ordered-p-tree-def}.
\nomenclature{$\Ep, \Epst$}{Expectation conditional on the ordered $\vp$-tree $\cT_{m}^{\vp}$ and the tilted $\vp$-tree $\cT_{m}^{\vp,\star}$ respectively.}

Now since $\cJ^{\sss(m)}(v)\propto p_v \dA_{\sss(m)}(v)$, we see that the inner expectation in \eqref{eqn:n-eq-k-gphi} can be simplified as
\begin{equation}
\label{eqn:1038}
	\Epst\left[g_\phi^{\sss(k)}\left(\sigma(\vp)\cT_m^{\vp,\star}(\widetilde{\vV}_{k,k+\ell})\right)\right]= \frac{\Epst\left[\prod_{i=1}^k \dA_{\sss(m)}(V_i) g_\phi^{\sss(k)}\left(\sigma(\vp)\cT_m^{\vp,\star}({\vV}_{k,k+\ell})\right) \right]}{[\Epst(\dA_{\sss(m)}(V_1)]^k},
\end{equation}
where $\vV_{k, k+\ell} = (V_1, V_2, \ldots V_{k+\ell})$, and $V_i$ are i.i.d. with distribution $\vp$. Since $\cT_m^{\vp,\star}$ is sampled according to a tilted $\vp$-tree distribution, combining \eqref{eqn:n-eq-k-gphi}, and \eqref{eqn:1038}, we get the following result:
\begin{lemma}\label{lem:form-g-phi}
Fix $k\geq 0$. Then
\begin{align}\label{eqn:33}
&\E\left[\Phi\left(\sigma(\vp){\cG}_m^{\md}(\vp,a)\right)\ind\set{N_{\sss(m)}^\star =k}\right]\\
&\hskip30pt = C_m\E\left[\frac{\Ep\left[\left(\prod_{i=1}^k \dA_{\sss(m)}(V_i)\right) g_\phi^{\sss(k)}\left(\sigma(\vp)\cT_m^{\vp}({\vV}_{k,k+\ell})\right)  \right]}{[\E_{\vp}(\dA_{\sss(m)}(V_1)]^k} L(\cT_{m}^{\vp})\ind\set{N_{\sss(m)} =k} \right]+o(1), 	\notag
\end{align}
	where $C_m = \set{\E(L(\cT_m^{\vp}))}^{-1}$, and $L$ is the tilt as in \eqref{eqn:ltpi-def}. Further, conditional on $\cT_m^{\vp}$, $N_{\sss(m)}$ has a Poisson distribution with mean $\Lambda_{\sss(m)}(\cT_m^{\vp})= a\E_{\vp}(\dA_{\sss(m)}(V))$ as in \eqref{eqn:Lambda-tree}, where $V$ has distribution $\vp$ independent of $\cT_m^{\vp}$.
\end{lemma}
This formula will be the starting point to prove \eqref{eqn:999}. Recall from \eqref{eqn:lt-it-barl-decomp} that the tilt $L(\cdot) = \ch{\dI(\cdot)} \barL(\cdot)$, where $\ch{\dI(\cdot)}$ has a messy form given by \eqref{eqn:dit-bound}. We have already seen in \eqref{eqn:lvp-bound} that under Assumption \ref{ass:vp-a}, $\dI(\cdot)\leq C$ for a constant $C$ all $m\geq 1$. The following lemma coupled with dominated convergence theorem
will now imply that we can replace $L$ with $\barL$ in Lemma \ref{lem:form-g-phi} and in all the subsequent analysis below:

\begin{lemma}\label{lem:it-go-one}
	Under Assumption \ref{ass:vp-a},  $\dI(\cT_m^{\vp}) \probc 1$ as $m\to\infty$.
\end{lemma}

 \noindent {\bf Proof:} By \eqref{eqn:dit-bound} we have
 $1\leq \dI(\cT_m^{\vp}) \leq \exp(a \sum_{(k,l)\in E(\cT_m^{\vp})} p_k p_l)$. Thus it is enough to show that $a\E(\sum_{(k,l)\in E(\cT_m^{\vp})} p_k p_l) \to 0$.
 Now for $k\neq l\in [m]$, write $\set{k\leadsto l}$ for the event in which $l$ is a child of $k$ in $\cT_m^{\vp}$. Then standard properties of $\vp$-trees \cite[Section 6.2]{pitman-random-mappings} implies that for $k\neq l_1\neq l_2\in [m]$
 \begin{equation}
 \label{eqn:p-tree-connection}
 	\pr(k\leadsto l_1) = p_k, \qquad \pr(k\leadsto l_1 \mbox{ and } k\leadsto l_2) = p_k^2.
 \end{equation}
 Thus
 \[a\E\left(\sum_{(k,l)\in E(\cT_m^{\vp})} p_k p_l\right) = a\sum_{k=1}^m p_k \sum_{l\neq k}p_l p_k \leq a\sum_{k=1}^m p_k^2 = a[\sigma(\vp)]^2\to 0, \]
 as $m\to\infty$ by Assumption \ref{ass:vp-a}. \qed



Write $\Et$ for expectation conditional on $\icrt$ and the random variables $U_j^{\sss(i)}$ that encode the order on $\icrt$, i.e.,
\begin{align*}
\Et(\cdot):=\E\left(\cdot\ \big|\ \icrt,\ \mvU\right),
\end{align*}
\nomenclature{$\Et$}{Expectation conditional on $\icrt$ and the random variables $U_j^{\sss(i)}$ that encode the order on $\icrt$.}
and note that $\E\left[\Phi\left({\cG}_{\infty}(\mvtheta, \gamma)\right)\ind\set{N_{\sss(\infty)}^\star =k}\right]$ has an expression similar to \eqref{eqn:33}. Indeed, from the construction of ${\cG}_{\infty}(\mvtheta, \gamma)$ given in Section \ref{sec:tilt-icrt}, it follows that
\begin{align}\label{eqn:33A}
&\E\left[\Phi\left({\cG}_{\infty}(\mvtheta, \gamma)\right)\ind\set{N_{\sss(\infty)}^\star =k}\right]\\
&\hskip30pt = C_{\infty}\E\left[\frac{\Et\left[\left(\prod_{i=1}^k \dA_{\sss(\infty)}(V_i^{\sss(\infty)})\right) g_\phi^{\sss(k)}\left(\icrt({\vV}_{k,k+\ell}^{\sss(\infty)})\right)  \right]}{\left[\Et(\dA_{\sss(\infty)}(V_1^{\sss(\infty)})\right]^k} L_{\sss(\infty)}(\icrt, \mvU)\ind\set{N_{\sss(\infty)} =k} \right], 	\notag
\end{align}
where
\begin{inparaenuma}
\item $\dA_{\sss(\infty)}(\cdot)$ is as defined in \eqref{eqn:limit-da}
\item $L_{\sss(\infty)}(\icrt, \mvU)$ is as in \eqref{eqn:ltheta-def},
\item $C_{\infty}=[\E L_{\sss(\infty)}(\icrt, \mvU)]^{-1}$,
\item $V_i^{\sss(\infty)}$ are i.i.d. random variables sampled from $\icrt$ using the mass measure $\mu$,
\item ${\vV}_{k,k+\ell}^{\sss(\infty)}=(V_1^{\sss(\infty)},\hdots, V_{k+\ell}^{\sss(\infty)})$,
\item $\icrt({\vV}_{k,k+\ell}^{\sss(\infty)})$ is the tree spanned by the root of $\icrt$ and ${\vV}_{k,k+\ell}^{\sss(\infty)}$, viewed as an element of $\vT_{0, k+\ell}^{\ast}$ by declaring the leaf values to be $\dA_{\sss(\infty)}(V_j^{\sss(\infty)})$ and the root-to-leaf measures to be $Q_{V_j}^{\sss(\infty)}(\cdot)$ as in \eqref{eqn:right-end-prob-inft}, and
\item conditional on $(\icrt, \mvU)$, $N_{\sss(\infty)}$ has a Poisson distribution with mean
\begin{align*}
\Lambda_{\sss(\infty)}:=\gamma\int_{y\in\icrt} \dA_{\sss(\infty)}(y)\mu(dy)=\Et\left[\dA_{\sss(\infty)}(V_1^{\sss(\infty)})\right].
\end{align*}
\end{inparaenuma}

Finally, observe that $L_{\sss (m)}(\cdot)=\dI_{\sss (m)}(\cdot)\barL_{\sss (m)}(\cdot)$ where $\bar{L}_{\sss(m)}(\vt)=\exp(a\Ep[\dA_{\sss(m)}(V_1^{\sss(m)})])$, and recall that $a\sigma(\vp)\to \gamma$ (Assumption \ref{ass:vp-a}) and ${L}_{\sss(m)}(\cT_m^{\vp})$ is uniformly integrable (Proposition \ref{prop:tilt-tightness}). Therefore, combining Lemma \ref{lem:it-go-one}, Lemma \ref{lem:form-g-phi} and \eqref{eqn:33A} with Theorem \ref{thm:jt-convg-E} stated below yields \eqref{eqn:999} and thus completes the proof of Theorem \ref{thm:connected-convg}.

 \begin{thm}
 	\label{thm:jt-convg-E}
	For each $k\geq 0$,
	\begin{align}\label{eqn:jt-convg}
&\left(\Ep\left[\frac{\dA_{\sss(m)}(V_1^{\sss(m)})}{\sigma(\vp)}\right],  ~\Ep\left[\left(\prod_{i=1}^k \frac{\dA_{\sss(m)}(V_i^{\sss(m)})}{\sigma(\vp)}\right) g_\phi^{\sss(k)}\left(\sigma(\vp)\cT_m^{\vp}({\vV}_{k,k+\ell}^{\sss(m)})\right)  \right]\right) \stackrel{d}{\longrightarrow}\\
&\hskip100pt		\left(\Et\left[\dA_{\sss(\infty)}(V_1^{\sss(\infty)})\right],~ \Et\left[\left(\prod_{i=1}^k {\dA_{\sss(\infty)}(V_i^{\sss(\infty)})}\right) g_\phi^{\sss(k)}\left(\icrt({\vV}_{k,k+\ell}^{\sss(\infty)})\right)  \right]\right).\notag
	\end{align}
 \end{thm}
 The proof of this theorem is accomplished via the following two theorems for which we need to set up some notation. Fix $I\geq 0$ and $J\geq 1$. We will assume that $\cT_m^{\vp}$ has been constructed via the birthday construction (see Section \ref{sec:ptree-const}). This construction gives rise to an unordered $\vp$-tree. To obtain an ordered $\vp$-tree from this, let $\cD_{\sss(m)}(i)$ denote the set of children of $i$ in the $\vp$-tree for every vertex $i$. Generate i.i.d. uniform random variables $\mvU_{\sss(m)}(i):=\set{U_{\sss(m),i}(v): v\in \cD_{\sss(m)}(i)}$, independent across $v\in \cT_m^{\vp}$. Think of these as ``ages'' of the children and arrange the children from left to right in decreasing order of their ages. We can construct the function $\dA_{\sss(m)}(\cdot)$ as in \eqref{eqn:amv-def} once this ordering has been defined.

 Now recall that the right hand side of \eqref{eqn:p-tree-joint-sample} tells us how to sample $J$ i.i.d. points $(V_1^{\sss(m)}, \ldots, V_J^{\sss(m)})$ from distribution $\vp$ and the corresponding spanning subtree $\cT_J^{\cB}$ from the tree  using the repeat time sequence $\set{R_k^{\sss(m)}: k\geq 1}$. Thus, by the $J$th repeat time $R_J$, we would have sampled all $J$ vertices $V_i^{\sss(m)} = Y_{R_{i}-1}$. 
 View $\cT_J^{\cB}$ as a tree with edge lengths and marked vertices as follows:
 \begin{inparaenuma}
    \item rescale every edge to have length $\sigma(\vp)$;
    \item relabel $V_j$ as $j+$ and the root as $0+$;
 	\item mark only those vertices $i\leq I$ which occur in $\cT_J^{\cB}$;
    \item for all $1\leq j \leq J$, set the leaf values to be $\dA_{\sss(m)}(V_j)/\sigma(\vp)$, and assign the measure $\nu^{\sss(m)}_j:= Q^{\sss(m)}_{V_j}$ as defined in \eqref{eqn:right-end-pt-prob} to the path connecting the root to $V_j$, i.e, to the path $[0+, j+]$ .
 \end{inparaenuma}
  \begin{defn}\label{def:tree-rij-m}
  	Fix $I\geq 0, J\geq 1$ and consider the tree constructed as above.
   Set $r_{IJ}^{\sss(m)}=\cR_{IJ}^{\sss(m)}=\partial$ if some $j+$ is not a leaf or if some leaf has been multiply labeled. Otherwise, write $r_{IJ}^{\sss(m)}\in \vT_{IJ}$ for the tree with edge lengths and at most $I$ labelled hubs, namely where we retain information in (a) and (b) above.   Write $\cR_{IJ}^{\sss(m)}\in \vT_{IJ}^*$ for the tree where we retain all information (a)-(d) above, namely the leaf values $\dA_{\sss(m)}(V_j)$ and the root-to-leaf probability measures $Q_{V_j}^{\sss(m)}(\cdot)$ in addition to (a) and (b).
  \end{defn}
\nomenclature[rIJ]{$r_{IJ}^{\sss(m)},\cR_{IJ}^{\sss(m)}$}{Spanning subtrees obtained from the birthday construction of $\vp$-trees and retaining specific set of information. See Definition \ref{def:tree-rij-m}. See Section \ref{sec:p-tree-ICRT-def} for corresponding objects for ICRTs.  }
  Now recall the tree $\cR_{IJ}^{\sss(\infty)}$ defined in Section \ref{sec:p-tree-ICRT-def} using the limit ICRT $\icrt$. The main ingredients in the proof of Theorem \ref{thm:jt-convg-E} are the following two theorems:

  \begin{thm}\label{thm:rij-convg}
  	Under Assumption \ref{ass:vp-a}, $\cR_{IJ}^{\sss(m)} \convd \cR_{IJ}^{\sss(\infty)}$ as $m\to\infty$ for every fixed $I\geq 0$ and $J\geq 1$. This convergence is with respect to the topology defined on $\vT^*_{IJ}$ in Section \ref{sec:space-of-trees}.
  \end{thm}

  The second result we will need is as follows.  Recall the function $g_\phi^{(k)}$ on $\vT_{I,(k+\ell)}^*$ as in \eqref{eqn:gphi-def}.

  \begin{thm}\label{thm:g-phi}
  	Fix $I\geq 0$,  $k\geq 0$, $\ell\geq 2$ and a bounded continuous function $\phi$ on $\bR^{\ell^2}$. Then the function $g_{\phi}^{\sss(k)}$ is continuous on $\vT_{I, (k+\ell)}^*$.
  \end{thm}

 \noindent{\bf Proof of Theorem \ref{thm:jt-convg-E}:} Assuming Theorems \ref{thm:rij-convg} and \ref{thm:g-phi}, let us now show how this completes the proof. Getting a handle directly on the conditional expectations as required in Theorem \ref{thm:jt-convg-E} is a little tricky. Naturally, conditional on $\cT_m^{\vp}$, repeated sampling of vertices and calculating sample averages should give a good idea of the conditional expectations (and the same for the limit object $\icrt$). This is made precise in the following simple lemma whose proof we leave to the reader.

\begin{lemma}\label{lem:approx-m-inf}
	Suppose $\vX^{\sss(m)} := (X^{\sss(m),1},X^{\sss(m),2})$ with $m\in \set{1,2,\ldots,}\cup\set{\infty}$ is a sequence of $\bR^2$-valued random variables such that for each fixed $r\geq 1$, there exist random variables $\vX_r^{\sss(m)}:=(X_r^{\sss(m),1} , X_r^{\sss(m),2})$ such that the following hold:
	\begin{enumeratei}
		\item There exists a constant $C <\infty$ such that for any $m\in \set{1,2,\ldots,}\cup\set{\infty}$, $r\geq 1$ and $\eps>0$,
		\[\max_{s=1,2}\ \pr\left(|X^{\sss(m),s} - X_r^{\sss(m),s}| > \eps\right) \leq \frac{C}{\eps^2 r}. \]
		\item For each fixed $r\geq 1$, $\vX_r^{\sss(m)}\convd \vX_r^{\sss(\infty)}$.
	\end{enumeratei}
	Then $\vX^{\sss(m)}\convd \vX^{\sss(\infty)}$.
\end{lemma}
  We will apply this lemma with the random variables that arise in Theorem \ref{thm:jt-convg-E}. That is, we set
  \begin{equation}
  \label{eqn:xm-xinf}
  	X^{\sss(m),1}:= \Ep\left[\frac{\dA_{\sss(m)}(V^{\sss(m)})}{\sigma(\vp)}\right], \qquad X^{\sss(\infty),1}:= \Et\left[\dA_{\sss(\infty)}(V^{\sss(\infty)})\right],
  \end{equation}
 and similarly define $X^{\sss(m),2}$ and $X^{\sss(\infty),2}$ to be the second coordinates in the display \eqref{eqn:jt-convg}. To define $\vX_r^{\sss(m)}$, we proceed as follows. For each fixed $r\geq 1$, sample a collection of $J_r:= [r+(k+\ell)r]$ points all i.i.d. $\vp$ from $\cT_m^{\vp}$ and think of them as $r$ individuals points-$(V_1^{\sss(m)}, V_2^{\sss(m)},\ldots, V_r^{\sss(m)})$, and $r$ $(k+\ell)$ dimensional vectors-$\vV^{\sss(m),i}_{k,k+\ell}:= (V_{i1}^{\sss(m)},\ldots, V_{i(k+\ell)}^{\sss(m)})$ for $1\leq i\leq r$. Define
 \[H_{\phi}^{\sss(m)}(i):= \prod_{j=1}^k \frac{\dA_{\sss(m)}(V_{ij}^{\sss(m)})}{\sigma(\vp)} g_\phi^{\sss(k)}\left(\sigma(\vp)\cT_m^{\vp}({\vV}_{k,k+\ell}^{\sss(m),i})\right),\quad\text{for}\quad 1\leq i\leq r. \]
 For $m=\infty$, sample as above $J_r$ points using the mass measure $\mu$ from $\icrt$ and define
 \[H_{\phi}^{\sss(\infty)}(i):= \prod_{j=1}^k {\dA_{\sss(\infty)}(V_{ij}^{\sss(\infty)})} g_\phi^{\sss(k)}\left(\icrt({\vV}_{k,k+\ell}^{\sss(\infty),i})\right),\quad\text{for}\quad 1\leq i\leq r. \]
 Now define
 \begin{align} \label{eqn:xr-def}
 &X_r^{\sss(m),1}:= \frac{\sum_{i=1}^r\dA_{\sss(m)}(V_i^{\sss(m)})}{r\sigma(\vp)}\ \text{ for }\ m\in\set{1, 2,\hdots},
 \qquad X_r^{\sss(\infty),1}:= \frac{\sum_{i=1}^r\dA_{\sss(\infty)}(V_i^{\sss(\infty)})}{r},\text{ and}\notag\\
 &\hskip80pt X_r^{\sss(m),2}:= \frac{\sum_{i=1}^r H_{\phi}^{\sss(m)}(i)}{r}\ \text{ for }\ m\in\set{1, 2,\hdots}\cup\set{\infty}.
 \end{align}
 Let $\vX_r^{\sss(m)}:=(X_r^{\sss(m),1},X_r^{\sss(m),2})$ for $m\in\set{1, 2,\hdots}\cup\set{\infty}$. To complete the proof of the theorem, we have to check the two conditions of Lemma \ref{lem:approx-m-inf}. Let us check condition (i) of Lemma \ref{lem:approx-m-inf} for the first coordinate. The second coordinate can be handled in an identical fashion.

 Applying Chebyshev's inequality conditional on $\cT_m^{\vp}$ and then taking expectations, we get
 \[\pr(|X^{\sss(m),1} - X_r^{\sss(m),1}|>\eps)\leq (\eps^2 r)^{-1}\E(\var_{\vp}(\dA_{\sss(m)}(V_1)/\sigma(\vp)))=: (\eps^2 r)^{-1} C_{\sss(m)}, \mbox{ say }, \]
  where $\var_{\vp}$ defined analogously to $\Ep$ is the conditional variance operator. Obviously \[\E\left(\var_{\vp}\left(\frac{\dA_{\sss(m)}(V_1)}{\sigma(\vp)}\right)\right)
  \leq \var\left(\frac{\dA_{\sss(m)}(V_1)}{\sigma(\vp)}\right)
  \leq \E\left(\left(\frac{\dA_{\sss(m)}(V_1)}{\sigma(\vp)}\right)^2\right).\]
 From the argument given below \eqref{eqn:sum-of-wts-active}, it follows that $\|\dA_{\sss(m)}\|_{\infty} \leq ||F^{\exec,\vp}||_{\infty}$. Hence Lemma \ref{lem:tail-fexec-inf} implies that $\sup_m\ C_{\sss(m)} < \infty$. This verifies (i) of the lemma.

  Let us now verify condition (ii) of the lemma. Writing this out explicitly,  we have to show for each fixed $r\geq 1$,
  \begin{equation}
  \label{eqn:to-show-lemma-ii}
	\left(\frac{\sum_{i=1}^r\dA_{\sss(m)}(V_i^{\sss(m)})}{r\sigma(\vp)}, \frac{\sum_{i=1}^r H_{\phi}^{\sss(m)}(i)}{r}\right)
  	\convd \left(\frac{\sum_{i=1}^r\dA_{\sss(\infty)}(V_i^{\sss(\infty)})}{r}, \frac{\sum_{i=1}^r H_{\phi}^{\sss(\infty)}(i)}{r}\right).
  \end{equation}
To this end, for each $m\in \set{1,2,\ldots}\cup \set{\infty}$, consider the subtree spanning the $J_r$ points $(V_i^{\sss(m)})_{1\leq i\leq r}, (\vV_{k,k+l}^{\sss(m),i})_{1\leq i\leq r}$, viewed as an element of $\vT_{IJ}^*$ as in Definition \ref{def:tree-rij-m}. Using Theorem \ref{thm:rij-convg} and continuity of the function $g_\phi^{(k)}$ from Theorem \ref{thm:g-phi}, we get
 \[\left(\left(\frac{\dA_{\sss(m)}(V_i)}{\sigma(\vp)}\right)_{1\leq i\leq r}, \left(H_{\phi}^{\sss(m)}(i)\right)_{1\leq i\leq r}\right)\convd \left(\left({\dA_{\sss(\infty)}(V_i)}\right)_{1\leq i\leq r}, \left(H_{\phi}^{\sss(\infty)}(i)\right)_{1\leq i\leq r}\right)\]
 with respect to weak convergence on $\bR^{2r}$, which in turn implies \eqref{eqn:to-show-lemma-ii}. This completes the verification of the conditions of Lemma \ref{lem:approx-m-inf} and thus the proof of Theorem \ref{thm:jt-convg-E}. \qed
\medskip

The rest of this section proves Theorems \ref{thm:rij-convg} and \ref{thm:g-phi}.

\medskip
\noindent{\bf Proof of Theorem \ref{thm:rij-convg}:} The proof will rely on a truncation argument that is qualitatively similar to Lemma \ref{lem:approx-m-inf}. Fix a truncation level $R\geq 1$.    Recall the definition of $\dA_{\sss(m)}(v)$ from \eqref{eqn:amv-def} which kept track of the contribution of {\bf all} right children of individuals $i$ on the path $[\rho,v]$. We will look at a truncated version of this object where we keep track of the potential contributions of only the first $R$ vertices. More precisely let
   \begin{equation}
   \label{eqn:I-amv-def}
   	\dA_{\sss(m)}^{\sss R}(v):= \sum_{\sss\substack{i\in [\rho,v]\\i\leq R}} \sum_{j\in [m]} p_j \ind\set{j\in \RC(i,[\rho,v])}.
   \end{equation}
   Let $\dA_{\sss(\infty)}^R(\cdot)$ be the analogous modification of $\dA_{\sss(\infty)}(\cdot)$ defined in \eqref{eqn:limit-da}, i.e,
   \begin{equation}
   \label{eqn:I-amv-def-icrt}
   \dA_{\sss(\infty)}^R(v)=\sum_{i\leq R}\theta_{i}\left[\sum_{j\geq 1}U_j^{\sss(i)}\times\ind\set{v\in \cT_j^{\sss(i)}}\right].
   \end{equation}
    Similarly modify the ``second endpoint'' measure in \eqref{eqn:right-end-pt-prob} to keep track of only ancestors with labels $\leq R$, namely
	
	\begin{equation}
	\label{eqn:right-end-pt-prob-trunc}
		Q_{v}^{\sss(m), R}(y):= \begin{cases}
			\sum_{u} p_u \ind\set{u\in \RC(y,[\rho,v])}/\dA_{\sss(m)}^{\sss R}(v), & \text{ if } y\in [\rho,v] \text{ and } y\leq R,\\
			0, & \text{ otherwise}.
		\end{cases}
	\end{equation}

Note that this does not make sense if $\dA_{\sss(m)}^{\sss R}(v) = 0$, i.e., when there is {\bf no} vertex with label $\leq R$ on the path from the root to $v$. In this case we follow the convention of defining the measure to be the uniform probability measure on the line $[\rho,v]$. Define $Q_{v}^{\sss(\infty), R}(\cdot)$ on $\icrt$ in an analogous fashion.

Consider the tree $r_{IJ}^{\sss(m)}$ as in Definition \ref{def:tree-rij-m}, and assign to leaf $V_j$ the truncated measure $Q_{V_j}^{\sss(m), R}(\cdot)$ and leaf value $\dA_{\sss(m)}^{\sss R}(V_j)$ (instead of $Q_{V_j}^{\sss(m)}(\cdot)$ and $\dA_{\sss(m)}(V_j)/\sigma(\vp)$). We denote the resulting object (which is an element of $\vT_{IJ}^*$) by $\cR_{IJ}^{\sss(m),R}$. Similarly construct $\cR_{IJ}^{\sss(\infty),R}$.
\begin{prop}
	\label{prop:trunc-R-rij}
	The following hold:
	\begin{enumeratea}
		\item \label{prop:i} For all $R\geq 1$, $\cR_{IJ}^{\sss(m), R} \convd \cR_{IJ}^{\sss(\infty), R}$.
		\item \label{prop:iii} $\cR_{IJ}^{\sss(\infty), R} \convd \cR_{IJ}^{\sss(\infty)}$ as $R\to\infty$.
		\item \label{prop:ii} 
		
		 For any bounded continuous function $f:\vT_{IJ}^* \to \bR$,
		\[\limsup_{R\to\infty}\limsup_{m\to\infty}\left|\E(f(\cR_{IJ}^{\sss(m), R})) - \E(f(\cR_{IJ}^{\sss(m)}))\right|   = 0. \]
		
	\end{enumeratea}
\end{prop}

Assuming this proposition, we now complete the proof of Theorem \ref{thm:rij-convg}. Note that for any fixed bounded continuous function $f$ on $\vT_{IJ}^*$ and any truncation level $R\geq 1$, we have
\begin{align*}
	|\E(f(\cR_{IJ}^{\sss(\infty)})) - \E(f(\cR_{IJ}^{\sss(m)}))|\leq |\E(f(\cR_{IJ}^{\sss(\infty)})) - \E(f(\cR_{IJ}^{\sss(\infty),R}))|+|&\E(f(\cR_{IJ}^{\sss(\infty),R})) - \E(f(\cR_{IJ}^{\sss(m),R}))|\\
	&+|\E(f(\cR_{IJ}^{\sss(m),R})) - \E(f(\cR_{IJ}^{\sss(m)}))|.
\end{align*}
Now letting $m\to\infty$ and then letting $R\to \infty$ and using Proposition \ref{prop:trunc-R-rij} completes the proof.\qed

We next prove Proposition \ref{prop:trunc-R-rij}.
\subsection{Proof of Proposition \ref{prop:trunc-R-rij}.} We start with three preliminary lemmas. Recall that $\set{i\leadsto j}$ denotes the event that $j$ is a child of $i$ in $\cT_m^{\vp}$.

\begin{lem}
\label{lem:hub-weight-conc}
Under Assumption \ref{ass:vp-a}, for each fixed $i\geq 1$,
\[\left|\frac{\sum_{j\in [m]}p_j \ind\set{i\leadsto j}}{\sigma(\vp)} - \frac{p_i}{\sigma(\vp)}\right|\probc 0 \qquad \mbox{ as } m\to\infty. \]	 
\end{lem}

\noindent{\bf Proof:}  Recall from \eqref{eqn:p-tree-connection} that for fixed $i$, the collection of events $\set{\set{i\leadsto j}: j\neq i}$ are pairwise independent and have the same probability $p_i$. Thus
 \begin{align*}
 \E\left(\frac{\sum_{j\in [m]}p_j \ind\set{i\leadsto j}}{\sigma(\vp)}\right) = \frac{p_i}{\sigma(\vp)}\left(\sum_{j\neq i} p_j\right) = (1-p_i)\frac{p_i}{\sigma(\vp)},
 \end{align*}
and
\begin{align*}
\var\left(\frac{\sum_{j\in [m]}p_j \ind\set{i\leadsto j}}{\sigma(\vp)}\right) = \sum_{j\in [m]} \frac{p_j^2}{\sigma^2(\vp)}\var(\ind\set{i\leadsto j})\leq p_i.
\end{align*}
This completes the proof as \ch{$\max_{i\in [m]}p_i = p_1\to 0$} and $p_i/\sigma(\vp)\to\theta_i$ under Assumption \ref{ass:vp-a}. \qed

\begin{lemma}\label{lem:maxpj-to-i}
	Under Assumption \ref{ass:vp-a}, for each fixed $i\geq 1$,
	\[\max_{j: i\leadsto j} \frac{p_j}{\sigma(\vp)} \probc 0, \qquad \mbox{ as } m\to\infty.\]
\end{lemma}
\noindent{\bf Proof:} Fix $\eps>0$ and write
\[\cN_{\eps}(m):= \set{j: p_j\geq \sigma(\vp)\eps}\ \text{ and }\ n_{\eps}(m)=|\cN_{\eps}(m)|.  \]
Note that by Assumption \ref{ass:vp-a}, for every $\eps> 0$, $\set{n_\eps(m):m\geq 1}$ is a bounded sequence.  Further, \eqref{eqn:p-tree-connection} and Markov's inequality yield
\[\pr\left(\max_{j: i\leadsto j} \frac{p_j}{\sigma(\vp)}> \eps\right) \leq \sum_{j\in \cN_{\eps}(m)} p_i = n_\eps{(m)}p_i\to 0,\]
as \ch{$\max_{i\in [m]}p_i =p_1\to 0$}. \qed

\medskip
Recall that $\cD_m(i)$ is the set of children of vertex $i$ in $\cT_m^{\vp}$. For later use let $d_m(i):=|\cD_m(i)|$ denote the degree of $i$ in $\cT_m^{\vp}$. Note that Lemma \ref{lem:hub-weight-conc} together with the lemma just proven gives
\begin{equation}
\label{eqn:di-to-infty}
	d_m(i)\probc \infty, \qquad \mbox{ as }m\to\infty.
\end{equation}
\begin{lem}\label{lem:sup-q-excursion}
	For each fixed $m$, let $\vq(m):= (q_1,q_2,\ldots q_d)$ be a probability mass function with $q_i > 0$ for all $i,m$ where $d = d(m)\to \infty$ as $m\uparrow \infty$. Assume further that $q_{\max}:=\max_{i\in [d]} q_i\to 0 $ as $m\to\infty$. Let $\set{U_i^{\sss(m)}:1\leq i\leq d}$ be i.i.d. uniform random variables and consider the function
	\[W_m(t):= \sum_{i=1}^d q_i \ind\set{U_i^{\sss(m)} \leq t} - t, \qquad t\in [0,1].\]
Then $\sup_{t\in [0,1]} |W_m(t)| \probc 0$ as $m\to\infty$.
\end{lem}
\noindent{\bf Proof:} Recall the proof of Lemma \ref{lem:tail-fexec-inf} where we studied the tightness of the tilt. Then replacing $\vp$ in the proof by $\vq$, the quantity of interest is $\sup_{t\in [0,1]} |W_m(t)| = \sigma(\vq)\cR_1(m)$ where $\cR_1(m)$ is as defined in \eqref{eqn:1013} and $\sigma(\vq):=\sqrt{\sum_i q_i^2}$. Now \eqref{eqn:nts-r1} and \eqref{eqn:tail-r1} imply the existence of a constant $C$ (independent of $m$) such that for all $m$ and $x\geq e$,
\begin{align*}\pr(\sup_{t\in [0,1]} |W_m(t)| > x\sigma(\vq))\leq \exp(-C x\log(\log{x})).
\end{align*}
Since $\sigma(\vq)\leq \sqrt{q_{\max}} \to 0$ as $m\to \infty$, this completes the proof. \qed

\medskip

We now have all the ingredients for the proof of Proposition \ref{prop:trunc-R-rij}. We prove parts (a), (b) and (c) one by one.

\noindent{\bf Proof of Proposition \ref{prop:trunc-R-rij} \eqref{prop:i}:} Recall from Definition \ref{def:tree-rij-m} the tree $r_{IJ}^{\sss(m)}$ that contains all the edge lengths and hub information in $\cR_{IJ}^{\sss(m)}$ but ignores root-to-leaf measures and lead values $\dA_{\sss(m)}(\cdot)$.  By \cite[Corollary 15]{pitman-camarri} or \cite[Proposition 3]{aldous-pitman-entrance}, for fixed $J\geq 1$, we have
\begin{equation}
\label{eqn:rij-weak-prod}
	\left(r_{I^{\prime}J}^{\sss(m)}: I^{\prime}\geq 0\right) \convd \left(r_{I^{\prime} J}^{\sss(\infty)}: I^{\prime}\geq 0\right)
\end{equation}
with respect to the product topology on $\prod_{I^\prime \geq 0} \vT_{I^{\prime} J}$.
Using Lemma \ref{lem:hub-weight-conc}, Lemma \ref{lem:maxpj-to-i} and Skorohod embedding, we assume that we are working on a probability space that supports a sequence of unordered $\vp$-trees  $\set{\cT_m^{\vp,\unordered}:m\geq 1}$, sampled vertices $\set{V_j^{\sss(m)}:1\leq j\leq J, m\geq 1}$ using the associated sequence of probability mass functions $\set{\vp(m):m\geq 1}$, an ICRT $\icrt$, and sampled vertices $\set{V_j^{\sss(\infty)}:1\leq j\leq J}$ using the mass measure such that the following hold:
\begin{enumerateA}
	\item Convergence in \eqref{eqn:rij-weak-prod} happens almost surely:
	\begin{equation}
	\label{eqn:rij-as}
		\left(r_{I^{\prime}J}^{\sss(m)}: I^{\prime}\geq 0, \right) \convas \left(r_{I^{\prime} J}^{\sss(\infty)}: I^{\prime}\geq 0\right)\text{ as }m\to\infty
	\end{equation}
coordinatewise, where the underlying tree corresponding to $r_{I^{\prime}J}^{\sss(m)}$ is spanned by the root of $\cT_m^{\vp,\unordered}$ and $V_j^{\sss(m)},\ 1\leq j\leq J$.
\item Writing $s_m(i):= \sum_{v\in \cD_m(i)} p_v$ for the sum of weights of children of $i$ in $\cT_m^{\vp,\unordered}$, we have
\begin{align}\label{eqn:convergence-to-theta-i}
\left(\frac{s_m(i)}{\sigma(\vp)}\colon i\geq 1\right) \convas \left(\theta_i\colon i\geq 1\right)
\end{align}
coordinatewise. (We can assume that this holds because of Lemma \ref{lem:hub-weight-conc}.)
	\item For fixed hub $i\geq 1$ and $m\geq 1$, write
 \begin{equation}
	\label{eqn:qm-i}
		 q_{m,i}(v):= \frac{p_v}{s_m(i)} , \qquad v\in \cD_m(i), \qquad q_{m,i}^{\max}:=\max_{v\in \cD_m(i)} q_{m,i}(v).
	\end{equation}
	 Then we assume (using Lemma \ref{lem:maxpj-to-i} and \eqref{eqn:di-to-infty}) that for all $i\geq 1$
\[q_{m,i}^{\max}\convas 0\text{ and }d_m(i)\convas\infty.\]
\end{enumerateA}

Now, for each $z\in [m]$ and $i\geq 1$, if $i\in [\rho,z]$ (where $\rho=\rho_m$ is the root of $\cT_m^{\vp,\unordered}$), write $c(i;z)\in \cD_m(i)$ for the child of $i$ that is the ancestor of $z$. Next, construct a collection $\set{U_{m,i}(v):m\geq 1,i\geq 1, v\in [m]}$ of uniform$[0, 1]$ random variables on the same space
such that
\begin{enumeratea}
 \item $\set{\cT_m^{\vp,\unordered}, U_{m,i}(v)\colon i\geq 1, v\in [m]}$ are jointly independent for each $m\geq 1$; and
 \item for each $i\leq R$ and $j\leq J$ for which $i\in [\rho, V_j^{\sss(\infty)}]$, $U_{m,i}\left(c(i;V_j^{\sss(m)})\right)$ is a constant sequence (in $m$) eventually.
 \end{enumeratea}
 As described below Theorem \ref{thm:jt-convg-E}, we can use these uniform random variables to generate the sequence of ordered $\vp$-trees $\set{\cT_m^{\vp}}$ from $\set{\cT_m^{\vp,\unordered}}$ as follows: Let $\mvU_{m,i}:=\set{U_{m,i}(v): v\in \cD_{m}(i)}$. Think of these as ``ages'' of the children and arrange the children from left to right in decreasing order of their ages.

Once this ordering has been defined, we can construct the function $\dA_{\sss(m)}(\cdot)$ as in \eqref{eqn:amv-def}. In this case we can write this function explicitly in terms of the associated uniform random variables as follows.  Define
\begin{align}
	\dO_{\sss(m),i}(z)&:= \frac{1}{\sigma(\vp)} \ind\set{i\in [\rho,z]}\sum_{v\in \cD_m(i)} p_v \ind\set{U_{m,i}(v) < U_{m,i}(c(i;z)) } \notag\\
	&= \ind\set{i\in [\rho,z]}\frac{s_m(i)}{\sigma(\vp)} \sum_{v\in \cD_m(i)} q_{m,i}(v) \ind\set{U_{m,i}(v) < U_{m,i}(c(i;z))}. \notag \label{eqn:dom-def}
\end{align}
Then
\begin{equation}
\label{eqn:dam-redef}
\left(\sigma(\vp)\right)^{-1}\dA_{\sss(m)}^{\sss R}(z)= \sum_{i\leq R} \dO_{\sss(m), i}(z).
\end{equation}
Similarly, the root-to-leaf measure $Q_{\sss v}^{\sss(m), R}$ (recall \eqref{eqn:I-amv-def}) can also be expressed in terms of this function.

Now using \eqref{eqn:rij-as}, for every fixed hub $i\leq R$, $j\leq J$, and a.e. sample point $\omega$, one of the following two holds:
\begin{enumeratea}
	\item $i\notin [\rho, V_j^{\sss(\infty)}]$, in which case there exists $m= m(\omega)$ such that $i\notin [\rho, V_j^{\sss(m)}]$ for all $m> m(\omega)$.
	\item $i\in [\rho, V_j^{\sss(\infty)}]$, in which case there exists $m= m(\omega)$ such that $i\in [\rho, V_j^{\sss(m)}]$ for all $m > m(\omega)$. 	
\end{enumeratea}
When the latter happens, using Lemma \ref{lem:sup-q-excursion} together with \eqref{eqn:convergence-to-theta-i} and \eqref{eqn:qm-i}, we get
	\[\bigg|\dO_{\sss(m), i}(V_j^{\sss(m)}) - \theta_i U_{m, i}\left(c(i;V_j^{\sss(m)})\right)\bigg|\probc 0 \quad \mbox{as } m\to\infty. \]
By construction, $U_{m, i}\left(c(i;V_j^{\sss(m)})\right)$ is eventually constant in $m$ on the event $\set{i\in [\rho, V_j^{\sss(\infty)}]}$. This immediately implies convergence of the (scaled) truncated leaf values $\dA_{\sss(m)}^{\sss R}(V_j^{\sss(m)})/\sigma(\vp)$ (see \eqref{eqn:dam-redef}) for $1\leq j\leq J$, and similarly the truncated root to leaf measures $Q_{\sss V_j^{\sss(m)}}^{\sss(m), R}$ jointly with the convergence in \eqref{eqn:rij-as} and thus yields the convergence $\cR_{IJ}^{\sss(m), R} \convd \cR_{IJ}^{\sss(\infty), R}$. \qed

\medskip
\noindent{\bf Proof of Proposition \ref{prop:trunc-R-rij} \eqref{prop:iii}:} Recall from Section \ref{sec:p-tree-ICRT-def} that $\cR_{IJ}^{\sss(\infty)}$ is obtained by applying the stick-breaking construction to $[0,\eta_J]$, and leaf $j+$ in $\cR_{IJ}^{\sss(\infty)}$ corresponds to the vertex coming from $\eta_j$.  It is easy to see from the definition of $\dA_{\sss(\infty)}^{\sss R}$ and $Q_{v}^{\sss(\infty),R}$ that it suffices to prove
\begin{equation}
\label{eqn:to-show-R-infty}
	0\leq \cE_R^{\sss(1)}:= \sum_{j=1}^J (\dA_{\sss (\infty)}(\eta_j) - \dA_{\sss (\infty)}^{\sss R}(\eta_j)) \probc 0, \quad \mbox{ as } R\to \infty.\notag
\end{equation}
For every hub $i\geq 1$ and leaf $\eta_j$, write $\set{i\to \eta_j}$ if $\eta_j$ is a descendant of $i$ (namely $i\in [\rho,\eta_j]$).  Then note that
\[\cE_R^{\sss(1)}=\sum_{j=1}^J  \sum_{i=R+1}^{\infty}  \sum_{k=1}^{\infty} \theta_i U_{k}^{\sss(i)} \ind\set{\eta_j\in \cT_k^{\sss(i)}} \leq \sum_{j=1}^J \sum_{i=R+1}^\infty \theta_i \ind\set{i\to \eta_j} := \cE_R^{\sss(2)}.
 \]
Thus, it is enough to show that given $\eps > 0$, we can find $R= R(\eps) < \infty$ such that $\pr(\cE_{R}^{\sss(2)} > \eps)<\eps$. To this end, first choose $K_{\eps}$ large enough so that $\pr(\eta_J > K_\eps ) < \eps/2$, and then choose $R_\eps$ large enough so that
\[\frac{J K_{\eps}}{\eps}\sum_{R_{\eps}+1}^\infty \theta_i^2 < \eps/2. \]
Then note that
\begin{align*}
	\pr(\cE_{R_\eps}^{\sss(2)} > \eps) &\leq \pr(\eta_J > K_\eps) + \pr\left(J\sum_{i=R_\eps+1}^\infty \theta_i \ind\set{i \mbox{th hub appears before time } K_\eps } > \eps\right)\\
	&\leq \frac{\eps}{2}+ \frac{J}{\eps}\sum_{i=R_{\eps}+1}^\infty \theta_i \left(1-\exp(-\theta_i K_\eps)\right)  \\
	&< \eps\qquad \text{ by the choice of $R_{\eps}$},
\end{align*}
where the first term in the second inequality follows from the choice of $K_\eps$, while the second term comes from the stick-breaking construction of $\icrt$ using the countable collection of Poisson point processes. This completes the proof.  \qed

\medskip
\noindent{\bf Proof of Proposition \ref{prop:trunc-R-rij}\eqref{prop:ii}:} Recall that the tree $\cR_{IJ}^{\sss(m)}$ (and $\cR_{IJ}^{\sss(m),R}$) can be thought of as being made up of $2J+1$ coordinates:
\begin{enumeratea}
	\item One coordinate for the shape and edge length information along with the labels smaller than $I$ namely $r_{IJ}^{\sss(m)}$ (see Definition \ref{def:tree-rij-m}). Note that this is the same for both $\cR_{IJ}^{\sss(m)}$ and $\cR_{IJ}^{\sss(m),R}$.
	\item $J$ coordinates for the leaf values $\dA_{\sss(m)}(V_j)/\sigma(\vp)$ (resp. $\dA_{\sss(m)}^{\sss R}(V_j)/\sigma(\vp)$).
	\item $J$ coordinates for the measured metric spaces $\cM_j^{\sss(m)}:= ([\rho, V_j^{\sss(m)}],\ Q_{\sss V_j}^{\sss(m)})$ (resp. $\cM_j^{\sss(m), R}:= ([\rho, V_j^{\sss(m)}],\ Q_{\sss V_j}^{\sss(m), R})$).
\end{enumeratea}
Since $\vT_{IJ}^*$ assumes the product topology on these coordinates, it is enough to show the required estimate in Proposition \ref{prop:trunc-R-rij} \eqref{prop:ii} with functions of the form
\begin{equation}
\label{eqn:f-form}
	f(\vt, (a_i)_{1\leq j\leq J}, (M_j)_{1\leq j\leq J}):= F(\vt) \prod_{1\leq j\leq J} g_j(a_i) \prod_{1\leq j\leq J} h_j(M_j).\notag
\end{equation}
Here $\vt\in \vT_{IJ}$, $a_j\in \bR$ are associated leaf values and $M_j$ are the paths from the root to leaf $j$ with an associated probability measure and $f$, $g_j$ and $h_j$ are bounded \emph{uniformly continuous} functions on the spaces $\vT_{IJ}$, $\bR$ and $\sS$ (measured compact metric spaces) respectively.  To simplify notation, we will simply write this as $f(\vt)$.

Now we can go from $\cR_{IJ}^{\sss(m)}$ to $\cR_{IJ}^{\sss(m), R}$ by flipping one coordinate at a time. Thus writing
\[f^{(-i)}_1(\vt):= F(\vt) \prod_{\substack{1\leq j\leq J\\j\neq i}} g_j(a_i) \prod_{1\leq j\leq J} h_j(M_j), \qquad f^{(-i)}_2(\vt):= F(\vt) \prod_{1\leq j\leq J} g_j(a_i) \prod_{\substack{1\leq j\leq J\\j\neq i}} h_j(M_j), \]
we get
\begin{align}
	|\E(\ch{f(\cR_{IJ}^{\sss(m)})})-\E(f(\cR_{IJ}^{\sss(m), R}))|\leq \sum_{j=1}^J||f^{(-j)}_1||_{\infty}& \E\left(\left|g_j\left(\frac{\dA_{\sss (m)}(V_j)}{\sigma(\vp)}\right) - g_j\left(\frac{\dA_{\sss (m)}^{\sss R}(V_j)}{\sigma(\vp)}\right) \right|\right) \notag \\
	&+\sum_{j=1}^J||f^{(-j)}_2||_{\infty} \E(|h_j(\cM_j^{\sss(m)}) - h_j(\cM_j^{\sss(m), R}) |) \label{eqn:lindeberg}
\end{align}
Since $V_j$'s have been sampled in an i.i.d. fashion from $\vp$, it is enough to show that for any two bounded uniformly continuous functions  $h, g$ on $\bR$ and $\sS$ respectively,
\begin{equation}
\label{eqn:g-go-zero}
	\limsup_{R\to\infty}\limsup_{m\to\infty}\ \E\left(\left|g\left(\frac{\dA_{\sss (m)}(V_1^{\sss(m)})}{\sigma(\vp)}\right) - g\left(\frac{\dA_{\sss (m)}^{\sss R}(V_1^{\sss(m)})}{\sigma(\vp)}\right) \right|\right)  = 0,
\end{equation}
and
\begin{equation}
\label{eqn:h-go-zero}
	\limsup_{R\to\infty}\limsup_{m\to\infty}\ \E\left(|h(\cM_1^{\sss(m)}) - h(\cM_1^{\sss(m), R}) |\right)  = 0.
\end{equation}
Now consider the measured metric spaces $\cM_1^{\sss(m)} $ and $\cM_1^{\sss(m),R}$. As remarked above, they share the same metric space, namely the path $[\rho, V_1^{\sss(m)}]$.  The only difference is in the associated probability measures. Consider the natural correspondence $C=\set{(x,x): x\in [\rho,V_1^{\sss(m)}]}$ between $\cM_1^{\sss(m)} $ and $\cM_1^{\sss(m),R}$. Further, define a probability measure $\pi$ on $[\rho,V_1^{\sss(m)}] \times [\rho, V_1^{\sss(m)}]$ as
\begin{equation}
\label{eqn:pi-joint}
	\pi(i,i):= \begin{cases}
		{\sum_{u} p_u \ind\set{u\in \RC(i,[\rho,V_1^{\sss(m)}])}}/{\dA_{\sss(m)}(V_1^{\sss(m)})}, & \text{ if }i\in \ch{(\rho,V_1^{\sss(m)}]} \text{ and } i\leq R,\\
			\left[{\dA_{\sss(m)}(V_1^{\sss(m)}) - \dA_{\sss(m)}^{\sss R}(V_1^{\sss(m)}) }\right]/{\dA_{\sss(m)}(V_1^{\sss(m)})}, & \text{ if } i=\rho.
	\end{cases}
\end{equation}

Writing $\pi_1$ and $\pi_2$ for the marginals of $\pi$, we have, using the above choice of correspondence $C$ and of the measure $\pi$,
\begin{equation}
\label{eqn:qm-pi-tv}
	d_{\GHP}^{\point}(\cM_1^{\sss(m)}, \cM_1^{\sss(m), R})\leq \left(|| \pi_1 - Q_{\sss V_1}^{\sss(m)} ||+||\pi_2 - Q_{\sss V_1}^{\sss(m), R}  || \right)\leq \frac{2\left[{\dA_{\sss(m)}(V_1) - \dA_{\sss(m)}^{\sss R}(V_1) }\right]}{\dA_{\sss(m)}(V_1)}.
\end{equation}
Now suppose we show \eqref{eqn:g-go-zero}. Using part \eqref{prop:i} and part \eqref{prop:iii} of Proposition \ref{prop:trunc-R-rij}, we get $(\sigma(\vp))^{-1}\dA_{\sss(m)}(V_1^{\sss(m)}) \convd \dA_{\sss(\infty)}(V_1^{\sss(\infty)}) > 0$. Now using the bound in \eqref{eqn:qm-pi-tv} and uniform continuity of $h$, we see that \eqref{eqn:h-go-zero} is true. Hence it is enough to prove \eqref{eqn:g-go-zero}.

 Recall from Section \ref{sec:ptree-const} the construction of $V_1^{\sss(m)}$ and the tree simultaneously via the birthday construction, where $V_1^{\sss(m)}$ is obtained as the value before the first repeat time, namely $Y_{R_1-1}$.
    Fix $\eps>0$. By \cite[Theorem 4]{pitman-camarri}, under Assumptions \ref{ass:vp-a} we may choose $K_\eps$ large so that the first repeat time satisfies  $\pr(R_1 > K_\eps/\sigma(\vp)) < \eps$ for all $m\geq 1$. Next, by uniform continuity of $g$, choose $\delta\in (0,1)$ such that $|g(x) -g(y)| < \eps$ if $|x-y| < \delta$. Finally choose  $R$ large so that for all $m$,
   \[\frac{K_\eps^2}{\delta \wedge \eps} \sum_{i=R+1}^m \frac{p_i^2}{\sigma^2(\vp)} < \eps. \]
      First, by choice of $K_\eps$ and boundedness of $g$,
   \begin{equation}
   \label{eqn:am-am-rr}
\left|\E\left[g\left(\frac{\dA_{\sss(m)}(V_1^{\sss(m)})}{\sigma(\vp)}
\right)\right]-\E\left[g\left(\frac{\dA_{\sss(m)}(V_1^{\sss(m)})}{\sigma(\vp)}
\right)\ind\set{R_1 \leq \frac{K_\eps}{\sigma(\vp)}}\right]\right| \leq ||g||_{\infty} \eps,
   \end{equation}
   and a similar inequality holds true if we replace the functional $\dA_{\sss(m)}$ by $\dA_{\sss(m)}^R$.
   Next, writing
   \begin{equation}
   \label{eqn:def-eps-I-am}
	\cE_{\sss(m)}^{\sss(1)}(R):= \left|\E\left[\set{g\left(\frac{\dA_{\sss(m)}(V_1^{\sss(m)})}{\sigma(\vp)}
	\right)-g\left(\frac{\dA_{\sss(m)}^{\sss R}(V_1^{\sss(m)})}{\sigma(\vp)}
	\right)}\ind\set{R_1 \leq \frac{K_\eps}{\sigma(\vp)}}\right]\right|,\notag
   \end{equation}
   we have
   \begin{equation}
   \label{eqn:ami-diff-bound}
   	\cE_{\sss(m)}^{\sss(1)}(R) \leq \eps + 2||g||_{\infty} \pr\left(R_1 \leq \frac{K_\eps}{\sigma(\vp)}~,~ \frac{ \left(\dA_{\sss(m)}(V_1^{\sss(m)}) - \dA_{\sss(m)}^R(V_1^{\sss(m)})\right)}{\sigma(\vp)} \geq \delta\right)
   \end{equation}
   by our choice of $\delta$. The difference $\dA_{\sss(m)}(\ch{V_1^{\sss(m)}}) - \dA_{\sss(m)}^{\sss R}(\ch{V_1^{\sss(m)}})$ is a tricky object for which we will need a tractable upper bound. Recall that we have used $\cT_1^{\cB}$ for the birthday tree in \eqref{eqn:p-tree-joint-sample} constructed by time $R_1$. For every vertex $i\in \cT_1^{\cB}$, let $\cJ(i)$ be the first child of $i$ in the birthday construction (the first {\bf new, i.e., previously un-sampled} vertex sampled immediately after a prior sampling of $i$). This will be an empty set if $i$ is a leaf in the eventual full tree $\cT_m^{\vp}$. Recall that $\set{i\leadsto j}$ was used to denote the event that $j$ is a child of $i$ in $\cT_m^{\vp}$. Then note that
  \begin{equation}
  \label{eqn:am-aim-sup-bound}
  	\dA_{\sss(m)}(V_1^{\sss(m)}) - \dA_{\sss(m)}^{\sss R}(V_1^{\sss(m)}) \leq \sum_{i\geq R+1}\ind\set{i\in \cT_1^{\cB}}\sum_{j\in [m]}p_j\ind\set{i\leadsto j, j\neq \cJ(i)}.\notag
  \end{equation}
  Thus,
  \begin{align}
  	&\pr\Bigg(R_1 \leq \frac{K_\eps}{\sigma(\vp)}~,~ \frac{ \left(\dA_{\sss(m)}(V_1^{\sss(m)}) - \dA_{\sss(m)}^{\sss R}(V_1^{\sss(m)})\right)}{\sigma(\vp)} \geq \delta\Bigg)\notag\\
	&\hskip20pt\leq \frac{1}{\delta}\sum_{i=R+1}^m \sum_{j\in [m]}\frac{p_j}{\sigma(\vp)}\pr\left( i \mbox{ appears before }\frac{K_\eps}{\sigma(\vp)}, i\leadsto j, j\neq \cJ(i)\right)=:\cE_{\sss (m)}^{\sss (2)}(R).\label{eqn:1358}
  \end{align}
  For $i\neq j\in [m]$, define the event $E_{ij}:=\set{i \mbox{ appears before }\frac{K_\eps}{\sigma(\vp)}, i\leadsto j, j\neq \cJ(i)} $. Then for $E_{ij}$ to happen, the following needs to happen in the birthday construction: (a) There is an $0\leq r_1\leq K_\eps/\sigma(\vp)$ such that till time $r_1$, neither $i$ or $j$ have been sampled. (b) At time $r_1+1$ vertex $i$ is sampled. (c) There is an $r_2\geq 0$ such that in the times $[r_1+1, r_1+1+r_2]$ samples, $j$ does not appear. (d) Then at time $r_1+r_2+2$, vertex $i$ is sampled again. (e) In the next time step $r_1+r_2+3$ vertex $j$ is sampled. Therefore,
  \begin{equation}
  \label{eqn:eij-bound}
  	\pr(E_{ij})\leq \sum_{r_1=0}^{K_\eps/\sigma(\vp)} \sum_{r_2=0}^{\infty} (1-p_i-p_j)^{r_1} p_i(1-p_j)^{r_2} p_i p_j \leq p_i^2\frac{K_\eps}{\sigma(\vp)}.\notag
  \end{equation}
  Using this in \eqref{eqn:1358}, we get
  \begin{equation}
  \label{eqn:1420}
  	\cE_{\sss(m)}^{\sss (2)}(R)\leq \frac{K_\eps}{\delta}\sum_{R+1}^m \frac{p_i^2}{[\sigma(\vp)]^2} \leq \eps,\ \mbox{ by our choice of } R.
  \end{equation}
  Combining \eqref{eqn:am-am-rr}, \eqref{eqn:ami-diff-bound}, \eqref{eqn:1358} and \eqref{eqn:1420} now gives the following lemma which completes the proof of \eqref{eqn:g-go-zero} and thus the proof of part \eqref{prop:ii} of the proposition. \qed

  \begin{lemma}\label{lem:approx-am-ami}
  	Given $\eps> 0$ choose $K_\eps, \delta$ and $R$ as above. Then, for all $m\geq 1$,
\begin{align*}
\left|\E\left[g\left(\frac{\dA_{\sss(m)}(V_1^{\sss(m)})}{\sigma(\vp)}
\right)\right]- \E\left[g\left(\frac{\dA_{\sss(m)}^{\sss R}(V_1^{\sss(m)})}{\sigma(\vp)}
\right)\right] \right|\leq \eps(4||g||_{\infty}+1).
\end{align*}
  \end{lemma}

\noindent{\bf Proof of Theorem \ref{thm:g-phi}:} We now prove continuity of the function $g_{\phi}^{(k)}$ on the space $\vT_{I,(k+\ell)}^*$. In fact, we will give a quantitative estimate. Since we are assuming the discrete topology on the coordinate corresponding to the shape, without loss of generality we will work with two trees $\vt, \overline\vt\in \vT_{I,(k+\ell)}^*$ having the same shape. We need to distinguish the labels for the root and the leaves in the two trees; so write $0+$ (respectively $\overline{0}+$) for the root of $\vt$ (respectively $\overline\vt$) and write $\set{j+: 1\leq j\leq k+\ell}$ (respectively $\set{\overline{j}+: 1\leq j\leq k+l}$) for the collection of leaves in $\vt$ (respectively $\overline\vt$). Finally, let $\nu_j$ be the corresponding probability measure on the path $\cM_j:= [0+, j+]$ for $1\leq j\leq k$, and analogously let $\overline\nu_j$ be the probability measure on $\overline\cM_j:= [\overline{0}+, \overline{j}+]$. View these paths as pointed measured metric spaces pointed at the roots $0+$ and $\overline0+$ respectively. Now let $\eps_j:= d_{\GHP}^{\point}(\cM_j,\overline\cM_j)$, where $d_{\GHP}^{\point}$ is the pointed Gromov-Hausdorff-Prokhorov metric defined in Section \ref{sec:metric-convg}.

Write $L = {\ell \choose 2}$. Let $\phi:\bR_+^L\to \bR$ be a bounded continuous function. For $K >0$, let $\square(K) = [0,K]^L$, and for $\delta >0$, define
\begin{equation}
\label{eqn:osc-def}
	\osc_{\phi}(\delta, K):= \sup_{\substack{\vx,\vy\in \square(K)\\||\vx-\vy||_{\infty} < \delta} } |\phi(\vx) - \phi(\vy)|.\notag
\end{equation}
 Finally, define
 \begin{equation}
 \label{eqn:eps-tone-ttwo}
 	\eps:= 4\sum_{j=1}^k \eps_j + (k+1)\sum_{e} \big|l_e(\vt) - l_e(\overline\vt)\big|,
 \end{equation}
  where $l_e(\cdot)$ denotes the length of the edge $e$ and we have used the fact that both trees have the same shape. Write $\height(\vt)$ for the height of tree $\vt$ (not graph distance, rather in terms of maximal distance from the root when incorporating edge lengths).
  \nomenclature[height]{$\height(\vt)$}{Height of a tree $\vt$ with edge lengths incorporated into the distance. }
   The following proposition completes the proof of Theorem \ref{thm:g-phi}:
  \begin{prop}
	  \label{prop:gph-osc}
  	For two trees $\vt, \overline\vt\in \vT_{I,(k+\ell)}^*$ having the same shape, and with $\eps$ as in \eqref{eqn:eps-tone-ttwo},
	\[|g_{\phi}^{\sss(k)}(\vt) - g_{\phi}^{\sss (k)}(\overline\vt)|\leq 2\eps||\phi||_{\infty} + \osc_{\phi}\bigg(\eps\;,\; 2\height(\vt)+2\height(\overline\vt)\bigg).   \]
  \end{prop}
\noindent{\bf Proof:} For each $j\leq k$, choose a correspondence $C_j$ and a measure $\pi_j$ on the product space $[0+, j+]\times [\overline{0}+, \overline{j}+]$ such that the following conditions are met:
\begin{inparaenuma}
	\item \label{cond:i} $(0+,\overline{0}+)\in C_j$;
	\item \label{cond:ii} the distortion satisfies $\dis(C_j) <3\eps_j$;
	\item \label{cond:iii} the measure of the complement satisfies $\pi_j(C_j^c)< 2\eps_j$;
	\item \label{cond:iv} and finally
	\begin{equation}
	\label{eqn:tv-bd}
		||\nu_j - p_\ast\pi_j|| + ||\overline\nu_j - \overline p_\ast\pi_j|| < 2\eps_j,
	\end{equation}
\end{inparaenuma}
where $p_\ast\pi_j$ and $\overline p_\ast\pi_j$ are the marginals of $\pi_j$.
  Now sample $(X_j^{\star}, \overline X_j^{\star})\sim \pi_j$ from $[0+, j+]\times [\overline{0}+, \overline{j}+]$ independently for $1\leq j\leq k$. By \eqref{eqn:tv-bd}, we can couple $(X_j^{\star}, \overline X_j^{\star})$ with two random variables $X_j, \overline X_j$ (again independently for $1\leq j\leq k$) such that $X_j\sim \nu_j$ and $\overline X_j\sim \overline\nu_j$, and further
  \begin{equation}
  \label{eqn:xi-xi-st-neq}
  	\pr\left(X_j\neq X_j^{\star}\right)+\pr\left(\overline X_j\neq \overline X_j^{\star}\right)< 2\eps_j.
  \end{equation}
  Using conditions \eqref{cond:ii} and \eqref{cond:iii}, we get
  \begin{equation}
  \label{eqn:dt1-dt2}
  	\pr\left(\left|d_{\vt}\left(0+,X_j^{\star}\right) - d_{\overline\vt}\left(\overline0+,\overline X_j^{\star}\right)\right|> 3\eps_j\right) \leq 2\eps_j,
  \end{equation}
  where $d_{\vt}$ is the distance metric on tree $\vt$ which incorporates the edge lengths.
  \nomenclature[dt]{$d_{\vt}$}{Distance metric on tree $\vt$ which incorporates the edge lengths.}
    Now write $E$ for the following ``good event'':
  \begin{equation}
  \label{eqn:good-event}
  	E:=\bigcap_{j=1}^k \set{X_j = X_j^{\star},\; \overline X_j= \overline X_j^{\star}, \; \left|d_{\vt}\left(0+,X_j^{\star}\right) - d_{\overline\vt}\left(\overline0+,\overline X_j^{\star}\right)\right|\leq 3\eps_j }.\notag
  \end{equation}
  It follows from \eqref{eqn:xi-xi-st-neq} and \eqref{eqn:dt1-dt2} that
   \begin{equation}\label{eqn:P-E-bound}
  	\pr(E^c)\leq 4\sum_{j=1}^k\eps_j.
  \end{equation}

  Now we are going to create ``shortcuts" by gluing the leaves to the corresponding sampled points. Let $\vS$ (resp. $\overline\vS$) be the (random) metric space obtained by identifying each of the leaves $j+$ (resp. $\overline j+$) with $X_j$ (resp. $\overline X_j$) in $\vt$ (resp. $\overline\vt$) for $1\leq j\leq k$ and write $d_{\vS}$ (resp. $d_{\overline\vS}$) for the induced metric. Then by definition,
  \begin{align*}
  g_{\phi}^{\sss(k)}(\vt)= \E\left[\phi\bigg(d_{\vS}\left((k+i_1)+,\ {(k+i_2)}+ \right): 1\leq i_1< i_2\leq \ell\bigg)\right],
  \end{align*}
  and an analogous expression holds for $g_{\phi}^{\sss(k)}(\overline\vt)$. This gives
  \begin{align}
  \label{eqn:gphi-norm-bound}
&\big|g_{\phi}^{\sss(k)}(\vt)-g_{\phi}^{\sss(k)}(\overline\vt)\big|\leq \E\bigg(\bigg|\phi\big(d_{\vS}\left((k+i_1)+,\ {(k+i_2)}+ \right): 1\leq i_1< i_2\leq \ell\big)\notag\\
&\hskip120pt-\phi\left(d_{\overline\vS}\left(\overline{(k+i_1)}+,\ {\overline{(k+i_2)}}+ \right): 1\leq i_1<i_2\leq \ell\right)\bigg|\bigg).
  \end{align}

Consider the map from $\vt$ to $\overline\vt$ which takes every vertex to the corresponding vertex and points on each edge are mapped by linear interpolation (using the edge lengths) to points on the corresponding edge. Consider $a\in [0, j+]$ and let
$\overline a\in [\overline0, \overline j+]$ be the corresponding point in $\overline\vt$ for some $j\leq k$. Then note that
\begin{align}\label{eqn:a-a-bar}
\bigg|d_{\vt}\left(a, X_j\right)-d_{\overline\vt}\left(\overline a, \overline X_j\right)\bigg|
&\leq \bigg|d_{\vt}\left(0+, X_j\right)-d_{\overline\vt}\left(\overline0+, \overline X_j\right)\bigg|
+\bigg|d_{\vt}\left(0+, a\right)-d_{\overline\vt}\left(\overline0+, \overline a\right)\bigg|\notag\\
&\leq 3\eps_j+\sum_e |l_e(\vt)-l_{e}(\overline\vt)|
\end{align}
on the set $E$.

Now consider a shortest path in $\vS$ connecting $(k+i_1)+$ and ${(k+i_2)}+$. We can go from $\overline{(k+i_1)}+$ to $\overline{(k+i_2)}+$ by taking the same route in $\overline{\vS}$, i.e., by traversing the same edges and taking the same shortcuts in the same order. We make the following observations:
\begin{inparaenumi}
\item The difference between distance traversed while crossing the edge $e$ is $|l_e(\vt)-l_e(\overline\vt)|$.
\item By \eqref{eqn:a-a-bar}, on the set $E$, taking a ``shortcut" contributes at most $(3\eps_j+\sum_e |l_e(\vt)-l_{e}(\overline\vt)|)$ to the difference between distance traversed.
\end{inparaenumi}
Since we have to take at most $k$ shortcuts, we immediately get
\begin{align}\label{eqn:difference-distance-traversed}
d_{\overline\vS}\left(\overline{(k+i_1)}+,\ \overline{(k+i_2)}+\right)
\leq d_{\vS}\left({(k+i_1)}+,\ {(k+i_2)}+\right)+3\sum_{j=1}^k\eps_j+(k+1)\sum_e |l_e(\vt)-l_{e}(\overline\vt)|\notag
\end{align}
on the set $E$. By symmetry, a similar inequality holds if we interchange the roles of $\vS$ and $\overline\vS$. This observation combined with \eqref{eqn:P-E-bound} and \eqref{eqn:gphi-norm-bound} yields the result. \qed

\section{Proofs: Convergence in Gromov-weak topology}
\label{sec:proof-convg-gromov}
Recall from Proposition \ref{prop:generate-nr-given-partition} that conditional on the partition of the vertices $\set{\cV^{\sss(i)}\colon i\geq 1}$ into the connected components, the actual structure of the components of $\cG(\vx, t)$ can be generated independently as the connected graph $\tilde{\cG}_{|\cV^{\sss(i)}|}(a_n^{\sss(i)},\vp_n^{\sss(i)})$ where $a_n^{\sss(i)}, \vp_n^{\sss(i)}$ are as in Proposition \ref{prop:generate-nr-given-partition} and given $m, \vp,a$, $\tilde{\cG}_m(a,\vp)$ is the connected random graph model studied in the previous section. For Theorem \ref{thm:mc-main-1}, the time scale $t = t_n$ of interest in the expression of $a_n^{\sss(i)}$ is
\begin{equation}
\label{eqn:tn-def}
	t_n:= \lambda + \frac{1}{\sigma_2(\vx^{\sss(n)})},
\end{equation}
for fixed $\lambda \in \bR$. Let $\cN(\bR_+)$ denote the space of counting measures on $\bR_+$ equipped with the vague topology.\nomenclature[cn]{$\cN(\bR_+)$}{Space of counting measures on $\bR_+$ equipped with the vague topology.}
Define $\mvUpsilon_n^{\sss(i)}:= (p_v/\sigma(\vp), v\in \cV^{\sss(i)})$ and view $(a_n^{\sss(i)}\sigma(\vp_n^{\sss(i)}), \mvUpsilon_n^{\sss(i)})$ as a random element of $\dS:= \bR_+\times \cN(\bR_+)$ (equipped with the product topology).
\nomenclature[ds]{$\dS$}{Space $\bR_+\times \cN(\bR_+)$ equipped with the product topology.}
 Finally, define
\begin{align*}
\cP_n:= \left(\left(a_n^{\sss(i)}\sigma(\vp_n^{\sss(i)}), \mvUpsilon_n^{\sss(i)}\right)\colon i\geq 1\right)
\end{align*}
viewed as an element of $\dS^{\infty}$, again equipped with the product topology induced by a single coordinate $\dS$. Now given an infinite vector $\vc\in l_0$ recall the process $\bar{V}_\lambda^{\vc}(\cdot)$ as in \eqref{eqn:vc-reflection}, the corresponding excursions $\sZ(\lambda)$ as in \eqref{eqn:excursions-def} and the corresponding excursion lengths in \eqref{eqn:vz-excursions-def}. Finally recall the definitions of $\barg^{\sss(i)}, \mvtheta^{\sss(i)}$ from \eqref{eqn:barg-mvtheta-def-comps}. Writing these out explicitly, define
\begin{equation}
\label{eqn:cp-infty}
	\cP_{\infty}:= ((\barg^{\sss(i)}, \mvtheta^{\sss(i)})\colon i\geq 1)= \left(Z_i(\lambda)\sqrt{\sum_{v\in \cZ_i(\lambda)} c_v^2}\;,\;  \Bigg(\frac{c_j}{\sqrt{\sum_{v\in \cZ_i(\lambda)} c_v^2}}: j\in \cZ_i(\lambda)\Bigg)\colon i\geq 1\right).
\end{equation}

\begin{prop}
	\label{prop:weights-good}
	The following hold under Assumption \ref{ass:mc-assumptions}:
	\begin{enumeratei}
		\item For every $i\geq 1$, $\sigma(\vp_n^{\sss(i)}) \probc 0$ as $n\to\infty$.
		\item $\cP_n\convd \cP_{\infty}$ on $\dS^{\infty}$ as $n\to\infty$. Further for every fixed $i\geq 1$, almost surely,
	\begin{equation}
	\label{eqn:cz-good-prop}
		  \sum_{v\in \cZ_i(\lambda)} c_v= \infty.
	\end{equation}
	\end{enumeratei}
	
\end{prop}

\noindent{\bf Proof of Theorem \ref{thm:mc-main-1}:} We prove the theorem assuming Proposition \ref{prop:weights-good}. By an application of Skorohod embedding we may assume that we are working on a probability space where the convergence in Proposition \ref{prop:weights-good} happens almost surely.
In particular, in this space, Assumption \ref{ass:vp-a} is satisfied almost surely for $\vp_n^{\sss(i)}$ for any fixed $i\geq 1$.  Now an application of Theorem \ref{thm:connected-convg} completes the proof.

\qed

\subsection{Verification of weight assumptions in maximal components}
Here we give the proof of Proposition \ref{prop:weights-good}. To ease notation, we will throughout assume $\lambda=0$. The general case follows in an identical fashion, but this assumption simplifies notation. We will write $V^{\vc}$ instead of $V^{\vc}_0$ for the process in \eqref{eqn:vc-process-def} with $\lambda =0$ and simply write $\cC_i$ for $\cC_i([\sigma_2(\vx^{\sss(n)})]^{-1})$.

 We start by describing an exploration scheme (developed in \cite{aldous-crit}) which simultaneously constructs the graph $\cG_n(\vx,t)$ and a ``breadth first'' walk.
 This was carefully analyzed in \cite{aldous-limic} to prove Theorem \ref{thm:aldous-limic}.

For every ordered pair $(u,v)$, let $\eta_{u,v}$ be an exponential random variable with rate $tx_v$ (independent across ordered pairs). Note that there is a simple relation between the connection probabilities of $\cG_n(\vx,t)$ given by \eqref{eqn:mc-connection} and the above random variables given by:
\begin{equation}
\label{eqn:mc-con-eta}
	q_{uv}:= \pr( \eta_{uv} < x_u).
\end{equation}
  At each stage $i\geq 1$, we have a collection of active vertices $\cA(i)$, a collection of explored vertices $\cO(i)$ and a collection of unexplored vertices $\cU(i)= [n]\setminus \cA(i)\cup \cO(i)$.

  Initialize with $\cO(1) = \emptyset$ and $\cA(1) = \set{v(1)}$, where the first vertex $v(1)$ is chosen by size-biased sampling, namely with probability proportional to vertex weights $\vx$. When possible we will suppress dependence on $n$ to ease notation. Now let $\cD(v(1)):=\set{v: \eta_{v(1),v}\leq x_{v(1)}}$ denote the collection of ``children'' of $v(1)$ and note that by \eqref{eqn:mc-con-eta} this generates the right connection probabilities in $\cG_n(\vx,t)$. Think of the associated $\eta_{v(1),v}$ values (for vertices connected to $v(1)$) as ``birth-times'' of these connections in the interval $[0,x_{v(1)}]$ and label the corresponding vertices as $v(2), v(3), \ldots v(|\cD(v(1))|+1)$. Update the process as $\cO(2):= \set{v(1)}$, $\cA(2):= \cD(v(1))$ and $\cU(2) = \cU(1)\setminus \cD(v(1))$.

   Associate with this construction a breadth-first walk as follows:
  \begin{equation}
  \label{eqn:breadth-walk}
  	Z_n(0):=0, \qquad Z_n(u):= -u + \sum_{v} x_v\ind\set{\eta_{v(1),v}\leq u}, \qquad 0\leq u\leq x_{v(1)}.
  \end{equation}
Recursively for $i\geq 2$ let $T_{i-1}:= \sum_{j=1}^{i-1} x_{v(j)}$. At this ``time'' we will explore the unexplored neighbors of $v(i)$. By this time, there are $|\cU(i)|:= i-1+|\cA(i)|$ vertices that have either been explored or are active.   Let $\cD(v(i)):= \set{v\in \cU(i): \eta_{v(i)v}\leq x_{v(i)}}$ and again label these as $v(i+|\cA(i)|), v(i+|\cA(i)|+1), \ldots v(i+|\cA(i)|+|\cD(v(i))|-1)$ in increasing order of their $\eta_{v(i)v}$ values. Update $\cO(i+1) = \cO(i)\cup\set{v(i)}$, $\cA(i+1) = \cA(i)\cup\cD(v(i))\setminus \set{v(i)}$ and $\cU(i+1) = \cU(i)\setminus \cD(v(i))$. Again update the walk as
\begin{equation}
\label{eqn:breadth-walk-update}
	Z(T_{i-1}+u) = Z(T_{i-1}) - u + \sum_{v\in \cD(v(i))} x_v \ind\set{\eta_{v(i),v}\leq u}, \qquad 0\leq u\leq x_{v(i)}.
\end{equation}
After finishing a component (which happens when $\cA(i) = \emptyset$ for some $i\geq 2$), choose the next vertex to explore in a size-biased manner from the unexplored set $\cU(i)$. If $\cU(i)=\emptyset$, then we have finished constructing the partition of the graph into the connected components.

Now note the following important properties of this exploration:
\begin{enumeratea}
	\item The ordering $({v(1)}, {v(2)}, \ldots, {v(n)})$ is a size-biased reordering of the vertex set $[n]$.
	\item If we start a new component at some stage $i$ with vertex $v(i)$, and finish exploring the component at stage $j\geq i$, then the walk satisfies
	\[Z(T_j) = Z(T_{i-1}) - x_{v(i)}, \qquad Z(u)\geq Z(T_j) \mbox{ on } T_{i-1} < u < T_j.\]
	Thus, the size of the component of $v(i)$, $\sum_{l=i}^{j} x_{v(l)}$ is essentially the length of the excursion of the walk beyond past minima.
\end{enumeratea}

  As a starting point in proving Theorem \ref{thm:aldous-limic}, Aldous and Limic \cite{aldous-limic} show the following result. Their result is more general (incorporating the presence of a ``Brownian component'') but we state their result as applied to our setting.
\begin{prop}[{\cite[Proposition 9]{aldous-limic}}]
	\label{prop:aldous-limic-walk}
	Consider the process $\set{\bar{Z}_n(s):s\geq 0}$ defined by setting $\bar{Z}_n(s) := Z(s)/\sigma_2$. Then under Assumption \ref{ass:mc-assumptions} $\bar{Z}_n\convd V^{\vc}$ as $n\to\infty$.
\end{prop}
Using this result, Aldous and Limic \cite{aldous-limic} show that the corresponding maximal excursions beyond past minima of $\bar{Z}_n$ also converge to the maximal excursions beyond past minima of $V^{\vc}_\lambda$, namely the excursion lengths of the reflected process $\bar{V}_\lambda^{\vc}$ (see \eqref{eqn:vc-reflection}) from zero. A consequence of the proof of Theorem \ref{thm:aldous-limic} in \cite{aldous-limic} using Proposition \ref{prop:aldous-limic-walk} is the following result:
\begin{lemma}\label{lem:tight-K-find}
	Fix $K$ and let $\cE_n(K)$ be the time required for the above construction to explore the maximal $K$ components $\set{\cC_i:1\leq i\leq K}$. Then $\set{\cE_n(K):K\geq 1}$ is tight.
\end{lemma}
In other words, for any fixed $K\geq 1$, the maximal length excursions of $\bar{V}^{\vc}$ are found in finite time. Thus, even though the total weight of vertices $\sigma_1\to\infty$, when exploring the graph in size-biased fashion, under Assumption \ref{ass:mc-assumptions} one needs only a finite amount of ``time'' to find the maximal components. Here time is measured in terms of the weight of vertices already explored.
Now define
\begin{equation}
\label{eqn:qn-qneps-def}
	S_{n,2}(u)= \sum_{i: T_{i}\leq u} \left(\frac{x_{v(i)}}{\sigma_2}\right)^2, \qquad R_n^{\eps}(u):= \sum_{i: T_i\leq u} \frac{x_{v(i)}^2}{\sigma_2^2}\ind\set{x_{v(i)}< \sigma_2\eps}.
\end{equation}
Thus, $S_{n, 2}(t)$ is the normalized sum of \emph{squares} of vertex weights of vertices explored by time $t$ and $R_n^{\eps}$ is the normalized sum of these squares where we only retain explored vertices with weight at most $\eps\sigma_2$.
Using the same set of exponential random variables $\set{\xi_j:j\geq 1}$  that arose in the definition of the process $V^{\vc}$ in \eqref{eqn:vc-process-def} define a new process
\begin{equation}
\label{eqn:vc-squares}
	S_{\sss\infty,2}(u):= \sum_{j=1}^{\infty} c_j^2 \ind\set{\xi_j\leq u}.
\end{equation}
The same proof techniques as in \cite{aldous-limic} now implies the following. Since the ideas basically follow from \cite{aldous-limic} we only sketch the proof.
\begin{lem}\label{lem:sum-squares}
Assumption \ref{ass:mc-assumptions} implies the joint convergence of the processes $(\bar{Z}_n(\cdot), S_{n,2}(\cdot))\convd(V^{\vc}(\cdot), S_{\infty,2}(\cdot))$ as $n\to\infty$.
\end{lem}
\noindent{\bf Proof:} Fix $K\geq 1$, and for each $i\geq 1$, let $\xi_i^{(n)}$ denote the time when vertex $i$ is added to the collection of active vertices.  Now consider the $K+1$ dimensional stochastic process
\begin{align*}
\vY_n^K(s):=\left(\bar{Z}_n(s),\frac{x_1}{\sigma_2} \ind\set{\xi_1^{(n)}\leq s}, \ldots, \frac{x_K}{\sigma_2} \ind\set{\xi_K^{(n)}\leq s}  \right), \qquad s\geq 0.
\end{align*}
Write
\[\vY_{\infty}^K(s):= (V^{\vc}(s), c_1\ind\set{\xi_1\leq s}, \ldots,c_K\ind\set{\xi_K\leq s} ).\]
In the proof of Proposition \ref{prop:aldous-limic-walk}, Aldous and Limic showed that $\vY_n^K\convd \vY_{\infty}^K$ for every fixed $K\geq 1$.
Thus to complete the proof, it is enough to show, for every fixed $A>0$ and $\eta >0$, $\limsup_{\eps\to 0}\limsup_{n\to\infty} \pr(R_n^{\eps}(A)> \eta)=0$.
Now as described on \cite[Page 17]{aldous-limic}, we can couple $(\xi_1^{(n)}, \xi_2^{(n)}, \ldots, \xi_n^{(n)})$ with a sequence of independent exponential random variables $(\tilde{\xi}_1^{(n)}, \tilde{\xi}_2^{(n)}, \ldots, \tilde{\xi}_n^{(n)})$ with $\tilde{\xi}_j^{(n)}$ having rate $x_j/\sigma_2$ such that $\tilde{\xi}_j^{(n)}\leq \xi_j^{(n)}$. Now write
\begin{equation}
\label{eqn:tildrn-def}
	\tildR_{n}^{\eps}(t):=\sum_{j: x_j<\eps\sigma_2} \frac{x_j^2}{\sigma_2^2} \ind\set{\tilde{\xi}_j^{(n)}\leq t}.
\end{equation}
Then it is enough to show
\begin{equation}
\label{eqn:rn-enough}
	\limsup_{\eps\to 0}\ \limsup_{n\to\infty}\ \E(\tildR_n^{\eps}(A)) = 0,
\end{equation}
which is trivial since
\begin{equation}
\label{eqn:tildrn-bound}
	\E\left(\tildR_n^{\eps}(A)\right)\leq A\sum_{j: x_j\leq \eps \sigma_2} \left(\frac{x_j}{\sigma_2}\right)^3 \to A\sum_{j: c_j< \eps} c_j^3.
\end{equation}
We have used both \eqref{eqn:sigma2-sigma3} and \eqref{eqn:xj-to-cj} in the last convergence assertion. Thus, first letting $n\to\infty$ and then $\eps\to 0$ completes the proof.\qed

\medskip


We can now complete the proof of Proposition \ref{prop:weights-good}.  First, note that to prove \eqref{eqn:cz-good-prop}, it is enough to show that for any two rationals $r<s$, $\sum_{j} c_j\ind\set{r\leq \xi_j \leq s} = \infty$ almost surely where $\xi_j$ are the associated exponential rate $c_j$ random variables. This, however, is trivially true as $\sum_j c_j^2 =\infty$.

To prove the other assertions, define, for $i\geq 1$, the point processes $\Xi_n^{\sss(i)}:=\set{x_u/\sigma_2: u\in \cC_i}$, namely the rescaled vertex weights in the $i$th maximal component. Analogously define $\Xi_{\infty}^{\sss(i)} = \set{c_v: v\in \cZ_i}$, namely the collection of jumps in the $i$th largest excursion of $\bar{V}^{\vc}$. Let
\begin{equation}
\label{eqn:sn-def}
s_n^{\sss(i)} = \sum_{v\in \cC_i} \frac{x_v^2}{\sigma_2^2},\quad\text{ and }\ s_{\infty}^{\sss(i)}:= \sum_{v\in \cZ_i} c_v^2,	
\end{equation}
  for the normalized sum of squares of vertex weights in a component. 
Define
\[\tilde{\cP}_n:= \left(\left(\mass(\cC_i), s_n^{\sss(i)}, \Xi_n^{\sss(i)}\right), i\geq 1\right), \qquad \tilde{\cP}_{\infty}:= \left(\left(Z_i,s_{\infty}^{\sss(i)}, \Xi_{\infty}^{\sss(i)}\right), i\geq 1\right).\]
We will view these as random elements of $\tilde{\dS}^{\infty}$ where $\tilde{\dS}:= \bR^2\times \cN(\bR)$.
Lemma \ref{lem:tight-K-find} and Lemma \ref{lem:sum-squares} now imply the following:

\begin{lem}\label{lem:pp-length}
	As $n\to\infty$, $\tilde{\cP}_n \convd \tilde{\cP}_{\infty}$ on $\tilde{\dS}^{\infty}$.
\end{lem}
Expressing the functionals that arise in Proposition \ref{prop:weights-good} in terms of vertex weights in maximal components completes the proof. Indeed,
\[\sigma(\vp_n^{\sss(i)})= \frac{\sqrt{\sum_{v\in\cC_i} x_v^2}}{\sum_{v\in\cC_i} x_v}=\frac{\sigma_2\sqrt{s_n^{\sss(i)}}}{\mass(\cC_i)} \to 0,\]
as $n\to\infty$.
 %
The proof of $\cP_n\convd\cP_{\infty}$ is similar. \qed


\subsection{Gromov-weak convergence in Theorem \ref{thm:rank-one}.}\label{sec:rank-one-gromov-weak}
That convergence in \eqref{eqn:conergence-rank-one} holds with respect to Gromov-weak topology is an easy consequence of Theorem \ref{thm:mc-main-1}. Indeed, setting
\begin{align}\label{eqn:x-i-from-w-i}
x_i= n^{-\frac{\tau-2}{\tau-1}}w_i\ \text{ and }\
t_n=\frac{1}{\ell_n}\left(1+\lambda n^{-(\tau-3)/(\tau-1)}\right)n^{\frac{2(\tau-2)}{\tau-1}},
\end{align}
we can write $\NRnwl$ as the model $\cG(\mvx, t_n)$ where $\mvx=\mvx^{\sss(n)}:=(x_i: i\in[n])$. A direct computation will show that $\mvx^{\sss(n)}$ satisfies Assumption \ref{ass:mc-assumptions} with the entrance boundary $\vc^{\nr}$ defined in \eqref{eqn:def-c-nr}. Note also that
\begin{align*}
t_n-\frac{1}{\sigma_2(\mvx^{\sss(n)})}
&=\frac{n^{\frac{2(\tau-2)}{\tau-1}}}{\sum_{i\in[n]}w_i^2}\left(\frac{1}{\ell_n}\sum_{i\in[n]}w_i^2-1\right)
+\frac{\lambda}{\ell_n/n}.
\end{align*}
Under the assumptions of Theorem \ref{thm:rank-one}, $\ell_n/n\to\E W$ and $\sum_{i}w_i^2/n\to\E W^2=\E W$. Further, by \cite[Lemma 2.2]{SBVHVJL12},
\begin{align*}
\frac{1}{\ell_n}\sum_{i\in[n]}w_i^2=1+\zeta n^{-(\tau-3)/(\tau-1)}+o(n^{-(\tau-3)/(\tau-1)}),
\end{align*}
where $\zeta$ is as defined in \eqref{eqn:t-nr}. Combining these observations, we see that
\begin{align*}
t_n-(\sigma_2(\mvx^{\sss(n)}))^{-1}\to t^{\nr}_{\lambda}\ \text{ as }n\to\infty,
\end{align*}
where $t^{\nr}_{\lambda}$ is as in \eqref{eqn:t-nr}. Since $n^{(\tau-3)/(\tau-1)}\sigma_2(\mvx^{\sss(n)})\to\E W$, we conclude that $\vM_{\infty}^{\nr}(\lambda)$ defined in \eqref{eqn:def-M-nr} is the Gromov-weak limit of $n^{-(\tau-3)/(\tau-1)}\vM_n^{\nr}(\lambda)$, where $\vM_n^{\nr}(\lambda)$ is as in \eqref{eqn:def-M-n}.

\begin{rem}
Theorem \ref{thm:mc-main-1} is stated for a fixed $\lambda\in\bR$, but in the argument just given, we have to work with a sequence, namely $t_n-(\sigma_2(\mvx^{\sss(n)}))^{-1}$ converging to $t^{\nr}_{\lambda}$. This, however, does not make any difference. Indeed, the proof of \cite[Proposition 9]{aldous-limic} can be imitated to prove the same result in the setup where we have a sequence converging to $t$ instead of a fixed $t$, and no new idea is involved here. (In \cite[Lemma 27]{aldous-limic}, Aldous and Limic prove a similar result for the multiplicative coalescent. They do not, however, explicitly state the convergence of the associated process under the same assumption.)
\end{rem}

\section{Proofs: Convergence in Gromov-Hausdorff-Prokhorov topology}
\label{sec:proof-convg-dghp}

In this section, we improve Gromov-weak convergence in Theorem \ref{thm:rank-one} to Gromov-Hausdorff-Prokhorov convergence. To do so, we will rely on \cite[Theorem 6.1]{AthLorWin14} which gives a criterion for convergence in Gromov-Hausdorff-weak topology. We do not give the definition of Gromov-Hausdorff-weak topology and instead refer the reader to \cite[Definition 5.1]{AthLorWin14}. Convergence in Gromov-Hausdorff-weak topology implies convergence in Gromov-Hausdorff-Prokhorov topology when we are working with metric measure spaces having full support (i.e., the support of the measure is the entire metric space). This is true in our situation. Indeed, it is a trivial fact that $\cC_i(\lambda)$ has full support. Further, the mass measure on an inhomogeneous continuum random tree has full support which implies that the same is true for $M_i^{\nr}(\lambda)$.

Applying \cite[Theorem 6.1]{AthLorWin14} to our situation, we see that it is enough to prove the following lemma:

\begin{lemma}[Global lower mass-bound]
\label{lem:ghp-from-gw}
Let $\cC_i(\lambda)$ be the $i$th largest component of $\NRnwl$. Then the following assertion is true:
For each $i\geq 1$, $v\in[n]$ and $\delta>0$, let $B(v, \delta)$ denote the intrinsic ball (in $\NRnwl$) of radius $\delta n^{(\tau-3)/(\tau-1)}$ around $v$ and set
    \[\fm_i^{(n)}(\delta)=\inf\set{n^{-\frac{\tau-2}{\tau-1}}\sum_{j\in B(v, \delta)}w_j\ \bigg|\ v\in\cC_i(\lambda)}.\]
Then the sequence $\set{\left(\fm_i^{(n)}(\delta)\right)^{-1}}_{n\geq 1}$ is tight.
\end{lemma}

Lemma \ref{lem:ghp-from-gw} ensures compactness of the spaces $M_i^{\nr}(\lambda)$ which, in turn, implies compactness of the spaces $M_i^{\vc}(\lambda)$ when $\vc=(c_1, c_2,\hdots)$ is of the form \eqref{eqn:c-tau-lamb-def}, thus proving the first assertion in Theorem \ref{thm:mc-main-2-tau}.

Before moving on to the proof of Lemma \ref{lem:ghp-from-gw}, we state a result that essentially says that instead of looking at the largest components, we can work with the components of high-weight vertices. This observation will be used to prove the global lower-mass bound:

\begin{prop}\label{prop:enough-to-work-with-high-weight-vertices}
For every $\eps>0$ and $k\geq 1$, there exists $K=K(\eps, k, \lambda)>0$ such that
\[\pr\left([K]\cap\cC_i(\lambda)=\emptyset\text{ for some }1\leq i\leq k\right)\leq\eps.\]
\end{prop}

Proposition \ref{prop:enough-to-work-with-high-weight-vertices} follows trivially from \cite[Theorem 1.6 (a)]{SBVHVJL12} and \cite[Theorem 1.1]{SBVHVJL12}.
\subsection{Bound on size of $\eps n^{(\tau-3)/(\tau-1)}$-nets for the largest components}
For convenience, we set
	\begin{align}\label{eqn:def-eta-rho}
	\eta=(\tau-3)/(\tau-1)\qquad \text{ and }\qquad \rho=(\tau-2)/(\tau-1).
	\end{align}
\nomenclature{$\eta, \rho$}{Critical exponents.}
The purpose of this section is to prove a strong result (Proposition \ref{prop-diam-non-hubs} stated below) that gives control over the number of intrinsic balls of radius $\eps n^{\eta}$ needed to cover the largest components. This acts as a crucial ingredient in the proof of Lemma \ref{lem:ghp-from-gw} as well as the proof of the bound on the upper box-counting dimension.


\begin{prop}[Small diameter after removing high-weight vertices]
\label{prop-diam-non-hubs}
For every $\vep, \delta>0$,
and $N=N(\eps):=\eps^{-\delta-1/\eta}$,
	\begin{align}\label{eqn:P-E-exponential-bound}
	\pr\left(\diam\left(\NRnwl\setminus [N]\right)>\eps n^{\eta}\right)
	\leq c_{\delta}\exp\left(-C/\eps^{\delta\eta}\right),
	\end{align}
for all $n$ sufficiently large, a positive constant $c_{\delta}$ depending on $\delta$ and a universal constant $C>0$. Here $\NRnwl\setminus [N]$ denotes the graph obtained by removing all vertices with labels in $[N]$ and the edges incident to them from the graph $\NRnwl$.
\end{prop}


We continue to prove Proposition \ref{prop-diam-non-hubs}. Write
	\begin{align}\label{eqn:def-E-n}
	E_n=\{\diam(\NRnwl\setminus [N]) \leq \vep n^{\eta}\}.
	\end{align}
The proof consists of four steps. In the first step, we reduce the proof to the study of the height of mixed-Poisson branching processes. In the second step, we ensure that we can take $\lambda=0$, while in the third step, we study the survival probability of such critical infinite-variance branching processes. In the fourth and final step, we prove the claim.

\medskip

\noindent{\bf Comparison to mixed-Poisson branching processes.}
Let $\cC_{\res}(i)$ be the cluster of $i$ in the (restricted) random graph on the vertex set $[n]\setminus [i-1]$ with edge probabilities $q_{k\ell}(\mvw(\lambda))$ for $k,\ell\in [n]\setminus [i-1]$, where $q_{k\ell}(\mvw(\lambda))$ is as in \eqref{eqn:nr-connection}.

Note that the event $E_n^c$ implies the existence of $i>N$ such that the following happens:
\begin{inparaenuma}
\item The diameter of the component of $i$ in $\NRnwl\setminus [N]$ is bigger than $\vep n^{\eta}$.
\item No $j\in\set{N+1,\hdots, i-1}$ belongs to the component of $i$ in $\NRnwl\setminus [N]$.
\end{inparaenuma}
In particular, $\diam(\cC_{\res}(i))\geq \vep n^{\eta}$ for this $i$.
Thus,
	\begin{align}\label{eqn:E-union-bound}
	\pr(E_n^c)\leq \sum_{i>N} \pr\left(\diam(\cC_{\res}(i))>\vep n^{\eta}\right).
	\end{align}

Now the random graph $\NRnwl$ restricted to $[n]\setminus [i-1]$ is the Norros-Reittu random graph $\NR_n(\mvw^{\sss (i)}(\lambda))$, where $\mvw^{\sss (i)}(\lambda)=(w_j^{\sss(i)}(\lambda)\colon j\in[n])$, $w_j^{\sss(i)}(\lambda)=w_j(\lambda)\ell_n^{\sss(i)}/\ell_n$ for $j\in[n]\setminus [i-1]$ and $w_j^{\sss(i)}(\lambda)=0$ for $j\in[i-1]$, and $\ell_n^{\sss (i)}=\sum_{k=i}^n w_k$. Indeed, this follows from the simple observation
	\begin{align*}
	\frac{w_k^{\sss(i)}(\lambda)w_{\ell}^{\sss(i)}(\lambda)}{\sum_{r=i}^n w_{r}^{\sss(i)}(\lambda)}
	=\left(1+\frac{\lambda}{n^{\eta}}\right)\frac{w_k w_{\ell}}{\ell_n}.
	\end{align*}
Write $W_n^{\sss(i)}(\lambda)$ for a random variable whose distribution is given by $(n-i+1)^{-1}\sum_{j=i}^{n}\delta_{w_j^{\sss(i)}(\lambda)}$, and for any non-negative random variable $X$ with $\E X>0$, let $X^\circ$ be the random variable having the size-biased distribution given by
	\[
	\pr(X^\circ\leq x)=\E(X\ind_{\set{X\leq x}})/\E X.
	\]
\nomenclature{$X^\circ$}{Random variable having the size-biased distribution.}
We will use the following comparison to a mixed-Poisson branching process:
\begin{lemma}[Domination by a mixed-Poisson branching process]
\label{lem-dim-MPBP}
Fix $i\in [n]$ and consider $\NR_n(\mvw^{\sss (i)}(\lambda))$. Then, there exists a coupling of $\cC_{\res}(i)$ and a branching process where the root has a $\Poi(w_i^{\sss(i)}(\lambda))$ offspring distribution while every other vertex has a $\Poi((W_n^{\sss(i)}(\lambda))^\circ)$ offspring distribution such that in the breadth-first exploration of $\cC_{\res}(i)$ starting from $i$, each vertex $v\in \cC_{\res}(i)$ has at most the number of children as in the branching process.
\end{lemma}
\noindent{\bf Proof:} See \cite[Proposition 3.1]{NorRei06}.
\qed
\medskip

It immediately follows from Lemma \ref{lem-dim-MPBP} that
	\begin{align}\label{eqn:diam-c-geq-i-bound}
	\pr\left(\diam(\cC_{\res}(i))>\vep n^{\eta}\right)\leq \pr\left(\he({\overline{T}_n}^{\sss (i)}(\lambda))>\vep n^{\eta}/2\right),
	\end{align}
where $\overline{T}_n^{\sss (i)}(\lambda)$ is a mixed-Poisson branching process tree whose root has a $\Poi(w_i^{\sss(i)}(\lambda))$ offspring distribution and every other vertex has a $\Poi((W_n^{\sss (i)}(\lambda))^\circ)$ offspring distribution. As before, $\he(\vt)$ denotes the height of the tree $\vt$.

When $\he(\overline{T}_n^{\sss (i)}(\lambda))>\vep n^{\eta}/2$, at least one of the subtrees of the root needs to have height at least $\vep n^{\eta}/2$. Combining this observation with \eqref{eqn:E-union-bound} and \eqref{eqn:diam-c-geq-i-bound}, we get
	\begin{align}\label{eqn:E-c-upper-bound}
	\pr(E_n^c)\leq \sum_{i>N} \E\left[\Poi(w_i^{\sss(i)}(\lambda))\right]\pr\left(\he\left(T_n^{\sss (i)}(\lambda)\right)\geq \eps n^{\eta}/2\right)
	\leq \sum_{i>N} w_i^{\sss(i)}(\lambda) \pr\left(\he\left(T_n^{\sss (i)}(\lambda)\right)\geq \vep n^{\eta}/2\right),
	\end{align}
where $T_n^{\sss (i)}(\lambda)$ is a branching process tree where every vertex has a $\Poi((W_n^{\sss (i)}(\lambda))^\circ)$ offspring distribution.

We make the convention of writing $T_n^{\sss (i)}$, $W_n^{\sss (i)}$ etc.\ instead of $T_n^{\sss (i)}(0)$, $W_n^{\sss (i)}(0)$ etc. With this notation, it is easy to see that $W_n^{\sss (i)}(\lambda)\stackrel{d}{=}(1+\lambda n^{-\eta}) W_n^{\sss (i)}$ and hence $(W_n^{\sss (i)}(\lambda))^\circ\stackrel{d}{=}(1+\lambda n^{-\eta}) (W_n^{\sss (i)})^\circ$.

\medskip

\noindent{\bf The survival probability of mixed-Poisson branching processes.} We would like to compare our mixed-Poisson branching process with an offspring distribution that is independent of $n$.
For this, we rely on the following two lemmas:
\begin{lemma}[Mixed-Poisson branching processes of different parameters]
\label{lem-MPBP-par}
Let $T_n^{\sss (i)}$ and $T_n^{\sss (i)}(\lambda)$ be as above. Assume further that $\lambda\geq 0$. Then, for each $k\geq 1$,
	\begin{align*}
	\pr\left(\he(T_n^{\sss (i)}(\lambda))\geq k\right)\leq (1+\lambda n^{-\eta})^k \cdot\pr\left(\he(T_n^{\sss (i)})\geq k\right).
	\end{align*}
\end{lemma}

\noindent{\bf Proof:} We follow \cite[Proof of Lemma 3.4(1)]{HofNac12}.
Writing $\delta=1+\lambda n^{-\eta}$, we note that we can obtain $T_n^{\sss (i)}$ as a subtree of $T_n^{\sss (i)}(\lambda)$ by killing every child independently with probability $1-\delta^{-1}$. Write $\mathcal{A}$ for the event in which $\he(T_n^{\sss (i)}(\lambda))\geq k$ and no vertex in the leftmost path of length $k$ starting from the root in $T_n^{\sss (i)}(\lambda)$ is killed. Then
	\begin{align}\label{PA-prob}
    	\pr(\mathcal{A}) = \delta^{-k} \pr\left(\he(T_n^{\sss (i)}(\lambda))\geq k\right).\notag
	\end{align}
Indeed, the probability of the leftmost path surviving is precisely $1/\delta^k$. To finish the proof, note that $\mathcal{A}$ implies $\he(T_n^{\sss (i)})\geq k$, so that
	\begin{align*}
	\pr\left(\he(T_n^{\sss (i)})\geq k\right)\geq \pr(\mathcal{A})=\delta^{-k} \pr\left(\he(T_n^{\sss (i)}(\lambda))\geq k\right),
	\end{align*}
which is the desired inequality.
\qed


\begin{lemma}[Stochastic bound by $n$-independent variable]
\label{lem-MPBP-ub}
Under Assumption \ref{ass:density}, the random variable $(W_n^{\sss (i)})^\circ$ is stochastically upper bounded by $W^{\circ}$ where $W\sim F$, i.e.,  $(W_n^{\sss (i)})^\circ\stackrel{\mathrm{st}}{\leq} W^{\circ}$.
\end{lemma}


\noindent{\bf Proof:} First we make the following elementary observation: if $a_1, a_2, b_1, b_2$ are positive numbers such that
\begin{equation*}
\frac{a_1}{b_1}\leq\frac{a_2}{b_2},\ \text{ then }\ \frac{a_1}{b_1}\leq\frac{a_1+a_2}{b_1+b_2}\leq\frac{a_2}{b_2}.
\end{equation*}
Repeated application of the above will yield the following simple inequality: if $\set{a_n}_{n\geq 1}$ and $\set{b_n}_{n\geq 1}$ are sequences of positive numbers satisfying
\begin{align}\label{eqn:122}
\frac{a_1}{b_1}\leq\frac{a_2}{b_2}\leq\frac{a_3}{b_3}\leq\hdots,\ \text{ then }\ \frac{a_1}{b_1}\leq\frac{a_1+a_2}{b_1+b_2}\leq\frac{a_1+a_2+a_3}{b_1+b_2+b_3}\leq\hdots.
\end{align}

Recall that $\iota$ denotes the leftmost point of the support of $F$, and note that from \eqref{eqn:wi-construc} it follows that $\int_{w_j}^{\infty}f=j/n$, $j=1, 2, \hdots, n$ (note also that $w_n=\iota$). Define the function $h_n: [\iota,w_1)\to (\iota,\infty)$ by
$
\int_y^{h_n(y)}f=1/n.
$
This immediately implies
\begin{align}\label{eqn:123}
f(h_n(y))h_n'(y)=f(y).
\end{align}
Let $g_n: [\iota,w_1)\to (0,\infty)$ be given by
\begin{align*}
g_n(y)=\frac{y}{\int_y^{h_n(y)}uf(u)\ du}.
\end{align*}
A direct computation and an application of \eqref{eqn:123} yields
\begin{align*}
\bigg(\int_y^{h_n(y)}uf(u)\ du\bigg)^2 g_n'(y)=\int_y^{h_n(y)}uf(u)\ du-yf(y)\left(h_n(y)-y\right).
\end{align*}
Since $uf(u)$ is non-increasing on $[\iota, \infty)$ under Assumption \ref{ass:density}, we conclude that $g_n'(y)\leq 0$ on $(\iota, w_1)$. Thus, $g_n(\cdot)$ is non-increasing on $[\iota, w_1)$. By right continuity, we can define $g_n(w_1)=w_1/(\int_{w_1}^{\infty}uf(u)\ du)$. Since $w_n\leq w_{n-1}\leq\hdots\leq w_1$, we conclude that $g_n(w_1)\leq g_n(w_2)\leq\hdots\leq g_n(w_n)$. Clearly $h_n(w_j)=w_{j-1}$ for $j=2,\hdots,n$. Thus
\begin{align*}
\frac{w_1}{\int_{w_1}^{\infty}uf(u)\ du}\leq
\frac{w_2}{\int_{w_2}^{w_1}uf(u)\ du}\leq
\hdots\leq
\frac{w_n}{\int_{\iota}^{w_{n-1}}uf(u)\ du}.
\end{align*}
Now an application of \eqref{eqn:122} gives
\begin{align*}
\frac{w_1+w_2+\hdots+w_k}{\int_{w_k}^{\infty}uf(u)\ du}\leq\frac{w_1+w_2+\hdots+w_n}{\int_{\iota}^{\infty}uf(u)\ du},\ \ k=1,2,\hdots,n,
\end{align*}
which is equivalent to
\begin{align*}
\pr\big((W_n^{\sss (i)})^\circ\geq w_k\big)
=\frac{w_1+w_2+\hdots+w_k}{w_1+w_2+\hdots+w_n}\leq\frac{\int_{w_k}^{\infty}uf(u)\ du}{\int_{\iota}^{\infty}uf(u)\ du}
=\pr\big(W^\circ\geq w_k\big),\  \ k=1,2,\hdots,n.
\end{align*}
This concludes the proof. \qed	
\medskip

We continue to study the survival probability of mixed-Poisson branching processes with infinite variance offspring distribution:
\begin{lemma}[Survival probability of infinite-variance MPBP]
\label{lem-surv-prob-MPBP}
Let $T$ denote a mixed-Poisson branching process tree with offspring distribution $\Poi(W^{\circ})$. Then, there exists a constant $c_{\sss \ref{lem-surv-prob-MPBP}}$ such that for all $m\geq 1$,
	\begin{align*}
	\pr(\he(T)\geq m)\leq c_{\sss \ref{lem-surv-prob-MPBP}} m^{-1/(\tau-3)}.
	\end{align*}
\end{lemma}

\noindent{\bf Proof:} This is a well-known result. We sketch the proof briefly for completeness. Recall the following facts about $W^\circ$:
\begin{inparaenuma}
\item$\expec[W^\circ]=\nu=1$ and
\item for $x\rightarrow \infty$, $\pr(W^\circ>x)=c x^{-(\tau-2)}(1+o(1)).$
\end{inparaenuma}
By the {\em Otter-Dwass formula}, which describes the distribution of the total progeny of a branching process (see \cite{Dwas69} for the special case when the branching process starts with a single individual, \cite{Otte49} for the more general case, and \cite{HofKea07} for a simple proof based on induction), we have
	\begin{equation}\label{eqn:OD-form}
	\pr(|T|=k)=\frac{1}{k} \pr\left(\sum_{i=1}^k X_i = k-1\right),\notag
	\end{equation}
where $X_i$ are i.i.d. random variables distributed as $W^\circ$.
By \cite[Proposition 2.7]{Hofs09a}, in our situation, $\prob (\sum_{i=1}^k X_i = k-1)\leq ck^{-1/(\tau-2)}$, so that
	\begin{align}\label{eqn:OD-form-simp}
	\pr(|T|= k)\leq ck^{-(\tau-1)/(\tau-2)}\qquad \text{ and }\qquad \pr(|T|\geq k)\leq ck^{-1/(\tau-2)}.
	\end{align}
Take $k=m^{(\tau-2)/(\tau-3)}$ in the second inequality in \eqref{eqn:OD-form-simp} to get
	\begin{align*}
	\pr(\he(T)\geq m)\leq c m^{-1/(\tau-3)}+\prob\left(\he(T)\geq m, |T|\leq m^{(\tau-2)/(\tau-3)}\right),
	\end{align*}
where $|T|$ denotes the total number of vertices in $T$. We condition on the size $|T|$ and write
	\begin{align}\label{eqn:condition-on-|T|}
	\prob\big(\he(T)\geq m, |T|\leq m^{(\tau-2)/(\tau-3)}\big)
	&=\sum_{k=1}^{m^{(\tau-2)/(\tau-3)}} \pr\left(\he(T)\geq m\ \big|\ |T|=k\right)
	\pr(|T|=k)\notag\\
	&\leq c \sum_{k=1}^{m^{(\tau-2)/(\tau-3)}}
	\pr\left(\he(T)\geq m\ \big|\ |T|=k\right)k^{-\frac{\tau-1}{\tau-2}}.
	\end{align}
By \cite[Theorem 4]{Kort15}, there exists a $\kappa>1$ such that, uniformly for $u\geq 1$,
	\begin{align*}
	\pr\left(\he(T)\geq u k^{(\tau-3)/(\tau-2)}\ \big|\ |T|=k\right)\leq \e^{-a u^{\kappa}}.
	\end{align*}
Combining this with \eqref{eqn:condition-on-|T|}, we get
	\begin{align*}
	\pr\left(\he(T)\geq m, |T|\leq m^{(\tau-2)/(\tau-3)}\right)
	\leq \sum_{k=1}^{m^{\frac{\tau-2}{\tau-3}}}	\exp\left(-a \left(mk^{-\frac{\tau-3}{\tau-2}}\right)^{\kappa}\right)k^{-\frac{\tau-1}{\tau-2}}=\Theta\left(m^{-1/(\tau-3)}\right),
	\end{align*}
as required.
\qed
\medskip

\noindent
{\bf Proof of Proposition \ref{prop-diam-non-hubs}:} Clearly
\begin{align*}
	\pr\left(\he\left(T_n^{\sss (i)}\right)\geq m\right)
	\leq \E\left[\left(W_n^{\sss (i)}\right)^\circ\right]\pr\left(\he\left(T_n^{\sss (i)}\right)\geq m-1\right)
	=:\nu_n^{\sss (i)}\pr\left(\he\left(T_n^{\sss (i)}\right)\geq m-1\right),
\end{align*}
where
\begin{align}\label{eqn:nu-n-i}
	\nu_n^{\sss (i)}=\E\left[\left(W_n^{\sss (i)}\right)^\circ\right]
	=\frac{\sum_{j\geq i} (w_j^{\sss(i)})^2}{\sum_{j\geq i} w_j^{\sss(i)}}
	=\left(\frac{\ell_n^{\sss(i)}}{\ell_n}\right)
	\frac{\sum_{j\geq i} w_j^2}{\sum_{j\geq i} w_j}=\frac{\sum_{j\geq i} w_j^2}{\ell_n}.
\end{align}
Iterating this $\vep n^{\eta}/4$ times, we get 
\begin{align}\label{eqn:iterate}
	\pr\left(\he\left(T_n^{\sss (i)}\right)\geq\vep n^{\eta}/2\right)
	&\leq (\nu_n^{\sss (i)})^{\vep n^{\eta}/4} \pr\left(\he\left(T_n^{\sss (i)}\right)\geq\vep n^{\eta}/4\right)\\
	&\leq (\nu_n^{\sss (i)})^{\vep n^{\eta}/4}\pr\left(\he\left(T\right)\geq\vep n^{\eta}/4\right)
	 \leq  (\nu_n^{\sss (i)})^{\vep n^{\eta}/4}\times c_{\sss \ref{lem-surv-prob-MPBP}} \left(\frac{4}{\eps}\right)^{1/(\tau-3)}\frac{1}{n^{1/(\tau-1)}},\notag
\end{align}
where the second inequality is a consequence of Lemma \ref{lem-MPBP-ub} and the last step follows from Lemma \ref{lem-surv-prob-MPBP}.

Substituting the estimate \eqref{eqn:iterate} into \eqref{eqn:E-c-upper-bound} leads to
	\begin{align}\label{eqn:97}
	\pr(E_n^c)\leq c
	\eps^{-1/(\tau-3)} n^{-1/(\tau-1)}\left(1+\frac{\max\set{\lambda, 0}}{n^{\eta}}\right)^{1+\eps n^{\eta}/2} \sum_{i>N} w_i (\nu_n^{\sss (i)})^{\vep n^{\eta}/4}
	\end{align}
for some constant $c$. Here we have used Lemma \ref{lem-MPBP-par} and the simple fact that $w_i^{\sss (i)}\leq w_i$.

Next, note that it is an easy consequence of \eqref{eqn:tau-c-def} that there exist constants $c', c''>0$ such that for all $i\in[n]$,
	\begin{align}\label{eqn:98}
	w_i \leq c'\left(\frac{n}{i}\right)^{1/(\tau-1)}\qquad \text{ and }\qquad \sum_{j=1}^{i}w_j^2\geq c''\sum_{j=1}^{i}\left(\frac{n}{i}\right)^{2/(\tau-1)}.
	\end{align}
Further, \cite[Lemma 2.2]{SBVHVJL12} implies that $\nu_n^{\sss(1)}<1$ for large $n$. Hence, for every $i\geq 2$,
	\begin{align*}
 	\nu_n^{\sss (i)}=\nu_n^{\sss(1)}-\frac{1}{\ell_n}\sum_{j=1}^{i-1}w_j^2
	\leq 1-C n^{-\eta} i^{\eta}\leq \exp\left(-C n^{-\eta} i^{\eta}\right)
	\end{align*}
for some $C>0$. Here, we have used the second inequality in \eqref{eqn:98}. Combining this estimate with \eqref{eqn:97} and the first inequality in \eqref{eqn:98}, we end up with
	\begin{align*}
	\pr(E_n^c)\leq C'\eps^{-1/(\tau-3)} \sum_{i>N} i^{-1/(\tau-1)}\exp\left(- C\eps i^{\eta}/4\right)
	\end{align*}
for some $C'>0$. Taking $N=\eps^{-\delta-1/\eta}$, we arrive at
	\begin{align}\label{eqn:99}
	\pr(E_n^c)&\leq C'\eps^{-1/(\tau-3)}N^{-1/(\tau-1)} \sum_{i>N} \exp\left(- C\eps i^{\eta}/4\right)\notag\\
	&\leq C'\eps^{\delta/(\tau-1)}\sum_{k=0}^{\infty}\sum_{i=N2^k}^{N 2^{k+1}-1}\exp\left(- C\eps i^{\eta}/4\right)\notag\\
	&\leq C'\eps^{\delta/(\tau-1)}N\sum_{k=0}^{\infty}2^k\exp\left(- C\eps (N2^k)^{\eta}/4\right).
	\end{align}
Note that $\eps N^{\eta}=\eps^{-\delta\eta}$. A little more work after plugging this into \eqref{eqn:99} will lead to \eqref{eqn:P-E-exponential-bound}.
\qed


\subsection{Proof of global lower-mass bound}
In this section, we complete the proof of Lemma \ref{lem:ghp-from-gw}.
We start with some preliminaries:
\begin{lemma}[Weight of size-biased reordering]
\label{lem-sb-reor}
Let $\pi_v(1)=v$ and $(\pi_v(i): i\in[n]\setminus\{1\})$ be a size-biased reordering on $[n]\setminus \{v\}$ where the size of vertex $v'$ is proportional to $w_{v'}$ for $v'\in[n]\setminus\set{v}$. Then, for every \ch{$k=o(n)$}, there exists a $J>0$ such that
	\begin{align*}
	\pr\bigg(\exists v\colon \sum_{\ch{i=1}}^k w_{\pi_v(i)}\leq k/2\bigg)\leq n\e^{-Jk}.
	\end{align*}
\end{lemma}

\noindent{\bf Proof:} See \cite[Proof of Lemma 5.1]{SBVHVJL12}. \qed
\medskip

Recall the definitions of $\eta$ and $\rho$ from \eqref{eqn:def-eta-rho}. Recall that for $v\in [n]$, $B(v,\delta)$ denotes the intrinsic ball (in $\NRnwl$) around $v$ or radius $\delta n^{\eta}$. We will use the following bound on the weight of balls:

\begin{lemma}[Weights of balls around high-weight vertices cannot be too small]
\label{lem-weights}
For every $\eps>0$ and $i\geq 1$, there exist $n_{i, \eps}$ large and $\delta_{i,\eps}>0$ such that for all $n\geq n_{i,\eps}$ and $\delta\in(0, \delta_{i,\eps}]$,
	\begin{align}\label{weight-balls-lb}
	\pr\left(\sum_{j\in B(i,\delta)}w_j \leq \left(\frac{c_{\sss F}}{2i}\right)^{1/(\tau-1)} \frac{\delta n^{\rho}}{2}\right)
	\leq n\exp\left(-\frac{J \delta n^{\rho}}{i^{1/(\tau-1)}}\right) +\frac{\eps}{2^{i}}.
	\end{align}
\end{lemma}

\noindent{\bf Proof:}
We rely on a cluster exploration used in \cite{SBVHVJL12} which we describe next.
We denote by $(Z_l(i))_{l\geq 0}$ the exploration process of $\cluster(i)$, the cluster containing $i$, starting from $i$, in the breadth-first search, where $Z_0(i)=1$ and where $Z_1(i)$ denotes the number of potential neighbors of the initial vertex $i$.
\nomenclature{$\cluster(i)$}{Component containing node $i$.}
The variable $Z_l(i)$ has the interpretation of the number of potential neighbors of the first $l$ explored potential vertices in the cluster whose neighbors have not yet been explored. As a result, we explore by taking one vertex of the `stack' of size $Z_l(i)$, drawing its mark and checking whether it is a real vertex, followed by drawing its number of potential neighbors. Thus, we set $Z_0(i)=1, Z_1(i)=\Poi(w_i)$, and note that, for $l\geq 2$, $Z_l(i)$ satisfies the recursion relation
    \eqn{
    \label{Zl-recur}
    Z_l(i)=Z_{l-1}(i)+X_l-1,\nn
    }
where $X_l$ denotes the number of potential neighbors of the $l$th potential vertex that is explored, where $X_1=X_1(i)=\Poi(w_i)$. More precisely, when we explore the $l$th potential vertex, we start by drawing its mark $M_l$ in an i.i.d.\ way with distribution
    \eqn{
    \label{dist-mark}
    \prob(M=m)=w_m/\ell_n, \quad
    1 \leq m \leq n.\nn
    }
When we have already explored a vertex with the same mark as the one drawn, we turn the status of the vertex to be explored to inactive, the potential vertex does not become a real vertex, and we proceed with the next potential vertex. When, instead, it receives a mark that we have not yet seen, then the potential vertex becomes a real vertex, its mark $M_l\in [n]$ indicating to which vertex in $[n]$ the $l$th explored vertex corresponds, so that $M_l\in \cluster(i)$. We then draw $X_l=\Poi(w_{M_l})$, and $X_l$ denotes the number of potential vertices incident to the real vertex $M_l$. Again, upon exploration, these potential vertices might become real vertices, and this occurs precisely when their mark corresponds to a vertex in $[n]$ that has not appeared in the cluster exploration so far. We call the above procedure of drawing a mark for a potential vertex to investigate whether it corresponds to a real vertex
a \emph{vertex check}. Let
    \eqn{
    \label{Ztn-def}
    \ZZ_t^{\sss(n)}(i)=n^{-1/(\tau-1)} Z_{\lceil tn^{\rho}\rceil}(i)\ \text{ for }\ t>0.\nn
    }
Then, by imitating the techniques used in the proof of \cite[Theorem 2.4]{SBVHVJL12}, we obtain
	\eqn{
	\label{Ztn-conv}
	 (\ZZ_t^{\sss(n)}(i))_{t>0}\convd  (\SS_t(i))_{t>0}.\nn
	}
(\cite[Theorem 2.4]{SBVHVJL12} states the result for $i=1$. However the exact same proof goes through for any $i\geq 2$.) The limiting process $(\SS_t(i))_{t>0}$ is defined as follows: Let
	\begin{align}\label{eqn:def-a-b}
	a=c_{\sss F}^{1/(\tau-1)}/\expec[W]\qquad \text{ and }\qquad b=b(i)=(c_{\sss F}/i)^{1/(\tau-1)}.
	\end{align}
We let $(\II_i(t))_{i\geq 1}$ denote independent increasing indicator processes defined by
    \begin{equation}\label{Iit-def}
    \II_i(s)=\ind\set{\Exp(a i^{-1/(\tau-1)})\in[0,s]},\qquad s\geq 0,
    \end{equation}
so that
	\begin{align*}
    	\pr\left(\II_i(s)=0 \ \forall s\in[0,t]\right)=\exp\left(-at/i^{1/(\tau-1)}\right).
	\end{align*}
Here $\big(\Exp(a i^{-1/(\tau-1)})\big)_{i\geq 1}$ are independent exponential random variables with rates $a i^{-1/(\tau-1)}$.
Then we define
	\begin{align}\label{eqn:SSt-def-abc}
    	\SS_t(i)=b-abt+ct+\sum_{j\neq i}^{\infty} \frac{b}{j^{1/(\tau-1)}}\left[\II_j(t)- \frac{at}{j^{1/(\tau-1)}}\right]
	\end{align}
for all $t\geq 0$, where $c=\lambda+\zeta-ab$ and $\zeta$ is as in \eqref{eqn:t-nr}.
We call $(\SS_t)_{t\geq 0}$ a \emph{thinned L\'evy process}.
\nomenclature{$(\SS_t)_{t\geq 0}$}{Thinned L\'evy process.}

Let $\HH_n^{\sss(i)}(u)$ denote the hitting time of $u$ of the process $(\ZZ_t^{\sss(n)}(i))_{t> 0}$. Then, by \cite[Corollary 3.4]{SBVHVJL12}, $\HH_n^{\sss(i)}(u)\convd \HH_{\SS(i)}(u),$ the hitting time of $u$ of the process $(\SS_t(i))_{t>0}$.
This implies the existence of a $B_{\vep,i}$ (independent of $n$) and an integer $n_{i,\eps}$ such that
	\begin{align}\label{eqn:hitting-time-lower-bound}
	\pr\left(\HH_n^{\sss(i)}\big((c_{\sss F}/2i)^{1/(\tau-1)}\big)\leq B_{\vep,i}\right)
	\leq \vep 2^{-i}\ \text{ for }\ n\geq n_{i,\eps},
	\end{align}
since the limiting process $(\SS_t(i))_{t>0}$ starts from $(c_{\sss F}/i)^{1/(\tau-1)}$ and takes a positive amount of time to reach $(c_{\sss F}/2i)^{1/(\tau-1)}$.

Let $|B(i,r)|$ denote the number of vertices in $B(i,r)$. Let $\delta_{\vep,i}$ be so small that
	\eqn{
	\label{delta-vep-i}
	(c_{\sss F}/2i)^{1/(\tau-1)}\delta_{\vep,i}<B_{\vep,i}.
	}
Then we claim that for all $\delta\in (0,\delta_{\vep,i}]$,
	\eqn{
	\label{|B|-H}
	\prob\Big(|B(i,\delta )|\leq (c_{\sss F}/2i)^{1/(\tau-1)} \delta n^{\rho}\Big)
	\leq \prob\Big(\HH_n^{\sss(i)}\big((c_{\sss F}/2i)^{1/(\tau-1)}\big)\leq B_{\vep,i}\Big).
	}
That \eqref{|B|-H} holds can be seen as follows. For
$|B(i,\delta)|\leq (c_{\sss F}/2i)^{1/(\tau-1)} \delta n^{\rho}$ to occur, there has to exist some $j\in[1, \delta n^{\eta}]$ such that the number of vertices at distance $j$ from $i$ is smaller than $(c_{\sss F}/2i)^{1/(\tau-1)} \delta n^{\rho}/(\delta n^{\eta})$, i.e,
\begin{align}\label{eqn:hit-before-time}
	\min_{1\leq j\leq \delta n^{\eta}}|\partial B(i, jn^{-\eta})|
	\leq (c_{\sss F}/2i)^{1/(\tau-1)}n^{1/(\tau-1)}.
\end{align}
Now the number of vertices at distance $j$ from $i$ is precisely the number of vertices in generation $j$ of the breadth-first exploration process, and hence this number (scaled by $n^{\rho}$) appears in the function $\ZZ_t^{\sss(n)}(i)$. Thus, \eqref{eqn:hit-before-time} implies that $(\ZZ_t^{\sss(n)}(i))_{t>0}$ has to hit $(c_{\sss F}/2i)^{1/(\tau-1)}$ before we have finished exploring up to generation $\delta n^{\eta}$, i.e., we must have that
	\eqn{
	\HH_n^{\sss(i)}\big((c_{\sss F}/2i)^{1/(\tau-1)}\big)
	\leq \frac{|B(i,\delta)|}{n^{\rho}}\leq \left(\frac{c_{\sss F}}{2i}\right)^{1/(\tau-1)}\delta<B_{\eps, i},\nn
	}
where the last inequality holds by \eqref{delta-vep-i} and because $\delta\in(0, \delta_{\eps,i}]$.

Combining \eqref{eqn:hitting-time-lower-bound} and \eqref{|B|-H}, we conclude that for all $\delta\in (0,\delta_{\vep,i}]$ and $n\geq n_{i,\eps}$,
	\eqn{
	\label{|B|-small-bd}
	\prob\Big(|B(i,\delta)|\leq (c_{\sss F}/2i)^{1/(\tau-1)} \delta n^{\rho}\Big)
	\leq \vep 2^{-i}.
	}
This explains the second term in \eqref{weight-balls-lb}.

To see what happens when $|B(i,\delta)|\geq (c_{\sss F}/2i)^{1/(\tau-1)} \delta n^{\rho}$, recall that
the vertices appear in a size-biased fashion in our exploration process. Hence
	\eqan{
	\label{bd-large-balls}
	&\pr\left(\sum_{j\in B(i,\delta)}w_j\leq \left(\frac{c_{\sss F}}{2i}\right)^{1/(\tau-1)}\frac{\delta n^{\rho}}{2},
	|B(i,\delta)|\geq\left(\frac{c_{\sss F}}{2i}\right)^{1/(\tau-1)} \delta n^{\rho}\right)\\
	&\hskip30pt\leq \pr\left(\sum_{j=1}^{\delta n^{\rho}\left(c_{\sss F}/(2i)\right)^{1/(\tau-1)}} w_{\pi_i(j)}\leq
	\left(\frac{c_{\sss F}}{2i}\right)^{1/(\tau-1)}\frac{\delta n^{\rho}}{2}\right)\nn\\
    	&\hskip30pt\leq n\exp\left(-\frac{J \delta n^{\rho}}{i^{1/(\tau-1)}}\right)\nn,
	}
by Lemma \ref{lem-sb-reor}. Combining \eqref{|B|-small-bd} and \eqref{bd-large-balls} proves the claim.
\qed
	
\begin{lemma}\label{lem-intrin-ball}
For $v\in [n]$, let $\cC(v)$ denote the component of $v$ in $\NRnwl$.
Then for every fixed $i\geq 1$ and $\vep_1,\vep_2>0$, there exist $\xi=\xi_{\vep_1,\vep_2}^{\sss (i)}>0$ and an integer ${\bar n}={{\bar n}_{\vep_1,\vep_2}}^{\sss (i)}$ such that
	\begin{align*}
	\prob\left(\min_{v\in \cC(i)}
	\Big(\sum_{j\in B(v, \vep_1)}w_j\Big)\leq \xi n^{\rho}\right)\leq \vep_2\ \text{ for }\ n\geq\bar n.
	\end{align*}
\end{lemma}

\noindent{\bf Proof:} Recall Proposition \ref{prop-diam-non-hubs}, and choose $N_{\vep_1,\vep_2}$ and $n_{\eps_1,\eps_2}$ large so that
	\eqn{\label{eqn:44}
	\pr\left(\diam\left(\NRnwl\setminus [\ch{N_{\vep_1,\vep_2}}]\right) \leq\vep_1  n^{\eta}/2\right)\geq 1-\vep_2
	}
for all $n\geq n_{\eps_1, \eps_2}$. Let
	\[
	F_1=\set{\diam(\NRnwl\setminus [\ch{N_{\vep_1,\vep_2}}]) \leq \eps_1 n^{\eta}/2}\qquad \text{ and }
	\qquad F_2=\set{\diam(\cC(i))>\eps_1 n^{\eta}/2}.
	\]
Clearly, on the set $F_1\cap F_2$,
	\begin{align}\label{eqn:45}
	\min_{v\in \cC(i)}\set{
	\sum_{j\in B(v, \vep_1 )}w_j}\geq \min_{k\in [N_{\vep_1,\vep_2}]}
	\set{\sum_{j\in B(k, \vep_1/2)}w_j}.
	\end{align}
Recall the definition of $\delta_{\vep,i}$ in \eqref{delta-vep-i}, and let
	\eqn{
	\label{eta-def}
	\Delta_{\vep_1,\vep_2}=\vep_1\wedge \left(\delta_{\vep_2,1}\wedge \cdots \wedge \delta_{\vep_2,N_{\vep_1,\vep_2}}\right)/2.\nn
	}
Then \eqref{eqn:45} implies
\begin{align*}
\min_{v\in \cluster(i)} \sum_{j\in B(v, \vep_1)}w_j
\geq \min_{k\in [N_{\vep_1,\vep_2}]}\sum_{j\in B(k, \Delta_{\vep_1,\vep_2})}w_j\nn
\end{align*}
on the set $F_1\cap F_2$. Hence, for all $n\geq n_{\vep_1,\vep_2}$,
\begin{align}\label{eqn:47}
	&\prob\left(F_1\cap F_2\cap\set{\min_{v\in \cluster(i)}
	\set{\sum_{j\in B(v, \vep_1)}w_j}\leq \left(\frac{c_{\sss F}}{2N_{\vep_1,\vep_2}}\right)^{1/(\tau-1)}\frac{\Delta_{\vep_1,\vep_2}n^{\rho}}{2}}\right)\\
	&\hskip15pt\qquad\leq \sum_{k=1}^{N_{\vep_1,\vep_2}}
	\prob\Big(
	\sum_{j\in B(k, \Delta_{\vep_1,\vep_2})}w_j\leq
	\Big(\frac{c_{\sss F}}{2N_{\vep_1,\vep_2}}\Big)^{1/(\tau-1)} \frac{\Delta_{\vep_1,\vep_2}}{2} n^{\rho}\Big)\nn\\
	&\hskip15pt\qquad\leq \sum_{k=1}^{N_{\vep_1,\vep_2}}
	\left(n\exp\left(-\frac{J \Delta_{\vep_1,\vep_2}
	n^{\rho}}{N_{\vep_1,\vep_2}^{1/(\tau-1)}}\right)+\frac{\vep_2}{2^{k}}\right)\nn\\
	&\hskip15pt\qquad\leq n^2\exp\left(-\frac{J \Delta_{\vep_1,\vep_2}n^{\rho}}{N_{\vep_1,\vep_2}^{1/(\tau-1)}}\right)+\vep_2,\nn
\end{align}
where the second inequality is a consequence of Lemma \ref{lem-weights}.

Next, on the set $F_1\cap F_2^c$,
\begin{align*}
\sum_{j\in B(v,\eps_1)}w_j=\sum_{j\in \cC(i)}w_j
\end{align*}
for any $v\in\cC(i)$. Further, by \cite[Theorem 1.4]{SBVHVJL12}, $n^{-\rho}\sum_{j\in \cC(i)}w_j$ converges in distribution to a positive random variable. Hence, there exists $\xi_{\eps_2}^{\sss(i)}>0$ such that
	\begin{align}\label{eqn:48}
	&\limsup_{n\to\infty}\ \prob\left(F_1\cap F_2^c\cap\set{\min_{v\in \cC(i)}
	\Big(\sum_{j\in B(v, \vep_1)}w_j\Big)\leq \xi_{\eps_2}^{\sss(i)} n^{\rho}}\right)\\
	&\hskip50pt\leq \limsup_{n\to\infty}\
	\pr\left(\sum_{j\in \cC(i)}w_j\leq \xi_{\eps_2}^{\sss(i)} n^{\rho}\right)\leq\eps_2.\nn
	\end{align}
The result follows upon combining \eqref{eqn:44}, \eqref{eqn:47} and \eqref{eqn:48}.
\qed
\medskip

We are now ready for the proof of Lemma \ref{lem:ghp-from-gw}:\\
\noindent
{\bf Proof of Lemma \ref{lem:ghp-from-gw}:}  Using Proposition \ref{prop:enough-to-work-with-high-weight-vertices}, for any $i\geq 1$ and $\eps>0$, we can choose $K$ such that
	\begin{align}\label{eqn:65}
	\pr\left(\cC_i(\lambda)\cap [K]=\emptyset\right)\leq\eps/2.
	\end{align}
By Lemma \ref{lem-intrin-ball}, we can choose $\xi>0$ and an integer $\bar n$ such that
	\begin{align}\label{eqn:66}
	\prob\left(\min_{v\in \cC(k)}
	\Big(\sum_{j\in B(v, \delta)}w_j\Big)\leq \xi n^{\rho}\right)\leq \vep/(2K)
	\end{align}
for all $n\geq\bar n$ and $k\in [K]$. Combining \eqref{eqn:65} and \eqref{eqn:66}, we see that
	\begin{align*}
	\pr\left(\left(\fm_i^{(n)}(\delta)\right)^{-1}>1/\xi\right)\leq \eps\ \text{ for }\ n\geq\bar n,
	\end{align*}
which yields the desired tightness.
\qed

\section{Proofs: Fractal dimension}
\label{sec:fractal-proofs}

In this section, we prove the assertions about the box-counting dimension. Throughout this section, $C,C'$ will denote universal constants whose values may change from line to line.

We first prove a similar result for the component  of $j$, $\cC(j)$. Consider $\cC(1)$, and as usual, view $\cC(1)$ as a metric measure space via the graph distance and by assigning mass $p_v:=w_v/(\sum_{\ell\in\cC(1)}w_{\ell})$ to vertex $v\in\cC(1)$. Set $\vp:=(p_v: v\in\cC(1))$. Now note that conditional on the vertex set of $\cC(1)$, $\cC(1)$ has the same distribution as the graph $\tilde{\cG}_m(\vp, a)$ where $a=(1+\lambda n^{-\eta})(\sum_{j\in\cC(1)} w_j)^2/\ell_n$. Using \cite[Proposition 3.7]{SBVHVJL12} and \cite[Lemma 3.1]{SBVHVJL12}, it is easy to verify that the conditions in Assumption \ref{ass:vp-a} hold with this choice of $a$ and $\vp$. Thus, by Theorem \ref{thm:connected-convg}, $n^{-\eta}\cC(1)$ converges in Gromov-weak topology to a limiting space that we denote by $\cM(1)$. Further, the sequence $\set{n^{-\eta}\cC(1)}_{n\geq 1}$ satisfies the global lower mass-bound property by Lemma \ref{lem-intrin-ball}. Hence,
	\begin{align}\label{eqn:ghp-convergence-C-1}
	n^{-\eta}\cC(1)\convd \cM(1)
	\end{align}
with respect to the Gromov-Hausdorff-Prokhorov topology. By similar arguments, we can show that $n^{-\eta}\cC(j)\convd \cM(j)$ with respect to the Gromov-Hausdorff-Prokhorov topology for any $j\geq 1$ and an appropriate (random) compact metric measure space $\cM(j)$. In Section \ref{sec-ub-dim}, we identify the upper box-counting dimension, and in Section \ref{sec-lb-dim} the lower box-counting dimension.

\subsection{Upper bound on the Minkowski dimension}
\label{sec-ub-dim}
The key ingredient in the proof is the following lemma:

\begin{prop}
\label{lem:dim-C-1}
Write $\pi=(\tau-2)/(\tau-3)$. Then for every $j\geq 1$,
	\begin{align*}
	\pr\left(\odim(\cM(j))>\pi\right)=0.
	\end{align*}

\end{prop}

\noindent{\bf Proof:} For simplicity, we work with $j=1$. The proof is similar for any $j\geq 2$. Recall that $\cN(\cM, \delta)$ denotes the minimum number of open balls of radius $\delta$ needed to cover the compact space $\cM$. Write
	\begin{align}
	\label{defs-bc}
	\fN_{\sss(\infty)}(\eps):=\cN(\cM(1), \eps)\qquad \text{ and }
	\qquad \fN_{\sss(n)}(\eps):=\cN(\cC(1), \eps n^{\eta}).
	\end{align}
Since the convergence in \eqref{eqn:ghp-convergence-C-1} holds with respect to the Gromov-Hausdorff topology, for every $x, \eps>0$,
	\begin{align}
	\label{eqn:61}
	\pr\left(\fN_{\sss(\infty)}(2\eps)>x\right)\leq\limsup_n\ \pr\left(\fN_{\sss(n)}(\eps)>x\right).
	\end{align}

Fix an arbitrary $\delta>0$ and, for any $\eps>0$, define
	\begin{align}\label{eqn:62}
	x_{\eps}:=\eps^{-\delta-\pi},\qquad \ u_{\eps}:=|\log\eps|,\qquad \
	\delta':=\frac{\delta}{2}\left(\frac{\tau-1}{\tau-2}\right),\qquad \text{ and }\qquad N(\eps)=\eps^{-\delta'-1/\eta}.
	\end{align}
Let $E_n$ be the event defined in \eqref{eqn:def-E-n}. Clearly, on the event $E_n\cap\set{\fN_{\sss(n)}(\eps)>x_{\eps}}$, any $v\in\cC(1)$ is within distance $\eps n^{\eta}$ from a point in $\cC(1)\cap [N(\eps)]$. Hence,
	\begin{align}\label{eqn:63}
	\pr\left(\fN_{\sss(n)}(\eps)>x_{\eps}\right)
	\leq \pr(E_n^c)+\pr\left(|\cC(1)\cap [N(\eps)]|\geq x_{\eps}\right),
	\end{align}
and, by Proposition \ref{prop-diam-non-hubs},
	\begin{align}\label{eqn:64}
	\limsup_n\ \pr(E_n^c)\leq c_{\delta'}\exp\left(-C/\eps^{\delta'\eta}\right).
	\end{align}
It remains to bound $\pr\left(|\cC(1)\cap [N(\eps)]|\geq x_{\eps}\right)$. To this end, note that by \cite[Proposition 3.7]{SBVHVJL12},
	\begin{align}\label{eqn:608}
	|\cC(1)\cap [N(\eps)]|\convd\sum_{q=1}^{N(\eps)}\cI_q\left(\cH_{\cS(1)}(0)\right),
	\end{align}
where $\cI_q(\cdot)$ and $\cH_{\cS(1)}(\cdot)$ are as defined around \eqref{eqn:SSt-def-abc}. Further, \cite[Theorem 1.4]{remcokliemjohan} implies the existence of positive constants $A_1$ and $A_2$ such that
	\begin{align}\label{eqn:609}
	\pr\left(\cH_{\cS(1)}(0)>u_{\eps}\right)\leq A_1\exp(-A_2 u_{\eps}^{\tau-1}).
	\end{align}
Combining \eqref{eqn:63},\eqref{eqn:64}, \eqref{eqn:608} and \eqref{eqn:609}, we conclude that, for any $u_{\eps}>0$,
	\begin{align}\label{eqn:67}
	\limsup_n\ \pr\left(\fN_{\sss(n)}(\eps)>x_{\eps}\right)
	\leq
	c_{\delta'}\exp\left(-C\eps^{-\delta'\eta}\right)
	+A_1\exp(-A_2 u_{\eps}^{\tau-1})
	+\pr\left(\sum_{q=1}^{N(\eps)}\cI_q\left(u_{\eps}\right)\geq x_{\eps}\right).
	\end{align}
Now $\cI_q\left(u_{\eps}\right)$ are i.i.d.\ Bernoulli random variables with
	\[
	\pr(\cI_q\left(u_{\eps}\right)=1)=1-\exp\left(-au_{\eps}/q^{1/(\tau-1)}\right)=:p_q,
	\]
where $a$ is as in \eqref{eqn:def-a-b}. Choose $s>0$ small so that $\e^s-1\leq 2s$. Clearly
	\begin{align*}
	\E\exp\left(s\cI_q\left(u_{\eps}\right)\right)&=1+p_q\left(\e^s-1\right)
	\leq\exp\left(p_q\left(\e^s-1\right)\right)\leq\exp(2sp_q).
	\end{align*}
Hence, there exists a constant $A_3>0$ such that
	\begin{align}\label{eqn:68}
	\pr\left(\sum_{q=1}^{N(\eps)}\cI_q\left(u_{\eps}\right)\geq x_{\eps}\right)
	&\leq \exp\left(-sx_{\eps}+2s\sum_{q=1}^{N(\eps)}p_q\right)\\
	&\leq \exp\left(-sx_{\eps}+2sA_3 u_{\eps} N(\eps)^\rho\right)\notag\\
	&=\exp\left(-sx_{\eps}+2sA_3u_{\eps}\eps^{-\frac{\delta}{2}-\pi}\right).\notag
	\end{align}
Combining \eqref{eqn:61}, \eqref{eqn:67} and \eqref{eqn:68}, we see that $\sum_{k=1}^{\infty}\pr\left(\fN_{\sss(\infty)}(2/k)>k^{\delta+\pi}\right)<\infty.$ Since $\delta>0$ was arbitrary, we conclude that
	\begin{align*}
	\limsup_k\ \frac{\log \left(\fN_{\sss(\infty)}(2/k)\right)}{\log(k/2)}\leq \pi\qquad a.s.
	\end{align*}
By sandwiching $\eps$ between $2/(k-1)$ and $2/k$, we get the desired upper bound on $\odim(\cM(1))$. \qed
\bigskip

\noindent{\bf Proof of upper bounds in \eqref{eqn:dim-nr} and \eqref{eqn:dim-mc}:} We only give the proof of \eqref{eqn:dim-nr}. This will imply \eqref{eqn:dim-mc} because of \eqref{eqn:def-M-nr}. Fix $i\geq 1$ and let
	\begin{align*}
	K_n:=\min\set{j\in[n]\colon \ j\in\cC_i(\lambda)}.
	\end{align*}
By Proposition \ref{prop:enough-to-work-with-high-weight-vertices}, $K_n$ is tight. By passing to a subsequence if necessary, we can assume that we are working on a space where
	\begin{align*}
	\left(n^{-\eta}\vM_n^{\nr}(\lambda), K_n\right)
	\to\left(\vM_{\infty}^{\nr}(\lambda), K_{\infty}\right)\qquad a.s.
	\end{align*}
for some (integer-valued) random variable $K_{\infty}$. Then
	\begin{align*}
	\pr\left(\odim\left(M_i^{\nr}(\lambda)\right)>\pi\right)
	=\sum_{j=1}^{\infty}\pr\left(\odim\left(M_i^{\nr}(\lambda)\right)>\pi,\ K_{\infty}=j\right).
\end{align*}
By Proposition \ref{lem:dim-C-1}, $\pr\left(\odim\left(M_i^{\nr}(\lambda)\right)>\pi,\ K_{\infty}=j\right)=0$
for every $j\geq 1$, and hence
\begin{align}\label{eqn:fff-1}
	\odim\left(M_i^{\nr}(\lambda)\right)\leq\pi\ \text{ a.s.}
\end{align}
This completes the proof of the upper bound on the Minkowski dimension. \qed
\bigskip

\subsection{Lower bound on the Minkowski dimension}
\label{sec-lb-dim}

We next extend the argument for the upper bound to prove a lower bound on the Minkowski dimension of $\cM(j)$.
As in \eqref{eqn:61},
	\begin{align}
	\label{eqn:61b}
	\pr\left(\fN_{\sss(\infty)}(\eps/2)<x\right)\leq\limsup_n\ \pr\left(\fN_{\sss(n)}(\eps)<x\right).
	\end{align}
Recall the definitions in \eqref{defs-bc}, and for an arbitrary $\delta>0$ and $\eps>0$, adapt \eqref{eqn:62} to
	\begin{align}\label{eqn:62b}
	\underline{x}_{\eps}:=\eps^{\delta-\pi},\qquad
	\delta':=\frac{\delta}{\pi}\left(1-h\right),
	\qquad
	\text{ and }
	\qquad
	\underline{N}(\eps)=\eps^{-(1-\delta')/\eta},
	\end{align}
where $\pi=(\tau-2)/(\tau-3)$ as in Proposition \ref{lem:dim-C-1}, and $h>0$ is sufficiently small so that
	\eqan{
	\label{eqn:power-eps-positive}
	\kappa_3:=2-\delta-(1-\delta')\left(\frac{3\tau-8}{\tau-3}\right)+\frac{\tau-2}{\tau-3}>0,\ \text{ and }\
	\kappa_4:=1-\delta-(1-\delta')\left(\frac{2\tau-5}{\tau-3}\right)+\frac{\tau-2}{\tau-3}>0.
	}
(A simple calculation will show that it is possible to choose $h>0$ small so that \eqref{eqn:power-eps-positive} holds whenever $\tau>3$.)

The main result in this section is the following estimate on $\fN_{\sss(n)}(\eps):=\cN(\cC(j), \eps n^{\eta})$:

\begin{prop}
\label{prop:dim-C-2}
There exist $\kappa>0$ and $c>0$ such that
	\eqn{
	\label{prob-lb-Md}
	\limsup_n\ \pr\big(\fN_{\sss(n)}(\eps)<\underline{x}_{\eps}\big)\leq c\eps^{\kappa}.
	}
Consequently, for every $j\geq 1$,
	\begin{align*}
	\pr\left(\udim(\cM(j))<\pi\right)=0.
	\end{align*}
\end{prop}
\bigskip

The rest of this section is devoted to the proof of Proposition \ref{prop:dim-C-2}. As in Section \ref{sec-ub-dim} and for simplicity, we work with $j=1$. The proof is similar for any $j\geq 2$.
Before starting with the proof, we collect some preliminaries. The proof below relies on two asymptotic bounds on $|\cC(1)|$. For this, we use
	\eqn{
	\label{eqn:401}
	\limsup_n\ \pr\big(n^{-\rho}|\cC(1)|\leq s)=\pr(\cH_{\cS(1)}(0)\leq s),
	}
where $\cH_{\cS(1)}(\cdot)$is defined around \eqref{eqn:SSt-def-abc}. Our main result on the lower tails of the distribution of $\cH_{\cS(1)}(0)$ is in the following lemma:

\begin{lemma}[Lower tails of $\cH_{\cS(1)}(0)$]
\label{lem-LT-H}
There exists $C>0$ such that
	\eqn{
	\label{lower-tail}
	\pr(\cH_{\cS(1)}(0)\leq s)\leq C s.
	}
\end{lemma}
\noindent{\bf Proof:} We note that
	\eqn{
	\pr\big(\cH_{\cS(1)}(0)\leq s\big)=\pr\big(\exists t\leq s\colon \SS_t(1)=0\big).
	}
We split
	\eqn{
    \label{eqn:343}
	\SS_t(1)=b-abt+ct+\frac{b}{a}\big(\RR_t-\DD_t\big),
	}
where, abbreviating $d_j=a/j^{1/(\tau-1)}$,
	\eqan{
	\RR_t =
	\sum_{j\geq 2}^{\infty} d_j\left[N_j(t)- d_jt\right],\ \text{ and }\
	\DD_t =\sum_{j\geq 2}^{\infty} d_j[N_j(t)-1]\indic{N_j(t)\geq 2}.
	}
Here $(N_j(t))_{t\geq 0}$ are independent rate $d_j$ Poisson processes. Thus, $(\RR_t)_{t\geq 0}$ is a L\'evy process, while $(\DD_t)_{t\geq 0}$ subtracts the multiple hits. When $b>0$ and $t\leq s$ with $s$ small,  and using that $\DD_s$ is non-decreasing,
	\eqn{
	\label{split-SS}
	\pr\big(\cH_{\cS(1)}(0)\leq s\big)\leq \pr\bigg(\inf_{t\in [0,s]} \RR_t\leq -a/4\bigg)+\pr\big(\DD_s\geq a/4\big).
	}
We start with the latter contribution. Since, for a Poisson random variable $Z$ with parameter $\lambda$,
	\begin{align*}
	\expec\big[(Z-1)\indic{Z\ge 2}\big]=\sum_{k\geq 2}(k-1) \e^{-\lambda}\frac{\lambda^k}{k!}
	=\lambda^2 \sum_{k\geq 2}\frac{1}{k}\e^{-\lambda}\frac{\lambda^{k-2}}{(k-2)!}
	\leq \frac{\lambda^2}{2},
	\end{align*}
we have
	\begin{align*}
	\pr\big(\DD_s\geq a/4\big)\leq \frac{4}{a}\expec[\DD_s]\leq \frac{2}{a}\sum_{j\geq 2}d_j (d_js)^2
	=\frac{2s^2}{a}\sum_{j\geq 2}d_j^3.
   \end{align*}
For the first term in \eqref{split-SS}, we use Doob's $L^2$-inequality to bound
	\begin{align*}
	\pr\bigg(\inf_{t\in [0,s]} \RR_t\leq -a/4\bigg)\leq \frac{16}{a^2} \expec[\RR_s^2]
	=\frac{16}{a^2} \sum_{j\geq 2}^{\infty} d_j^2\Var(N_j(s))
	=\frac{16}{a^2} \sum_{j\geq 2}^{\infty} d_j^2(d_js)=\frac{16 s}{a^2} \sum_{j\geq 2}^{\infty} d_j^3,
	\end{align*}
so that \eqref{lower-tail} follows.
\qed
\medskip

\begin{lemma}[Cluster weight convergence]
\label{lem-conv-weights-clusters}
For a set of vertices $A\subseteq [n]$, let $w(A)=\sum_{a\in A} w_a$ denote its weight. Then, as $n\rightarrow \infty$, for every $j\geq 1$, $\expec[n^{-\rho}w(\cluster(j))]$ remains uniformly bounded as $n\rightarrow \infty$, where $\rho$ is as in \eqref{eqn:def-eta-rho},.
\end{lemma}

\noindent{\bf Proof:} Fix $K\geq 0$ so large that
	\eqn{
	\label{nu-K-bd}
	\nu_n^{\sss(K+1)}:=\frac{1}{\ell_n}\sum_{i\in [n]\setminus [K]} w_i^2\leq
	\frac{1-n^{-\eta}}{1+|\lambda| n^{-\eta}}.
	}
This is possible, since $\ell_n/n\rightarrow \expec[W]$, while
	\eqan{
	\label{nu-K-bd-1}
	\frac{1}{\ell_n}\sum_{i\in [n]\setminus [K]} w_i^2&=\nu_n-\frac{1+o(1)}{\expec[W]n}\sum_{i\leq K} w_i^2
	\leq \nu_n-Cn^{-1+2/(\tau-1)} \sum_{i\leq K} i^{-2/(\tau-1)}\\
	&\leq 1-C'n^{-\eta} K^{\eta},\nn
	}
where we have used \eqref{eqn:98} in the second step to lower bound $\sum_{i\leq K} w_i^2$. We write $\cluster(A)=\bigcup_{a\in A} \cluster(a)$. Then, for $j\leq K$, we bound
	\eqn{
	w(\cluster(j))\leq w(\cluster([K])).
	}
We next investigate $\expec[w(\cluster([K]))]$. Note that
	\eqn{
	\label{pij-bd}
	p_{ij}=1-\e^{-(1+\lambda n^{-\eta}) w_iw_j/\ell_n}\leq (1+\lambda n^{-\eta})
	\frac{w_iw_j}{\ell_n},
	}
where $\ell_n=\sum_{l\in[n]}w_l$ is the total weight. Thus, for any $A\subseteq [n]$,
	\eqn{
	\label{dist-Aj-bd-1}
	\prob\big(\dist(A,j)=l\big)\leq \sum_{a\in A}\sum_{i_1,\ldots, i_{l-1}\in [n]}^* \prod_{s=1}^l p_{i_{s-1},i_s},
	}
where $i_0=a, i_l=j$ and the sum is over distinct vertices not in $A$. Using the bound on $p_{i,j}$ and performing the sum over $i_1,\ldots, i_{l-1}$, we obtain that
	\eqan{
	\label{dist-Aj-bd-2}
	\prob(\dist([K],j)=l)&\leq (1+\lambda n^{-\eta})\sum_{a\in [K]}\frac{w_a w_j}{\ell_n} \bigg((1+\lambda n^{-\eta})\nu_n^{\sss (K+1)}\bigg)^{l-1}\\
	& =(1+\lambda n^{-\eta})w([K])\frac{w_j}{\ell_n} \bigg[(1+\lambda n^{-\eta})\nu_n^{\sss (K+1)}\bigg]^{l-1}.
	}
By \eqref{nu-K-bd},
	\eqn{
	(1+\lambda n^{-\eta})\nu_n^{\sss (K+1)}\leq 1-n^{-\eta}.
	}
As a result, for large $n$,
	\eqan{
	\expec\big[w(\cluster([K]))\big]&\leq 2w([K])\Big[1+\sum_{j\in [n]\setminus [K]} \frac{w_j^2}{\ell_n}\sum_{l\geq 1}
	(1-n^{-\eta})^{l-1}\Big]\\
	&\leq 4w([K])\Big[1+\sum_{l\geq 0}
	(1-n^{-\eta})^{l}\Big]=8w([K]) n^{\eta}.\nn
	}
Since, by an argument similar to \eqref{nu-K-bd-1},
	\eqn{
	w([K])\leq CK^{\rho} n^{1/(\tau-1)},
	}
we arrive at
	\eqn{
	\expec\big[w(\cluster([K]))\big]\leq CK^{\rho} n^{\eta+1/(\tau-1)}=CK^{\rho} n^{\rho}.
	}
This completes the proof.
\qed

\medskip

By \eqref{eqn:401} and Lemma \ref{lem-LT-H},
	\eqn{
	\limsup_n\ \pr\big(|\cC(1)|\leq \eps^{\delta h/2}n^\rho\big)\leq C\eps^{\delta h/2}.
	}
We conclude that
	\eqan{
	\label{split-Min-dim}
	\limsup_n\ \pr\big(\fN_{\sss(n)}(\eps)<\underline{x}_{\eps}\big)
	&\leq
	\limsup_n\ \pr\Big(\{\fN_{\sss(n)}(\eps)<\underline{x}_{\eps}\}\cap\{|\cC(1)|> \eps^{\delta h/2}n^{\rho}\}\Big)
	+C\eps^{\delta h/2}.
	}

We now study the event in \eqref{split-Min-dim}.
We note that $\fN_{\sss(n)}(\eps)\geq X^{\sss(n)}(\vep)$, which is defined as
	\eqn{
	\label{lb-box-counting}
	X^{\sss(n)}(\vep)=1+\sum_{i=2}^{\underline{N}(\vep)} \indic{i\in \cluster(1)}
	\indic{\dist_{\sss \cluster(1)}(i, [i-1])>4\vep n^{\eta}},
	}
where $\dist_{\sss \cluster(1)}(A,B)$ is the graph distance between the sets of vertices $A\cap\cC(1)$ and $B\cap\cC(1)$.
Indeed, we start counting in the order $i\geq 1$, and determine whether an extra ball is needed to cover vertex $i$ after we have covered the vertices in $[i-1]\cap\cC(1)$. The first contribution in \eqref{lb-box-counting} comes from the ball that covers vertex $1$.

Use inclusion-exclusion to write $X^{\sss(n)}(\vep)$ as
	\eqn{
	X^{\sss(n)}(\vep)=X_1^{\sss(n)}(\vep)-X_2^{\sss(n)}(\vep),
	}
where
	\eqn{
	X_1^{\sss(n)}(\vep)=\sum_{i=1}^{\underline{N}(\vep)} \indic{i\in \cluster(1)},
	\qquad
	X_2^{\sss(n)}(\vep)=\sum_{i=2}^{\underline{N}(\vep)}\indic{i\in \cluster(1)}
	\indic{\dist_{\sss \cluster(1)}(i, [i-1])\leq 4\vep n^{\eta}}.
	}
Therefore,
	\begin{equation}\label{eqn:444}
	\pr\Big(\{\fN_{\sss(n)}(\eps)<\underline{x}_{\eps}\}
	\cap\{|\cC(1)|> \eps^{\delta h/2}n^{\rho}\}\Big)
	\leq
	\pr\Big(\{X_1^{\sss(n)}(\vep)<2\underline{x}_{\eps}\}
	\cap\{|\cC(1)|>\eps^{\delta h/2}n^{\rho}\}\Big)
	+\pr\big(X_2^{\sss(n)}(\vep)>\underline{x}_{\eps}\big).
	\end{equation}
We will show that the limsup as $n\rightarrow \infty$ of the first probability is bounded by $C\vep^{\kappa_1}$, and the limsup as $n\rightarrow \infty$ of the second by $C\vep^{\kappa_2}$ with $\kappa_1,\kappa_2>0$, so that  Proposition \ref{prop:dim-C-2} will follow with $\kappa=\min\{\delta h/2, \kappa_1, \kappa_2\}$.
\medskip

\paragraph{{\bf Analysis of $X_1^{\sss(n)}$.}}
It follows from \cite[Proposition 3.7]{SBVHVJL12} that
	\eqan{\label{eqn:0}
	\limsup_{n\rightarrow \infty}\ \pr\Big(\{X_1^{\sss(n)}(\vep)<2\underline{x}_{\eps}\}
	\cap\{|\cC(1)|>\eps^{\delta h/2}n^{\rho}\}\Big)
	\leq\pr\big(\underline{X}_1(\vep)\leq 2\underline{x}_{\eps}\big),
	}
where
\begin{align*}
\underline{X}_1(\vep):=\sum_{i=1}^{\underline{N}(\vep)} \mathcal{I}_i(\eps^{\delta h/2})
\end{align*}
is a sum of independent indicators with success probabilities $1-\exp\big(-d_i \eps^{\delta h/2}\big)$, $i=1,\hdots, \underline{N}(\vep)$ (recall \eqref{Iit-def}), with $d_j$ as defined right below \eqref{eqn:343}. Note that
    \begin{equation}\label{eqn:1}
	\expec\big[\underline{X}_1(\vep)\big]
    =\sum_{i=2}^{\underline{N}(\vep)}\prob\big(\mathcal{I}_i(\eps^{\delta h/2})=1\big)
	=\sum_{i=2}^{\underline{N}(\vep)}\bigg[1-\exp\big(-d_i \eps^{\delta h/2}\big)\bigg]
    \leq\sum_{i=2}^{\underline{N}(\vep)}d_i \eps^{\delta h/2}
	\leq C\eps^{\delta h/2}\underline{N}(\vep)^{\rho}.
	\end{equation}
Similarly, for small enough $\eps>0$,
   \eqan{\label{eqn:2}
	\expec\big[\underline{X}_1(\vep)\big]
	=\sum_{i=2}^{\underline{N}(\vep)}\bigg[1-\exp\big(-d_i \eps^{\delta h/2}\big)\bigg]
	\geq \tfrac{1}{2}\sum_{i=2}^{\underline{N}(\vep)}d_i \eps^{\delta h/2}
	\geq C'\eps^{\delta h/2}\underline{N}(\vep)^{\rho}\geq 3\underline{x}_{\eps}.
	}
Further, since $\underline{X}_1(\vep)$ is a sum of independent indicators,
	\eqan{\label{eqn:3}
	\Var(\underline{X}_1(\vep))\leq \expec[\underline{X}_1(\vep)].
	}
Combining \eqref{eqn:0}, \eqref{eqn:1}, \eqref{eqn:2}, and \eqref{eqn:3}, we get
	\eqan{
	\label{eqn:445}
	&\limsup_{n\rightarrow \infty}\ \prob\big(\{X_1^{\sss(n)}(\vep)\leq 2\underline{x}_{\eps}\}
	\cap \{|\cC(1)|> \eps^{\delta h/2}n^{\rho}\}\big)\\
	&\hskip30pt\leq \prob\big(\underline{X}_1(\vep)\leq 2\underline{x}_{\eps}\big)
      \leq \prob\Big(\big|\underline{X}_1(\vep)-\expec[\underline{X}_1(\vep)]\big|
	\geq \underline{x}_{\eps}\Big)
	\leq \underline{x}_{\eps}^{-2}\Var(\underline{X}_1(\vep))\nn\\
	&\hskip60pt\leq \underline{x}_{\eps}^{-2}\expec[\underline{X}_1(\vep)]
     \leq C\eps^{\kappa_1},\nn
	}
where $\kappa_1=2\pi-2\delta+\delta h/2-\pi(1-\delta')>0$ when $\delta>0$ is sufficiently small. This proves a bound on the first term on the right side of \eqref{eqn:444}.
\medskip

\paragraph{{\bf Analysis of $X_2^{\sss(n)}$.}}
We next give an upper bound on $\pr(X_2^{\sss(n)}(\vep)\geq\underline{x}_{\eps})$. We start with
	\eqn{
	\label{Markov-X2}
	\prob\big(X_2^{\sss(n)}(\vep)
	\geq \underline{x}_{\eps}\big)
	\leq
	\underline{x}_{\eps}^{-1}
	\expec\big[X_2^{\sss(n)}(\vep)\big].
	}
Further,
	\eqan{
    \label{eqn:495}
	\expec\big[X_2^{\sss(n)}(\vep)\big]
	&=\sum_{i=2}^{\underline{N}(\vep)}\prob\Big(i\in \cluster(1), \dist_{\sss \cluster(1)}(i, [i-1])\leq 4\vep n^{\eta}\Big).
	}
When $i\in \cluster(1)$ and $\dist_{\sss \cluster(1)}(i, [i-1])\leq 4\vep n^{\eta}$, there has to be $j\in [i-1]$ and $k\in [n]$ such that the three events
	\begin{enumerate}
	\item[(i)] $\{\dist(i, k)\leq 4\vep n^{\eta}\}$;
	\item[(ii)] $\{\dist(j, k)\leq 4\vep n^{\eta}\}$;
	\item[(iii)]  $\{k\in \cluster(1)\}$,
	\end{enumerate}
occur disjointly, where $\dist(i,j)$ denotes the graph distance in the random graph $\NRnw$. There are two cases depending on whether $k>\underline{N}(\vep)$ or $k\leq \underline{N}(\vep)$. When $k\leq \underline{N}(\vep)$, we can ignore the event $\{\dist(j, k)\leq 2\vep n^{\eta}\}$. This gives, for $2\leq i\leq \underline{N}(\vep)$,
	\eqan{
	&\prob\Big(i\in \cluster(1), \dist_{\sss \cluster(1)}(i, [i-1])\leq 4\vep n^{\eta}\Big)\\
	&\qquad \leq\sum_{j=1}^{\underline{N}(\vep)} \sum_{k>\underline{N}(\vep)}
	\prob\Big(\{\dist(i, k)\leq 4\vep n^{\eta}\}\circ
	\{\dist(j, k)\leq 4\vep n^{\eta}\}\circ \{k\in \cluster(1)\}\Big)\nn\\
	&\qquad\qquad+\sum_{k=1}^{\underline{N}(\vep)}
	\prob\Big(\{\dist(i, k)\leq 4\vep n^{\eta}\}\circ
	\{k\in \cluster(1)\}\Big),\nn
	}
where, for two increasing events $A,B$, we write $A\circ B$ for the event that $A$ and $B$ occur disjointly.

By the BK inequality, we bound
	\eqan{
	\prob\Big(i\in \cluster(1), \dist_{\sss \cluster(1)}(i, [i-1])\leq 4\vep n^{\eta}\Big)
	&\leq \sum_{j=1}^{\underline{N}(\vep)} \sum_{k>\underline{N}(\vep)}
	\prob\big(\dist(i, k)\leq 4\vep n^{\eta}\big)
	\prob\big(\dist(j, k)\leq 4\vep n^{\eta}\big)
	\prob\big(k\in \cluster(1)\big)\nn\\
	&\qquad+\sum_{k=1}^{\underline{N}(\vep)}
	\prob\big(\dist(i, k)\leq 4\vep n^{\eta}\big)
	\prob\big(k\in \cluster(1)\big).
	}
Similar to \eqref{dist-Aj-bd-2}, we have
	\eqn{
	\prob\big(\dist(i,j)=l\big)\leq \bigg(1+\frac{\lambda}{n^{\eta}}\bigg)\frac{w_iw_j}{\ell_n} \nu_n(\lambda)^{l-1},
	}
where $\nu_n(\lambda)=(1+\lambda n^{-\eta})\nu_n$.
In our case, $\nu_n=1+O(n^{-\eta})$, so that, for $l\leq 4\vep n^{\eta}$,
	\eqn{
	\prob\big(\dist(i,j)\leq l\big)\leq Cl\frac{w_i w_j}{\ell_n}.
	}
Further,
	\eqan{
	\prob\big(k\in \cluster(1)\big)&=\indic{k=1} +\sum_{l\in [n]}
	\prob\big(\{l\in \cluster(1)\}\circ \{kl\text{ occupied}\}\big)\leq \indic{k=1} +
	\sum_{l\in [n]}\bigg(1+\frac{\lambda}{n^{\eta}}\bigg) \frac{w_k w_l}{\ell_n}
	\prob\big(l\in \cluster(1)\big)\nn\\
	&= \indic{k=1} +\bigg(1+\frac{\lambda}{n^{\eta}}\bigg) \frac{w_k}{\ell_n} \expec\big[w(\cluster(1))\big],\nn
	}
where we recall that $w(A)=\sum_{a\in A}w_a$ denotes the total weight of $A$. By Lemma \ref{lem-conv-weights-clusters}, $\expec[n^{-\rho}w(\cluster(1))]$ remains uniformly bounded as $n\rightarrow \infty$. We conclude that
	\eqan{
	&\prob\Big(i\in \cluster(1), \dist_{\sss \cluster(1)}(i, [i-1])\leq 4\vep n^{\eta}\Big)\\
	&\qquad \leq C\vep^2 n^{2\eta+\rho} \sum_{j=1}^{\underline{N}(\vep)} \sum_{k>\underline{N}(\vep)} \frac{w_i w_j w_k^3}{\ell_n^3}
	+C\vep n^{\eta+\rho} \sum_{k=2}^{\underline{N}(\vep)}
	\frac{w_i w_k^2}{\ell_n^2}
	+C\vep n^{\eta}
	\frac{w_i w_1}{\ell_n}\nn\\
	&\qquad\leq C'\vep^2 n^{2\eta+\rho-3+5/(\tau-1)}
	\sum_{j=1}^{\underline{N}(\vep)} \sum_{k>\underline{N}(\vep)} i^{-1/(\tau-1)} j^{-1/(\tau-1)} k^{-3/(\tau-1)}\nn\\
	&\qquad\qquad +C'\vep n^{\eta+\rho-2+3/(\tau-1)}
	\sum_{k=2}^{\underline{N}(\vep)}i^{-1/(\tau-1)} k^{-2/(\tau-1)}+C'\vep n^{\eta-1+2/(\tau-1)}i^{-1/(\tau-1)},\nn
	}
where the last step uses the first inequality in \eqref{eqn:98}. Note that
\[2\eta+\rho-3+5/(\tau-1)=\eta+\rho-2+3/(\tau-1)=\eta-1+2/(\tau-1)=0,\]
so that the powers of $n$ cancel. Combining the above with \eqref{eqn:495} leads to
	\eqan{
	\expec\big[X_2^{\sss(n)}(\vep)\big]
	&\leq C\vep^2 \Big(\sum_{j=1}^{\underline{N}(\vep)} j^{-1/(\tau-1)}\Big)^2
	\sum_{k>\underline{N}(\vep)} k^{-3/(\tau-1)}
	+C\vep\sum_{j=1}^{\underline{N}(\vep)} j^{-1/(\tau-1)}
	\sum_{k=1}^{\underline{N}(\vep)} k^{-2/(\tau-1)}.
	}
Note that
	\eqan{
	\sum_{j=1}^{N} j^{-p/(\tau-1)}&=O(N^{(\tau-1-p)/(\tau-1)})\ \text{ for }\ p=1,2,\ \text{ and }\
	\sum_{k>N} k^{-3/(\tau-1)}
	=O(N^{(\tau-4)/(\tau-1)}).
	}
Thus
	\eqan{
	\expec\big[X_2^{\sss(n)}(\vep)\big]
	&\leq C\vep^2
	\underline{N}(\vep)^{2(\tau-2)/(\tau-1)+(\tau-4)/(\tau-1)}
	+C\vep\underline{N}(\vep)^{(\tau-2)/(\tau-1)+(\tau-3)/(\tau-1)}\nn\\
	&=C\Big[\vep^2 \underline{N}(\vep)^{(3\tau-8)/(\tau-1)}
	+\vep \underline{N}(\vep)^{(2\tau-5)/(\tau-1)}\Big].
	}
Using \eqref{Markov-X2} and plugging in the values $\eta=(\tau-3)/(\tau-1), \pi=(\tau-2)/(\tau-3)$, we arrive at
	\eqan{
	\label{Markov-X2-rep}
	\limsup_{n\rightarrow \infty}\ \prob\big(X_2^{\sss(n)}(\vep)\geq \underline{x}_{\eps}\big)
	&\leq C
	\underline{x}_{\eps}^{-1}\Big[\vep^2 \underline{N}(\vep)^{(3\tau-8)/(\tau-1)}
	+\vep \underline{N}(\vep)^{(2\tau-5)/(\tau-1)}\Big]=C[\vep^{\kappa_3}+
	\vep^{\kappa_4}\big],
	}
where the exponents $\kappa_3$ and $\kappa_4$ are positive because of the choice of $\delta'$ (see \eqref{eqn:power-eps-positive}).

\medskip

\paragraph{{\bf Completion of the proof of Proposition \ref{prop:dim-C-2}:}}
Note that
\eqref{prob-lb-Md} follows upon combining \eqref{split-Min-dim}, \eqref{eqn:444}, \eqref{eqn:445}, and \eqref{Markov-X2-rep}.
Now fix $p>1/\kappa$, where $\kappa$ is as in \eqref{prob-lb-Md}. Then $\sum_{k=1}^{\infty}\pr\left(\fN_{\sss(\infty)}(1/k^p)<(2k)^{(\pi-\delta)p}\right)<\infty.$ Since $\delta>0$ was arbitrary, we conclude that
	\begin{align*}
	\liminf_k\ \frac{\log \left(\fN_{\sss(\infty)}(1/k^p)\right)}{\log(k^p)}
	\geq \pi\quad a.s.
	\end{align*}
By sandwiching $\eps$ between $1/(k-1)^p$ and $1/k^p$, we obtain the bound: $\udim(\cM(1))\geq\pi$ a.s. \qed
\medskip

\noindent{\bf Proof of \eqref{eqn:dim-nr} and \eqref{eqn:dim-mc}:} Proposition \ref{prop:dim-C-2} combined with an argument identical to the the one given right after the proof of Proposition \ref{lem:dim-C-1} yields the lower bound: $\udim\left(M_i^{\nr}(\lambda)\right)\geq\pi$ a.s. \eqref{eqn:dim-nr} follows once we combine this lower bound with \eqref{eqn:fff-1}, and \eqref{eqn:dim-mc} follows as a consequence of \eqref{eqn:def-M-nr}.
\qed

\section{Open problems}
\label{sec:disc-open}

In Theorem \ref{thm:mc-main-1}, we have considered a general entrance boundary $\vc\in l_0$. To study specific properties of the limit objects, we focused mainly on the special case $\vc=\vc(\alpha,\tau)$ as in \eqref{eqn:c-tau-lamb-def} and in this case, we have shown compactness and identified the box counting dimension in Theorem \ref{thm:mc-main-2-tau}. An important problem in this context is to establish necessary and sufficient conditions on $\vc$ that ensure compactness of the limiting spaces.

Another motivation for pursuing this problem comes from the following simple corollary of Theorem \ref{thm:mc-main-2-tau}: For any $i\geq 1$, consider the sequence $\mvtheta^{\sss(i)}$ as in \eqref{eqn:barg-mvtheta-def-comps}. Then $\cT_{(\infty)}^{\mvtheta^{\sss (i)}}$ is almost surely compact.
Similarly, compactness of $\cM(1)$ (as defined in \eqref{eqn:ghp-convergence-C-1}) implies compactness of the associated ICRT $\cT_{(\infty)}^{\overline\mvtheta}$ where $\overline{\mvtheta}=(\overline{\theta}_i\colon i\geq 1)$ is given by the following prescription: Let $q_k$ be such that
	\begin{align*}
	\sum_{q=1}^{q_k}\cI_q\left(\cH_{\cS(1)}(0)\right)=k,
	\end{align*}
where $\cI_q(\cdot)$ and $\cH_{\cS(1)}(\cdot)$ are as defined around \eqref{eqn:SSt-def-abc}. Define
	\begin{align*}
	\overline{\theta}_i= \frac{q_i^{-1/(\tau-1)}}{\left(\sum_{k=1}^{\infty}q_k^{-2/(\tau-1)}\right)^{1/2}}\
	\qquad
	\text{ for }\
	\qquad i\geq 1.
	\end{align*}

These can be thought of as ``annealed results," since $\mvtheta^{\sss(i)}$ and $\overline\mvtheta$ are random. No result is known in this direction without a prior distribution on $\mvtheta$, i.e., sufficient conditions on non-random $\mvtheta\in\Theta$ that ensure compactness of the tree $\icrt$ are not known. In \cite[Section 7]{AMP}, Aldous, Miermont and Pitman conjecture that boundedness of $\icrt$ for $\mvtheta\in\Theta$ is equivalent to $\int_1^{\infty}(\psi_{\mvtheta}(u))^{-1}du<\infty$, where $\psi_{\mvtheta}$, in our situation, is given by
\begin{align}\label{eqn:def-psi}
\psi_{\mvtheta}(u)=\sum_{i=1}^{\infty}\left(\exp(-u\theta_i)-1+u\theta_i\right).
\end{align}
This conjecture, however, is open to date. Our proof technique demonstrates a method of proving such annealed results via approximation by random graphs. Thus, classification of those $\vc\in l_0$ for which the spaces $M_i^{\vc}(\lambda)$ are compact will lead to a broad class of prior distributions on $\mvtheta$ for which $\icrt$ is compact.
\begin{problem}\label{prob:compactness}
Find necessary and sufficient conditions on $\vc$ that ensure compactness of the spaces $M_i^{\vc}(\lambda)$ for $i\geq 1$.
\end{problem}

Another related problem is to find the fractal dimensions of the limiting spaces. As a corollary to Theorem \ref{thm:mc-main-2-tau}, we get
\begin{align}\label{eqn:odim}
\dim\left(\cT_{(\infty)}^{\mvtheta^{\sss (i)}}\right)= (\tau-2)/(\tau-3)\ a.s.
\end{align}
where $\mvtheta^{\sss (i)}$ is as in \eqref{eqn:barg-mvtheta-def-comps} corresponding to $\vc$ of the form \eqref{eqn:c-tau-lamb-def}. Proposition \ref{lem:dim-C-1} and Proposition \ref{prop:dim-C-2} show that the assertion in \eqref{eqn:odim} remains true if we replace $\mvtheta^{\sss (i)}$ by $\overline\mvtheta$. Now, it is not hard to prove that
\begin{align*}
\inf_j\ \overline{\theta}_j j^{1/(\tau-2)}>0\quad a.s.\quad\text{and}\quad\sup_j\ \overline{\theta}_j j^{1/(\tau-2)}<\infty\quad a.s.
\end{align*}
It then follows that
\begin{align*}
\tau-2=\sup\set{a\geq 0\ :\ \lim_{u\to\infty}u^{-a}\psi_{\overline\mvtheta}(u)=\infty}
=\inf\set{a\geq 0\ :\ \lim_{u\to\infty}u^{-a}\psi_{\overline\mvtheta}(u)=0}\ a.s.,
\end{align*}
which in turn implies that both the Hausdorff dimension and the packing dimension of a $\psi_{\overline\mvtheta}$ L\'evy tree equal $(\tau-2)/(\tau-3)$ a.s. (see \cite{duquesne2005probabilistic, haas-miermont}). Using the analogy between ICRTs and L\'evy trees as in \cite[Section 7]{AMP}, it is natural to expect that the same is true for $\cT_{(\infty)}^{\overline\mvtheta}$ and hence for $\cM(1)$. This is the heuristic behind Conjecture \ref{conj:fractal-dimension}.

\begin{problem}\label{prob:fractal}
Prove Conjecture \ref{conj:fractal-dimension}.
\end{problem}

\vskip0.5cm

\section*{Acknowledgements}
The authors are indebted to Gr\'egory Miermont for many enlightening discussions about inhomogeneous continuum random trees. SS thanks ENS Lyon for hospitality and accommodation during visits.
The authors thank Igor Kortchemski for drawing their attention to his recent preprint \cite{Kort15}.
The authors also thank an anonymous referee who pointed out a number of issues in an earlier version which significantly improved the readability of the manuscript.
SB has been partially supported by NSF-DMS grants 1105581, 1310002,  160683, 161307 and SES grant 1357622.
RvdH and SS have been supported in part by the Netherlands Organisation
for Scientific Research (NWO) through the Gravitation Networks
grant 024.002.003. In addition, RvdH has been supported by VICI grant 639.033.806 and SS has been supported by a CRM-ISM fellowship.

\bibliographystyle{plain}
\bibliography{critical-3-4}

\small
\printnomenclature[3.1cm]
\normalsize

\end{document}